\documentclass[12pt]{article}
\usepackage{amssymb}
\input epsf
\usepackage{epic}
\renewcommand{\arraystretch}{1.3}

\catcode`\@=11
\def\marginnote#1{}

\newcount\hour
\newcount\minute
\newtoks\amorpm
\hour=\time\divide\hour by60 \minute=\time{\multiply\hour by60
\global\advance\minute by-\hour}
\edef\standardtime{{\ifnum\hour<12 \global\amorpm={am}%
        \else\global\amorpm={pm}\advance\hour by-12 \fi
        \ifnum\hour=0 \hour=12 \fi
        \number\hour:\ifnum\minute<10 0\fi\number\minute\the\amorpm}}
\edef\militarytime{\number\hour:\ifnum\minute<10
0\fi\number\minute}

\def\draftlabel#1{{\@bsphack\if@filesw {\let\thepage\relax
      \xdef\@gtempa{\write\@auxout{\string
          \newlabel{#1}{{\@currentlabel}{\thepage}}}}}\@gtempa \if@nobreak
    \ifvmode\nobreak\fi\fi\fi\@esphack} \gdef\@eqnlabel{#1}}
    \def\@eqnlabel{}
\def\@vacuum{}
\def\draftmarginnote#1{\marginpar{\raggedright\scriptsize\tt#1}}

\def\draft{
  \oddsidemargin -.5truein
  \def\@oddfoot{\footnotesize \sl preliminary draft \hfil
    \rm\thepage\hfil\sl\today\quad\militarytime}
  \let\@evenfoot\@oddfoot \overfullrule 3pt
    \let\label=\draftlabel
    \let\marginnote=\draftmarginnote
  \def\@eqnnum{(\theequation)\rlap{\kern\marginparsep\tt\@eqnlabel}%
    \global\let\@eqnlabel\@vacuum}

  }

\newtheorem{theorem}{Theorem}
\newtheorem{proposition}[theorem]{Proposition}
\newtheorem{lemma}[theorem]{Lemma}
\newtheorem{cor}[theorem]{Corollary}

\newtheorem{defin}{Definition}
\newtheorem{remark}{Remark}

\newtheorem{constr}{Construction}

\let\Bbb=\mathbb
\let\wtd=\widetilde
\def\XX{{\hbox{\scriptsize${{\hbox{\tiny$\times$}}\atop{
\hbox{\tiny$\times$}}}$}}}
\newcommand{\ORD}[1]{\XX{#1}\XX}

\renewcommand{\theequation}{\thesection.\arabic{equation}}

\newcommand{\newsection}{
\setcounter{equation}{0}
\setcounter{theorem}{0}
\setcounter{con}{0}
\setcounter{defin}{0}
\setcounter{remark}{0}
\setcounter{example}{0}
\section}
\def\appendix#1{
\addtocounter{section}{1} \setcounter{equation}{0}
\renewcommand{\thesection}{\Alph{section}}
\section*{Appendix \thesection\protect\indent
#1}
}
\newcommand{\cpict}[3]{
\dimen1=#1\advance\dimen1 by-\hsize\divide\dimen1 by-2
\vtop to #2{
\noindent\hskip\dimen1{\special{em:graph #3.bmp}}
\vfil}\hskip-2cm
}

\newcommand{\tr}{\,{\rm tr}\,}
\def\e{{\,e}\,}

\def\eps{\varepsilon}

\def\be{\begin{equation}}
\def\ee{\end{equation}}
\def\bea{\begin{eqnarray}}
\def\eea{\end{eqnarray}}

\def\RR{{\Bbb R}}

\begin{document}
\setlength{\unitlength}{1.5mm}

\begin{center}
{\Large On Quantizing Teichm\"uller and Thurston theories}\\
\vspace{4mm}
\end{center}
\begin{center}
\vspace{4mm}
{L. Chekhov${}^{{\rm a)}}$\footnote{E-mail: chekhov@mi.ras.ru.}
 and R. C. Penner${}^{{\rm b)}}$%
\footnote{E-mail: rpenner@math.usc.edu.} }\\
\vspace{18pt}
${}^{{\rm a)}}${\it Steklov Mathematical Institute,}
\\ {\it  Gubkina 8, 117966, GSP--1, Moscow, Russia} \\[.2cm]
${}^{{\rm b)}}${\it Departments of Mathematics and
Physics/Astronomy}\\{\it University of Southern California,}
\\ {\it Los Angeles, CA 90089 USA}
\end{center}
\vskip 0.9 cm
\begin{abstract}
\noindent {\bf Abstract} In earlier work, Chekhov and Fock have given a quantization of Teichm\"uller space as a Poisson
manifold, and the current paper first surveys this material adding further mathematical
and other detail, including the underlying geometric work by Penner on classical Teichm\"uller theory.
In particular, the earlier quantum ordering solution is found to essentially agree with
an ``improved'' operator ordering given by serially traversing general edge-paths on a graph in the underlying
surface. Now, insofar as Thurston's sphere of projectivized foliations of compact support provides a useful
compactification for Teichm\"uller space in the classical case, it is natural to consider corresponding limits
of appropriate operators to provide a framework for
studying degenerations of quantum hyperbolic structures.  After surveying the required background
material on Thurston theory and ``train tracks'', the current paper continues to give a quantization of
Thurston's boundary in the special case of the once-punctured torus, where there are already
substantial analytical and combinatorial challenges.  Indeed, an operatorial version of continued fractions
as well as the improved quantum ordering are required to prove existence of these limits.
Since Thurston's boundary for the once-punctured torus is a topological circle, the main new result
may be regarded as a quantization of this circle.
There is a discussion of quantizing Thurston's boundary spheres for higher genus surfaces
in closing remarks.
\end{abstract}

\renewcommand{\thefootnote}{\arabic{footnote}}
\setcounter{footnote}{0}

\newsection{Introduction}
One manifestation of Teichm\"uller space in contemporary
mathematical physics
is as the Hilbert space and
the algebra of observables for three-dimensional (3D) quantum gravity since
E.~Verlinde and H.~Verlinde~\cite{VV} have argued that
the classical phase space
of Einstein gravity in a 3D manifold is the Teichm\"uller space of its
boundary.  (Analogously, the classical phase space for 3D
Chern--Simons theory is the moduli space of flat connections on the
boundary; this theory was quantized in~\cite{FockRosly,FR2}.)
Teichm\" uller space possesses its canonical (Weil--Petersson)
Poisson structure, whose symmetry group is the mapping class group
of orientation-preserving homeomorphisms modulo isotopy.  The algebra of
observables is the collection of geodesic length functions of geodesic representatives
of homotopy classes of essential closed curves together with its natural mapping class group action.

Given this Poisson structure, one can turn to the problem of quantizing it,
thereby obtaining a variant of the quantum 3D gravity description.
According to the correspondence principle: (1) the algebra of observables of
the corresponding quantum theory is the noncommutative deformation of the
$*$-algebra of functions on it governed by the Poisson
structure; (2) the Hilbert spaces of the theory are the representation
spaces of these $*$-algebras; and (3) the symmetry group acts on the algebra of
observables by automorphisms.  Under the assumption that the quantization of a
Poisson manifold exists and is unique, to solve this problem it
suffices to construct a family of $*$-algebras, which depend
on the quantization
parameter $\hbar$, and an action of the mapping class group on this family
by outer automorphisms, and to show that the algebra and the action
thus constructed
reproduces the
classical algebra, the classical action, and the classical Poisson structure in the
limit $\hbar \rightarrow 0$. This program has been successfully performed in
\cite{ChF} and \cite{kashaev}, and we describe it in Secs.~2 and~3 of this
paper. Actually, the problem that was solved in these papers differs
slightly from the original formulation because the methods of \cite{ChF}
are suited for describing only {\it open\/} surfaces (surfaces with nonempty
boundary, components of which can, however, reduce to a puncture).
The corresponding Teichm\"uller space has a degenerate
Weil--Petersson Poisson structure, while the mapping class group acts as symmetry
group. We describe the deformation quantization of the corresponding
Teichm\"uller space, the action of the mapping class group by
outer automorphisms,
the representations of the observable algebra, and the induced action of the
mapping class group on the representation space following
\cite{ChF}.

As in
\cite{VV}, the representation space of the observable algebra can also be
interpreted as the space of conformal blocks of the Liouville conformal field
theory. This program is under current development (see~\cite{Teich}).
Our construction can be therefore interpreted
as the construction of the conformal block spaces and the mapping
class group actions for this CFT.

The key point of the quantization procedure is constructing quantum mapping class
group transformations that define in a consistent way the morphisms between quantum
$*$-algebras simultaneously preserving the quantum geodesic algebra. The main
mathematical ingredient of the construction is a version of the quantum dilogarithm
by L.~D.~Faddeev~\cite{Faddeev}. We interpret the corresponding five-term relation
as the only nontrivial relation in a certain groupoid that has the mapping class
group as maximal subgroup.  A similar construction has been made independently and
simultaneously by R.~M.~Kashaev~\cite{kashaev}. The key difference between these two
constructions lies in the dimensions of the Poisson leaves of the two theories:
given a graph with $v$ three-valent vertices, $e$ edges, and $f=s$ faces suitably embedded in a
surface $F_g^s$ of genus $g$ with $s$ ideal boundary components, the genus is given
by Euler's formula $v-e+f=2-2g$.  In Kashaev's approach, there are $2e$ variables
and $v+f$ central elements (at each vertex and at each face), so the Poisson leaf
dimension is $2e-v-f=8g-8+3s$, while in our approach, there are $e$ variables and
$f$ central elements, so the Poisson dimension $e-f=6g-6+2s$ exactly coincides with
the dimension of the Teichm\"uller space of Riemann surfaces of genus~$g$ with $s$
punctures. The approach of \cite{ChF} is thus appropriate for describing 2D
topological theories while that of \cite{kashaev} is suitable for describing Liouville
field theory as a lattice theory (for instance, the Liouville field central charge
can be calculated, see \cite{Kashaev2}).

One of the mathematical tools employed in the quantization \cite{ChF} is the
decorated Teichm\"uller theory \cite{Penn1}, and the relevant aspects are briefly reviewed
in Section~2. In effect in the classical case, to each edge of a trivalent
``fatgraph'' $\Gamma$ (i.e., a graph plus a cyclic ordering of the half-edges about each vertex)
embedded as
a ``spine'' (i.e., a deformation retract of the surface)
$F=F_g^s$ is assigned a number $Z_\alpha$, where $\alpha$ here and below indexes the edges of $\Gamma$.  The
tuple $(Z_\alpha )$ gives global coordinates on an appropriate Teichm\"uller space
${\cal T}_H={\cal T}_H(F)$ of $F$ (as first studied by Thurston and later by Fock),
where any or all of the ``punctures'' of $F_g^s$ are permitted to be uniformized
instead as circular boundary components (and the subscript $H$ stands for
``holes''); the details are given in Section2.1.  The Weil--Petersson K\"ahler
two-form, which is known \cite{Penn2} in the $(Z_\alpha )$ coordinates, pulls back to a
degenerate two-form, i.e., to a degenerate Poisson structure on ${\cal T}_H$.
Furthermore, the action of the mapping class group $MC=MC(F)$ of $F$ on ${\cal
T}_H$, which is again known \cite{Penn2} in the $(Z_\alpha )$ coordinates as well as
combinatorially \cite{Penn3}, leaves invariant this Poisson structure, and this is the
Poisson manifold ${\cal T}_H$ with $MC$-action which has been quantized.

More explicitly still, for each edge $\alpha$, we may associate a pair of M\"obius
transformations $R_{Z_\alpha},L_{Z_\alpha}$ depending upon $Z_\alpha$ with the
following property. For any homotopy class $\gamma$ of geodesic in $F$ with
corresponding closed edge-path $P$ on $\Gamma$, consider the serial product
$P_\gamma$ of operators $R_{Z_\alpha},L_{Z_\alpha}$ taken in order as one traverses
$P$ in some orientation from some starting point, where one inserts the former (or latter,
respectively) if immediately after traversing edge $\alpha$, then $P$ turns right (or left) in
$\Gamma$; thus, the combinatorial geometry of serially traversing edges of $\Gamma$
as dictated by $P$ determines an ordered product $P_{\gamma}$ of matrices depending
upon $(Z_\alpha )$, and the length $l_\gamma$ of the geodesic representative of
$\gamma$ for the point of ${\cal T}_H$ determined by $(Z_\alpha )$ is given by
$G_\gamma =2~{\rm cosh}~(l_\gamma /2)=~|{\rm tr}~P_{\gamma}|$.  It is the
Poisson algebra of these geodesic functions $G_\gamma$, the algebra of
observables, which has been quantized.

In the quantum case (after passing to a suitable subspace on which the Poisson
structure is non-degenerate), standard techniques of deformation quantization
produce the appropriate Hilbert space ${\cal H}$ of the quantum theory, as well as
pairs of operators $R_{Z_\alpha},L_{Z_\alpha}$ on ${\cal H}^2$ for each edge
$\alpha$ of $\Gamma$.  Again, to the homotopy class of a geodesic $\gamma$ in $F$ or
its corresponding closed edge-path $P$ on $\Gamma$, we may assign the ordered
product $P_\gamma$ of these operators as dictated by the combinatorial geometry of
$P$ in $\Gamma$, whose trace $G_\gamma=~{\rm tr}~P_\gamma$ is the ``quantum geodesic
operator.'' The main point is to prove the invariance under the action of $MC$,
which is intimately connected with functional properties of the quantum dilogarithm
as was mentioned before.

In \cite{ChF} was proved the existence
and uniqueness of an appropriate ``proper quantum ordering'' of  operators that
enjoy $MC$-invariance as well as satisfy the standard physical requirements.
In Section~3.5, we observe
that the natural operatorial ordering given by the
combinatorial geometry of edge-paths in $\Gamma$ can be used
to derive this physically correct quantum ordering.  The improved
ordering is required in the subsequent quantization (discussed below).

In order to explain the further new results in this paper, we must recall aspects of
Thurston's seminal work on surface geometry, topology, and dynamics from the
1970-1980's, which is surveyed in Section 4.  Very briefly, Thurston introduced the
space ${\cal PF}_0={\cal PF}_0(F)$ of ``projective measured foliations of compact
support in $F$'' as a boundary for the Teichm\"uller space ${\cal T}={\cal T}(F)$,
where ``Thurston's compactification'' $\overline{\cal T}={\cal T}\cup {\cal PF}_0$
is a closed ball with boundary sphere ${\cal PF}_0$, where the action of $MC$ on
${\cal T}$ extends continuously to the natural action on $\overline{\cal T}$.
Furthermore, the sphere ${\cal PF}_0$ contains the set of all homotopy classes of
geodesics in $F$ as a dense subset, i.e., ${\cal PF}_0$ is an appropriate completion
of this set. (Unfortunately, the action of $MC$ on ${\cal PF}_0$ has dense orbits,
so the beautiful structure of Thurston's compactification does not descend in any
reasonable way, with the current state of understanding, to a useful structure on
the level of Riemann's moduli space.)
Thurston also devised an elegant graphical formalism for understanding ${\cal PF}_0$
using ``train tracks'', which are graphs embedded in $F$ with the further structure of a ``branched
one-submanifold" (cf. Section 4.2).  In effect, a maximal train track gives a chart
on the sphere ${\cal PF}_0$, and furthermore, the combinatorial expression for inclusion of charts
effectively captures the dynamics of the action of diffeomorphisms of $F$ (cf. Section 4.4).

Thus, Thurston theory arises as a natural
tool to understand degenerations in Teichm\"uller space and dynamics on $F$, and in light of remarks
above, its quantization should provide a natural tool for studying degenerations of
3D gravity or Liouville conformal field theory.  Many aspects of the survey of
classical Thurston surface theory (in Section 4) are
required for our subsequent quantization.

In any case, a mathematically natural problem armed with \cite{ChF} is to ``quantize
Thurston theory'': assign operators $H_\gamma$ on ${\cal H}$ to each homotopy class
$\gamma$ of geodesic in such a way that as $\gamma$ converges in ${\cal PF}_0$ to a
projectivized measured foliation $[{\cal F},\bar\mu ]$, then the corresponding
operators $H_\gamma$ converge in an appropriate sense to a well-defined operator
$H_{[{\cal F},\bar \mu ]}$ on ${\cal H}$.  Upon choosing a spine $\Gamma\subseteq F$
(i.e., an embedding, up to homotopy, of a graph whose inclusion is a homotopy equivalence),
any homotopy class of curve $\gamma$ may be essentially uniquely realized as an edge-path on $\Gamma$.
Define the ``graph length'' ${\rm g.l.} (\gamma )$ to be the total number of edges of $\Gamma$
traversed by this edge-path counted with multiplicities.

In the special case of the torus $F_1^1$,
we have succeeded here (in Section 5) in showing that the ratio of operators
$$
H_\gamma=~{\rm p.l.}(\gamma )/{\rm g.l.}(\gamma )
$$
converge
weakly to a well-defined operator as $\gamma$ converges in ${\cal PF}_0(F_1^1)$,
where the ``proper length'' is defined by
$${\rm p.l.}(\gamma )=\lim _{n\to\infty}{1\over n}~{\rm tr}~{\rm log}~2T_n({1\over 2} P_\gamma),$$
with $T_n$ the Chebyshev polynomials.  In fact in the classical case, ${\rm p.l.}(\gamma )$ agrees with
half the length of the geodesic homotopic to $\gamma$ in the Poincar\'e metric, and in the quantum case,
there is an appropriate operatorial interpretation, both of which are described in Section~5.1.
In particular, for several spines whose corresponding charts
cover the circle ${\cal PF}_0(F_1^1)$, the analysis involves rather intricate estimates.
This leads to a
natural operatorial quantization (in Section~5.3.4) of the standard simple continued fractions, which are
intimately connected with Thurston theory on $F_1^1$ (as we shall describe in
Section 4.5).  To complete the basic theory on the torus, one would like an intrinsic operatorial
description of the circle of unbounded operators we have constructed, as well as an
an explication of the mapping class group action on it, viz., Thurston's
classification of surface automorphisms.

This paper is organized as follows:  Section 2 covers classical Teichm\"uller theory and Section 3 the
Chekhov--Fock quantization; Section 4 surveys classical Thurston theory of surfaces, and Section 5 gives our
quantization of Thurston's boundary for the punctured torus.  Excerpts of Sections 2 and 3 are derived from an
earlier manuscript of Chekhov--Fock, and we strive to include further mathematical detail.  Section 4 surveys
aspects of a large literature on Thurston theory and train tracks, explicitly covering only what is required in
Section 5.  Section 5 should be regarded as work in physics in the sense that some of the formal calculations
depend upon manipulations of asymptotic spectral expansions for which there may be remaining mathematical
issues.  Closing remarks in Section 6 discuss the natural extension of these results to more complicated surfaces
as well as other related work. Appendix~A includes a novel proof of the required convergence in the classical
case for any surface (in a sense, a new proof of the existence of Thurston's compactification), which may yet be
useful in the quantum case.  Appendix~B contains an analysis of the Casimir operators in the Poisson algebra in
the $Z$-variables and the appropriate diagonalization of Poisson structure.

\vskip .2in

\noindent{\bf Acknowledgements}~We are indebted to L.~D.~Faddeev, V.~V.~Fock, M. Lapidus, and F.~Bonahon for
useful discussions. The work was partially supported by the RFFI Grant No.~01-01-00549 (L.Ch.), by the Program
Mathematical Methods in Nonlinear Dynamics (L.Ch.), and by the COBASE Project.

\vskip .2in

\newsection{Classical Teichm\"uller spaces}

To begin, we shall briefly recall the two related roles played by graphs in
Teichm\"uller theory as both aspects will be required here.

A {\it fatgraph} or {\it ribbon graph} is a graph $\Gamma$ together with a cyclic ordering on the half-edges
incident on each vertex, and we canonically associate to $\Gamma$ a surface $F(\Gamma )$ with boundary obtained
by ``fattening each edge of the graph into a band'' in the natural way; we shall tacitly require all vertices to
have valence at least three unless stated otherwise, and we shall call a fatgraph {\it cubic} if each vertex has
valence three.  To each homotopy class of homotopy equivalence $\iota :\Gamma \to F$ for some surface $F$, where
$\iota$ respects the orientation, there is a corresponding cell in ${\cal T}_g^s\times{\Bbb R}_{>0}^s$ as
explained in \cite{Penn4} in the hyperbolic setting \cite{Penn1} and in the conformal setting \cite{Streb1}.
Thus, a homotopy class of $\iota :\Gamma \to F$ is the name of a cell in the canonical cell decomposition.  We
shall sometimes suppress the mapping $\iota :\Gamma\to F$ and refer to $\Gamma$ itself as a {\it spine} of $F$,
where $\Gamma$ is identified with $\iota (\Gamma )\subseteq F$.  The cell decomposition is invariant under the
action of the mapping class group (induced by post-composition of $\iota$ with homeomorphisms), and this has
been an effective tool for studying Riemann's moduli space; for instance, we shall recall here the corresponding
presentation of the mapping class groups.

The second role of fatgraphs is exclusive to the hyperbolic setting, namely, fatgraphs provide a kind of
``basis'' for geometrically natural global parameterizations of Teichm\"uller space.  Specifically, fix a
homotopy class $\iota :\Gamma \to F$ as above, where we now demand that $\Gamma $ is cubic, and let $E=E(\Gamma
)$ denote the set of edges of $\Gamma $.  In several different contexts, one can naturally identify ${\Bbb
R}_{>0}^E$ with a suitable modification of an appropriate Teichm\"uller space of $F$.

For instance, for punctured surfaces, recall \cite{Penn1} that the lambda length of a pair of horocycles is
$\sqrt{2~e^\delta}$, where $\delta$ is the signed hyperbolic distance between the horocycles. Lambda lengths
give a global real-analytic parametrization of the decorated Teichm\"uller space as the trivial bundle
$\widetilde{\cal T}_g^s ={\cal T}_g^s\times{\Bbb R}_{>0}^s$ over Teichm\"uller space, where the fiber over a
point is the space of all $s$-tuples of horocycles in the surface, one horocycle about each puncture
(parameterized by hyperbolic length).

For another example, Thurston's shear
coordinates \cite{ThSh},\cite{Bon2} give global parameters not only on Teichm\"uller space
(cf. Section 2.1.1) but also on the related space of measured foliations (cf. Section \ref{scmf}).

On the level of Teichm\"uller space, the two global coordinate systems (lambda
lengths and shear coordinates) are closely related, and we choose to give the
exposition here principally in shear coordinates with the parallel lambda length
discussion relegated to a series of ongoing remarks.  On the other hand, certain
proofs of identities involving shear coordinates are easy calculations in lambda
lengths.

\subsection{Graph description of Teichm\"uller spaces}

\subsubsection{Global coordinates on Teichm\"uller space}

In addition to the Teichm\"uller space ${\cal T}_g^s$ and decorated Teichm\"uller
space $\widetilde{\cal T}_g^s$, we shall also require the following modification.
Given an open Riemann surface $F$ of finite topological type, a neighborhood of an
ideal boundary component is either an annulus or a punctured disk; in the former
case, the ideal boundary component will be called a ``true'' boundary component and
in the latter will be called a ``puncture.''  We shall study the latter as a
degeneration of the former with an elaboration
$$
{\cal T}\hbox{eich}(F)={\rm Hom}'\biggl (\pi _1(F), PSL_2({\Bbb R})\biggr )/PSL_2({\Bbb R})
$$
of the usual Teichm\"uller space, where ${\rm Hom}'$ denotes the space of all discrete faithful representations
with no elliptic elements, i.e., $|{\rm tr}~(\rho(\gamma))|\geq 2$ for all $\gamma\in\pi _1(F)$ for any
representation $\rho:\pi _1(F)\to PSL_2({\Bbb R})$.

Assume that $\iota :\Gamma \to F$ is a homotopy class of homotopy equivalence
and $\gamma\in \pi _1(F)$ is conjugate in $\pi _1(F)$ to
the boundary of a regular neighborhood of an ideal boundary component of
$F(\Gamma )$.
Thus, the ideal boundary component is a puncture if and only if $\rho (\gamma )$ is
a parabolic transformation, i.e., $|{\rm tr}~(\rho (\gamma ))|=2$.  For any
$\gamma\in\pi _1(F)$ with $|{\rm tr}~(\rho (\gamma ))|>2$, the underlying free
homotopy class of unbased curves contains a unique hyperbolic geodesic whose length
$l_\gamma $ is given by $G_\gamma =2{\rm cosh}(l _\gamma )=|{\rm tr}~\rho (\gamma
)|$, where $G_\gamma$ is called the {\it geodesic operator} and is constant on the
conjugacy class of $\gamma\in\pi _1(F)$. Furthermore, by definition,
\bea
\nonumber
{\cal T}_g^s=\{ [\rho ]\in{\cal T}\hbox{eich}(F): \forall\gamma\in\pi _1(F)~{\rm
freely~homotopic~into~the~boundary},\\\nonumber ~{\rm we~have}~ |{\rm tr}~(\rho
(\gamma ))|=2\}\hskip 2.7in
\\
\nonumber\subseteq {\cal T}\hbox{eich} (F)\hskip 4.0in\\\nonumber
\eea
Finally, define the space
${\cal T}_H(F)$ to be the $2^s$-fold cover of
${\cal T}\hbox{eich} (F)$ branched over ${\cal T}_g^s$, where the fiber is given
by the set of all orientations on the boundary components of $F$.

\begin{theorem}
Fix any spine $\Gamma \subseteq F$, where $\Gamma$ is a cubic fatgraph.  Then there is a real-analytic
homeomorphsim ${\Bbb R}^{E(\Gamma )}\to {\cal T}_H(F)$.  The hyperbolic length $l_\gamma$ of a true boundary
component $\gamma$ is given by $l_\gamma =|\sum Z_i|$, where the sum is over the set of all edges traversed by
$\gamma$ counted with multiplicity.  Furthermore, $\sum Z_i=0$ if and only if the corresponding ideal boundary
component is a puncture, so ${\cal T}_g^s\subseteq {\cal T}_H(F)$ is determined by $s$ independent linear
constraints.
\end{theorem}

The theorem is due to Thurston with a systematic study by Fock.  We shall not give a proof here (though there is
not a complete proof in the literature so far as we know), but we shall at least give the construction that
defines the homeomorphism ${\Bbb R}^{E(\Gamma )}\to {\cal T}_H(F)$.

The basic idea is to associate to each edge of $\Gamma $ an appropriate cross ratio.
To set this up, consider the topological surface $F^+\subseteq F$ obtained by
adjoining a punctured disk to each true boundary component of $F$.  The fatgraph
$\iota (\Gamma )\subseteq F\subseteq F^+$ is thus also a spine of $F^+$, and its
Poincar\'e dual in $F^+$ is an ideal triangulation $\Delta$ of $F^+$ (i.e., a
decomposition into triangles with vertices among the punctures).  In the universal
cover of $F^+$, each arc of $\Delta$ thus separates two complementary triangles
which combine to give a topological quadrilateral, and the basic idea is to
associate to each edge the cross-ratio of this quadrilateral.

To make this precise and describe the homeomorphism in the theorem, let $\alpha
=1,\ldots ,E=E(\Gamma )$ index the edges of $\Gamma $, and let $(Z_\alpha)$ denote a
point of  ${\Bbb R}^{E}$.
We associate the M\"obius transformation
\be
\label{XZ}
X_{Z_\alpha}=\left(
\begin{array}{cc} 0 & -\e^{Z_\alpha/2}\\
                \e^{-Z_\alpha/2} & 0\end{array}\right).
\ee
to the edge $\alpha$.  To explicate this definition, consider an ideal
quadrilateral in the hyperbolic plane triangulated by a diagonal into two
ideal triangles $T_1,T_2$.  We may conjugate in $PSL_2({\Bbb R})$ to arrange that the
vertices of $T_1$ are $0,-1,\infty$ and the vertices of $T_2$ are $0,-1,t$, where
$0<t<\infty$, and an appropriate cross ratio of the original quadrilateral is $t$.
Setting $Z_\alpha={\rm log}~t$ in the formula above, $X_{Z_\alpha}$ is the M\"obius
transformation interchanging $0,\infty$ and sending $-1$ to $t$, i.e., sending $T_1$
to $T_2$.  Notice that $X_{Z_\alpha}^2$ is the identity in $PSL_2({\Bbb R})$, so
$X_{Z_\alpha}$ also sends $T_2$ to $T_1$.

We also introduce the ``right'' and ``left'' turn matrices
\be
\label{R}
R=\left(\begin{array}{cc} 1 & 1\\ -1 & 0\end{array}\right), \qquad
L= R^2=\left(\begin{array}{cc} 0 & 1\\ -1 &
-1\end{array}\right),
\ee
and define the corresponding operators $R_Z$ and $L_Z$,
\bea
\label{Rz}
R_Z\equiv RX_Z&=&\left(\begin{array}{cc}
                \e^{-Z/2}&-\e^{Z/2}\\
                     0   &\e^{Z/2}
                     \end{array}\right),\\
\label{Lz}
L_Z\equiv LX_Z&=&\left(\begin{array}{cc}
                \e^{-Z/2}&   0\\
                 -\e^{-Z/2}&\e^{Z/2}
                     \end{array}\right).
\eea

Consider a closed oriented edge-path $P$ in $\Gamma $, where we assume that $P$
never consecutively traverses an oriented edge followed by its reverse, i.e., there
is no ``turning back''.  Choosing also an initial base point on $P$, we may imagine
the corresponding curve serially traversing the oriented edges of $\Gamma $ with
coordinates $Z_1,\ldots ,Z_n$ turning left or right from
$Z_{i}$ to
$Z_{i+1}$, for
$i=1,\ldots ,n$ (with the indices mod $n$ so that $Z_{n+1}=Z_1$). Assign to
$P$ the corresponding composition
\be
P_{Z_1\dots Z_n}=
L_{Z_n}L_{Z_{n-1}}R_{Z_{n-2}}\dots R_{Z_2}L_{Z_1},
\label{Pz}
\ee
where the
matrices $L_{Z_i}$ or $R_{Z_i}$ are inserted depending on
which turn---left or right---the path takes at the corresponding stage.

Fixing any base point, the assignment
$P\mapsto P_{Z_1,\ldots ,P_n}\in PSL_2({\Bbb R})$
gives rise to a representation $\rho\in{\cal T}_H(F)$, and this defines the
required map ${\Bbb R}^{E(\Gamma )}\to {\cal T}_H(F)$.  Furthermore summarizing
standard formulas and facts mentioned above, we have

\begin{proposition}\label{prop1}~
There is a one-to-one correspondence between the set of
conjugacy classes of elements of $\pi _1(F)$ and free homotopy classes
of closed oriented geodesics in $F$.  For any spine of $F$, each
free homotopy class is uniquely
represented by a cyclically defined closed edge-path $P$ with no turning back,
and the length of $\gamma$ is determined by
\be
\label{geod}
G_\gamma\equiv 2\cosh (l_\gamma/2)=|\tr P_{Z_1\dots Z_n}|.
\ee
\end{proposition}

By construction if $P$ corresponds to a boundary component $\gamma$ of $F(\Gamma   )$,
then the associated matrix
has the form $R_{Z_1}R_{Z_2}\dots R_{Z_n}$, or
$L_{Z_1}L_{Z_2}\dots L_{Z_n}$ depending on the orientation.
In this case, because all of the matrices $R_x$
($L_x$) are upper (lower) triangular, formula (\ref{geod}) gives
\be
\label{loophole}
l_{\gamma}=\left|\sum_{i=1}^{n}Z_i\right|,
\ee
where the sign of this sum gives the orientation of the boundary component, and
$l_\gamma =0$ corresponds to a puncture.  This proves the assertions about boundary
lengths.

The $Z$-coordinates (i.e., log cross ratios) are called {\it (Thurston) shear
coordinates} \cite{ThSh},\cite{Bon2} and can alternatively be defined by dropping perpendiculars from
each of the two opposite vertices to the diagonal $\alpha$ of a quadrilateral, and
measuring the signed hyperbolic distance $Z_\alpha$ along $\alpha$ between these two
projections.

Assume that there is an enumeration of the edges of $\Gamma $ and that
edge $\alpha$ has distinct endpoints. Given a spine $\Gamma$ of $F$, we may produce
another spine $\Gamma _\alpha$ of $F$ by contracting and expanding edge $\alpha$ of
$\Gamma $, the edge labelled $Z$ in Figure 1, to produce $\Gamma _\alpha$ as in the
figure; the fattening and embedding of $\Gamma_\alpha$ in $F$ is determined from
that of $\Gamma$ in the natural way. Furthermore, an enumeration of the edges of
$\Gamma $ induces an enumeration of the edges of $\Gamma _\alpha$ in the natural
way, where the vertical edge labelled $Z$ in Figure 1 corresponds to the horizontal
edge labelled $-Z$.  We say that $\Gamma _\alpha$ arises from $\Gamma$ by a
{\it Whitehead move} along edge $\alpha$.  We shall also write $\Gamma
_{\alpha\beta}=(\Gamma _\alpha )_\beta$, for any two indices $\alpha ,\beta$ of
edges, to denote the result of first performing a move along $\alpha$ and then along
$\beta$; in particular, $\Gamma _{\alpha \alpha}=\Gamma$ for any index $\alpha$.

\vspace{10pt}

\setlength{\unitlength}{1.5mm}%
\begin{picture}(50,27)(-12,48)
\thicklines
\put(28,70){\line( 1,-2){ 4}}
\put(32,62){\line( 1, 0){28}}
\put(60,62){\line( 1, 2){ 4}}
\put(60,62){\line( 1,-2){ 4}}
\put(32,62){\line(-1,-2){ 4}}
\thinlines
\put(18,62){\vector(-1, 0){  0}}
\put(18,62){\vector( 1, 0){ 5}}
\thicklines
\put(10,54){\line( 2,-1){ 8}}
\put(10,70){\line( 0,-1){16}}
\put(10,54){\line(-2,-1){ 8}}
\put( 2,74){\line( 2,-1){ 8}}
\put(10,70){\line( 2, 1){ 8}}
\put( 4,74){\makebox(0,0)[lb]{$A$}}
\put(16,74){\makebox(0,0)[rb]{$B$}}
\put(12,62){\makebox(0,0)[lc]{$Z$}}
\put(16,50){\makebox(0,0)[rt]{$C$}}
\put( 4,50){\makebox(0,0)[lt]{$D$}}
\put(30,54){\makebox(0,0)[lt]{$D - \phi(-Z)$}}
\put(62,54){\makebox(0,0)[rt]{$C+\phi(Z)$}}
\put(62,69){\makebox(0,0)[rb]{$B-\phi(-Z)$}}
\put(30,69){\makebox(0,0)[lb]{$A+\phi(Z)$}}
\put(47,64){\makebox(0,0)[cb]{$-Z$}}
\end{picture}

\centerline{\bf Figure~1-Whitehead move on shear coordinates}

\hspace{10pt}

\begin{proposition} {\rm \cite{ChF}}\label{propcase}
Setting $\phi (Z)={\rm log}(e^Z+1)$ and adopting the notation of Figure~1
for shear coordinates of nearby edges, the effect of a
Whitehead move is illustrated in the figure, viz.,
\be
W_Z\,:\ (A,B,C,D,Z)\to (A+\phi(Z), B-\phi(-Z), C+\phi(Z), D-\phi(-Z), -Z)
\label{abc}
\ee
In the various cases where the edges are not distinct
and identifying an edge with its shear coordinate in the obvious notation we have:
if $A=C$, then $A'=A+2\phi(Z)$;
if $B=D$, then $B'=B-2\phi(-Z)$;
if $A=B$ (or $C=D$), then $A'=A+Z$ (or $C'=C+Z$);
if $A=D$ (or $B=C$), then $A'=A+Z$ (or $B'=B+Z$).
\end{proposition}

\noindent{\sl Sketch of proof}~ Assume that $e$ is the diagonal of a
quadrilateral with consecutive sides $a,b,c,d$, where $e$ separates $a,b$ from
$c,d$.  Identifying an edge with its lambda length, the shear coordinate is given by
$Z={\rm log}~{{bd}\over{ac}}$, i.e., ${{bd}\over{ac}}$ is the required cross-ratio
\cite{Penn1}.  Furthermore, if $f$ is the lambda length of the other diagonal, then the
lambda lengths satisfy Ptolemy's relation $ef=ac+bd$ \cite{Penn1}, and the transformation
laws for shear coordinates in the proposition are readily derived from this
either in the surface for (\ref{abc}) or in the universal cover of the surface in the various cases.
~~~~~\hfill{\it q.e.d.}

\vskip .1in

Insofar as hyperbolic lengths of geodesics are well-defined invariants of homotopy classes
of curves in $F$, these lengths must be invariant under Whitehead moves, so we have the following

\begin{lemma} \label{lem1}
Transformation~{\rm(\ref{abc})} preserves
the traces of products over paths {\rm(\ref{geod})}.
\end{lemma}

\subsubsection{Weil--Petersson form}

${\cal T}_H(F)$ supports its canonical Weil--Petersson Poisson structure, which has a very
simple form in shear coordinates.

\begin{theorem}\label{th-WP}~{\rm\cite{Fock1}} In the coordinates $(Z_\alpha )$ on any fixed spine,
the Weil--Petersson bracket $B_{{\mbox{\tiny WP}}}$ is given by
\be
\label{WP-PB}
B_{{\mbox{\tiny WP}}}
= \sum_{v} \sum_{i=1}^{3}\frac{\partial}{\partial Z_{v_i}}\wedge
\frac{\partial}{\partial Z_{v_{i+1}}},
\ee
where the sum is taken over all vertices~$v$ and $v_i$, \
$i=1,2,3\ \hbox{\rm mod}\ 3$, are the labels of the cyclically ordered
half-edges incident on this vertex.
\end{theorem}

The proof~\cite{Fock1} relies on the independence of this form under Whitehead moves
as in \cite{Penn3}.  Indeed, the equivalent expression for the Weil--Petersson K\"ahler
two-form in the punctured case was first given in lambda length coordinates in \cite{Penn2}
starting from Wolpert's formula \cite{Wol}, and the formula  in shear coordinates follows
from direct calculation using the expression $Z={\rm log}~{{bd}\over{ac}}$. (There
is more to this geometrically, however, and one must show \cite{Bon3} that the same expression
is the Weil--Petersson form for surfaces with boundary.)

The set of Casimir functions is described by the following proposition, whose proof is given in
Appendix~B.

\begin{proposition}\label{prop12}
The center of the Poisson algebra {\rm(\ref{WP-PB})} is generated by
elements of the form $\sum Z_\alpha$, where the sum is over all edges
of $\Gamma $ in a boundary component of $F(\Gamma )$ and the sum is taken with multiplicity.
\end{proposition}

\subsubsection{Mapping class group description using graphs.}
Recall that the mapping class group $MC(F)$ of an open surface $F$ is the group of
homotopy classes of orientation-preserving homeomorphisms of $F$.  No special
constraints are imposed by the circle boundary components, i.e., a homeomorphism
must fix each boundary component only setwise, and the homotopies must likewise fix
each boundary component only setwise.  Thus, if $F$ has $b$ boundary component
circles, $p$ punctures, and genus $g$, then $MC(F)\approx MC(F_g^{b+p})$, so we
generally write $MC_g^s=MC(F)$ for any surface of genus $g$ with $s$ boundary
components.  In this section, we establish the combinatorial presentation of
$MC_g^s$ associated with the cell decomposition of decorated Teichm\"uller space.

Recall that a cell in the decomposition of $\widetilde{\cal T}_g^s$ is described by the
homotopy classes of an embedding $\iota :\Gamma \to F$ of a fatgraph $\Gamma $ as a
spine of $F$. $MC_g^s$ acts on the set of homotopy classes of such embeddings by
post-composition, and the cell decomposition of decorated Teichm\"uller space
$\widetilde{\cal T}_g^s$ descends to an orbifold cell decomposition of
$\widetilde {\cal M}_g^s=\widetilde{\cal T}_g^s/MC(F)$.

The {\it modular groupoid} $MG_g^s=MG(F)$ is the fundamental path groupoid of
$\widetilde{\cal M}_g^s$, and $MC_g^s$ arises as the subgroup of paths based at any
point.  Specifically, consider the dual graph ${\cal G}_g^s={\cal G} (F)$ of the
codimension-two skeleton of this decomposition of $\widetilde{\cal M}_g^s$ (where
there is one vertex for each top-dimensional cell, edges correspond to Whitehead
moves, and two-dimensional cells correspond to pairs of homotopic paths in the
one-skeleton which are homotopic to real endpoints in $\widetilde{\cal M}_g^s$.)
The fundamental path groupoid of ${\cal G}_g^s$ is the modular groupoid, and in
particular, $MC_g^s$ is the stabilizer in $MG_g^s$ of any vertex of ${\cal G}_g^s$.

We may think of a Whitehead move along edge $\alpha$ of fatgraph $\Gamma $ producing
another fatgraph $\Gamma _\alpha$ as an ordered pair $(\Gamma ,\Gamma _\alpha )$,
i.e., an oriented edge of ${\cal G}_g^s$.  Letting $[\Gamma _1,\Gamma _2 ]$ denote
the $MC_g^s$-orbit of a pair $(\Gamma _1,\Gamma _2 )$ by the diagonal action of the
mapping class group, the natural composition descends to a well-defined product
$$
[\Gamma _1,\Gamma _3]=[\Gamma _1,\Gamma _2][\Gamma _2,\Gamma _3].
$$

\begin{theorem} {\rm \cite{Penn3},\cite{Penn2}}\label{modgrpthm}
The modular groupoid $MG_g^s$ is generated by Whitehead moves and relabelings by
fatgraph symmetries. A complete list of relations in $MG_g^s$ is given by
relabelings under fatgraph symmetries together with the two following relations.

\leftskip .3in

\noindent{\bf Commutativity} \ \ \ If $\alpha$ and $\beta$ are two edges with no
common endpoints, then
$$
[\Gamma _{\alpha\beta}, \Gamma _\alpha][\Gamma _\alpha, \Gamma ] =
[\Gamma _{\alpha\beta},
\Gamma _\beta][\Gamma _\beta, \Gamma ].
$$

\noindent{\bf Pentagon} If $\alpha$ and $\beta$ share exactly one common endpoint, then
(see Figure 2 drawn for the dual graph)
$$
[\Gamma ,\Gamma _\alpha][\Gamma _\alpha,
\Gamma _{\beta\alpha}][\Gamma _{\beta\alpha},\Gamma _{\alpha\beta}][
\Gamma _{\alpha\beta},\Gamma _\beta][\Gamma _\beta,\Gamma ] =1.
$$

\leftskip=0ex

\noindent Furthermore, the expression in Theorem~\ref{th-WP} for the
Weil--Petersson form is invariant under Whitehead moves.
\end{theorem}

\centerline{\unitlength 0.45mm
\begin{picture}(0,0)(100,-145)
\put(0,0){\line(2,1){ 40}}
\put(40,20){\line(2,-1){ 40}}
\put(80,0){\line(-1,-2){ 20}}
\put(60,-40){\line(-1,0){40}}
\put(20,-40){\line(-1,2){ 20}}
\put(0,0){\line(1,0){80}}
\put(0,0){\line(3,-2){60}}
\put(20,13){\makebox(0,0)[rb]{$E_0$}}
\put(60,13){\makebox(0,0)[lb]{$D_0$}}
\put(8,-21){\makebox(0,0)[rc]{$A_0$}}
\put(72,-21){\makebox(0,0)[lc]{$C_0$}}
\put(40,-43){\makebox(0,0)[ct]{$B_0$}}
\put(40,3){\makebox(0,0)[cb]{$Y_0$}}
\put(29,-22){\makebox(0,0)[ct]{$X_0$}}
\put(45,-50){\vector(1,-2){5}}
\end{picture}
\begin{picture}(0,0)(20,-200)
\put(0,0){\line(2,1){ 40}}
\put(40,20){\line(2,-1){ 40}}
\put(80,0){\line(-1,-2){ 20}}
\put(60,-40){\line(-1,0){40}}
\put(20,-40){\line(-1,2){ 20}}
\put(0,0){\line(3,-2){60}}
\put(40,20){\line(1,-3){20}}
\put(20,13){\makebox(0,0)[rb]{$C_4$}}
\put(60,13){\makebox(0,0)[lb]{$B_4$}}
\put(8,-21){\makebox(0,0)[rc]{$D_4$}}
\put(72,-21){\makebox(0,0)[lc]{$A_4$}}
\put(40,-43){\makebox(0,0)[ct]{$E_4$}}
\put(53,-8){\makebox(0,0)[lc]{$X_4$}}
\put(29,-22){\makebox(0,0)[ct]{$Y_4$}}
\put(1,-29){\vector(-1,-1){10}}
\end{picture}
\begin{picture}(0,0)(-60,-145)
\put(0,0){\line(2,1){ 40}}
\put(40,20){\line(2,-1){ 40}}
\put(80,0){\line(-1,-2){ 20}}
\put(60,-40){\line(-1,0){40}}
\put(20,-40){\line(-1,2){ 20}}
\put(40,20){\line(1,-3){20}}
\put(40,20){\line(-1,-3){20}}
\put(20,13){\makebox(0,0)[rb]{$A_3$}}
\put(60,13){\makebox(0,0)[lb]{$E_3$}}
\put(8,-21){\makebox(0,0)[rc]{$B_3$}}
\put(72,-21){\makebox(0,0)[lc]{$D_3$}}
\put(40,-43){\makebox(0,0)[ct]{$C_3$}}
\put(53,-8){\makebox(0,0)[lc]{$Y_3$}}
\put(27,-8){\makebox(0,0)[rc]{$X_3$}}
\put(13,21){\vector(-1,1){10}}
\end{picture}
\begin{picture}(0,0)(-30,-65)
\put(0,0){\line(2,1){ 40}}
\put(40,20){\line(2,-1){ 40}}
\put(80,0){\line(-1,-2){ 20}}
\put(60,-40){\line(-1,0){40}}
\put(20,-40){\line(-1,2){ 20}}
\put(40,20){\line(-1,-3){20}}
\put(20,-40){\line(3,2){60}}
\put(20,13){\makebox(0,0)[rb]{$D_2$}}
\put(60,13){\makebox(0,0)[lb]{$C_2$}}
\put(8,-21){\makebox(0,0)[rc]{$E_2$}}
\put(72,-21){\makebox(0,0)[lc]{$B_2$}}
\put(40,-43){\makebox(0,0)[ct]{$A_2$}}
\put(27,-8){\makebox(0,0)[rc]{$Y_2$}}
\put(52,-23){\makebox(0,0)[ct]{$X_2$}}
\put(58,19){\vector(1,2){6}}
\end{picture}
\begin{picture}(0,230)(70,-65)
\put(0,0){\line(2,1){ 40}}
\put(40,20){\line(2,-1){ 40}}
\put(80,0){\line(-1,-2){ 20}}
\put(60,-40){\line(-1,0){40}}
\put(20,-40){\line(-1,2){ 20}}
\put(20,-40){\line(3,2){60}}
\put(0,0){\line(1,0){ 80}}
\put(20,13){\makebox(0,0)[rb]{$B_1$}}
\put(60,13){\makebox(0,0)[lb]{$A_1$}}
\put(8,-21){\makebox(0,0)[rc]{$C_1$}}
\put(72,-21){\makebox(0,0)[lc]{$E_1$}}
\put(40,-43){\makebox(0,0)[ct]{$D_1$}}
\put(52,-23){\makebox(0,0)[ct]{$Y_1$}}
\put(40,3){\makebox(0,0)[cb]{$X_1$}}
\put(82,-21){\vector(1,0){15}}
\end{picture}}
\vspace{5mm}
\centerline{\bf Figure~2-pentagon identity}
\vspace{5mm}

\noindent{\sl Proof of Theorem} The first parts are immediate consequences of the
cell decomposition. Specifically, by connectivity of $\widetilde{\cal T}_g^s$, any
two points can be joined by a smooth path, which we may put into general position
with respect to the codimension-one faces; this proves the first part. For the
second part, a homotopy between edge-paths in $\widetilde{\cal M}_g^s$ can likewise
be put into general position with respect to the codimension-two faces; there are
two possibilities for a pair of edges depending upon whether their vertices are
disjoint or not, corresponding respectively to the commutativity and pentagon
relations, proving the second part. The invariance of the expression for the
Weil--Petersson form under Whitehead moves is a direct calculation in lambda
lengths \cite{Penn4} or shear coordinates using (\ref{abc}).~~~~\hfill{\it q.e.d.}

\clearpage

\subsection{Poisson algebras of geodesic functions}

The algebra generated (with multiplication and with the
Weil--Petersson Poisson bracket) by the functions $\{G_\gamma\}$ (\ref{geod})
was first studied by W. Goldman~\cite{Gold}.

\subsubsection{Multicurves}

In the sequel, disjointly embedded families of geodesics will play a special role as
they constitute a basis for the algebra of observables in both the classical case
considered here and the quantum case discussed in Section 3. The homotopy class of such a
family is called a {\it multiple curve}.  A {\it multicurve} is  multiset based on the set of curves in a
multiple curve.

\begin{defin} \label{def1}{\rm
Consider the homotopy class of a finite collection
$C=\{ \gamma _1,\ldots ,\gamma _n\}$
of disjointly embedded (unoriented) simple closed curves $\gamma _i$ in a
topological surface $F$, where $C$ need {\sl not} be a mutiple curve. A {\it generalized multicurve}
(GMC) $\hat C$ in $F$
is a multiset based on $C$; one thus imagines $s_i\geq 1$ parallel copies of components
of $C$, or in other words, positive integral weights $s_i$ on each component of $C$,
where $s_i$ is the multiplicity of $\gamma _i$ in $\hat C$. Further, given a
hyperbolic structure on $F$, we associate to $\hat C$ the product $G_{\hat
C}=G_{\gamma_1}^{s_1}\dots G_{\gamma_n}^{s_n}$ of geodesic operators (\ref{geod}) of
all geodesics constituting a GMC; these operators Poisson commute in the classical
case since the components of $C$ are disjoint. In particular, a GMC containing a
contractible component (of length zero) is twice the GMC with this curve removed. }
\end{defin}

An edge-path on a spine
$\Gamma\subseteq F$ or its corresponding geodesic $\gamma$ in $F$
is said to be {\em graph simple\/} with respect to $\Gamma$ if it does not pass more
than once through any edge of $\Gamma$. Obviously, the set of graph simple geodesics
depends upon $\Gamma$ and is {\em not\/} invariant under Whitehead moves.
Nevertheless, this notion will be useful in what follows.

\subsubsection{Classical skein relation}
The trace relation $\tr(AB)+\tr(AB^{-1})-\tr A\cdot\tr B=0$ for
arbitrary $2\times 2$ matrices~$A$ and~$B$ with unit determinant allows one to
``disentangle'' any product of
geodesic functions,
i.e., express it uniquely as a finite linear combination of GMCs.
Introducing the additional factor $\#G$ to be the total number of components
in a GMC, we can then
uniformly present the classical skein relation as
\be
\setlength{\unitlength}{.8mm}%
\begin{picture}(90,20)(0,50)
\thicklines
\put(-20,60){\makebox(0,0){$(-1)^{\#G}$}}
\put(-10,50){\line(1, 1){9}}
\put( 1, 61){\line(1, 1){9}}
\put(-10,70){\line(1,-1){20}}
\put(18,60){\makebox(0,0){$+$}}
\put(30,60){\makebox(0,0){$(-1)^{\#G}$}}
\put(40,60){\oval(20,20)[r]}
\put(63,60){\oval(20,20)[l]}
\put(66,60){\makebox(0,0){$+$}}
\put(78,60){\makebox(0,0){$(-1)^{\#G}$}}
\put(95,48){\oval(20,20)[t]}
\put(95,72){\oval(20,20)[b]}
\put(110,60){\makebox(0,0){$=0,$}}
\end{picture}
\label{skeinclass}
\ee

\subsubsection{Poisson brackets for geodesic functions.}

Turning attention now to the Poisson structure, two geodesic functions Poisson
commute if the underlying geodesics are disjointly embedded. By the Leibnitz rule
for the Poisson bracket, it suffices to consider only ``simple'' intersections of
pairs of geodesics with respective geodesic functions
$G_1$ and $G_2$ of the form
\bea
\label{G1}
G_1&=&\tr^1\dots X_C^1R^1X_Z^1L^1X_A^1\dots,
\\
\label{G2}
G_2&=&\tr^2\dots X_B^2 L^2 X_Z^2R^2X_D^2\dots,
\eea
where the superscripts $1$ and~$2$ pertain to operators and traces in two different
matrix spaces.

The bracket between $X_C^1$ and $X_B^2$ possesses a simple $r$-matrix structure
\be
\{X_C^1, X_B^2\}=\frac{1}{4}(-1)^{i+j}(e^1_{ii}\otimes e^2_{jj})
X_C^1\otimes X_B^2,
\ee
where the ``elementary'' matrix $e_{ij}$ has entry
unity in its $i$th row and $j$th column and zero otherwise.
Direct calculations then give
\be
\label{Goldman}
\{G_1,G_2\}=\frac{1}{2}(G_{\mbox{\tiny H}}-G_{\mbox{\tiny I}}),
\ee
where $G_{\mbox{\tiny I}}$ corresponds to the geodesic that is obtained by
erasing the edge~$Z$ and joining together the edges ``$A$'' and ``$D$'' as
well as ``$B$'' and ``$C$'' in a natural way as illustrated in the middle diagram
in (\ref{skeinclass}); $G_{\mbox{\tiny H}}$ corresponds to the geodesic that
passes over the edge~$Z$ twice, so it has the form
$\tr \dots X_CR_ZR_D\dots$ $\dots X_BL_ZL_A\dots$
as illustrated in the rightmost diagram in (\ref{skeinclass}).
These relations were first obtained
in~\cite{Gold} in the continuous parametrization (the classical Turaev--Viro
algebra).

\vskip .2in

\noindent {\bf Torus Example}~~For the torus, ${\cal T}_H(F_1^1)$ has three
generators $X,Y,Z$, where
$$
\{X,Y\}=\{Y,Z\}=\{Z,X\}=2,
$$
corresponding to the combinatorially unique cubic spine and the Casimir element is $X+Y+Z$.
The geodesic functions for the three graph simple geodesics are
\bea
\nonumber &{}&G_X=\tr
LX_YRX_Z=\e^{-Y/2-Z/2}+\e^{-Y/2+Z/2}+\e^{Y/2+Z/2},
\\
\label{torus-classic} &{}&G_Y=\tr
RX_XLX_Z=\e^{-Z/2-X/2}+\e^{-Z/2+X/2}+\e^{Z/2+X/2},
\\
\nonumber &{}&G_Z=\tr
RX_YLX_X=\e^{-X/2-Y/2}+\e^{-X/2+Y/2}+\e^{X/2+Y/2}. \eea
Introducing the geodesic function
$$
\wtd{G_Z}=\tr L_YR_ZR_XL_Z=
\e^{-X/2-Y/2-Z}+\e^{X/2-Y/2}(\e^{-Z}+\e^Z+2)+\e^{X/2+Y/2+Z},
$$
obtained from $G_Z$ by a Whitehead move,
we find that
$\{G_X,G_Y\}={\wtd{G}_Z}/2-G_Z/2$, and because relation
(\ref{skeinclass}) implies that $G_XG_Y=G_Z+{\wtd G}_Z$, we have
\be
\{G_X,G_Y\}=\frac12 G_XG_Y-G_Z,\\
\label{GXGYclass}
\ee
plus the cyclic permutations in $X,Y,Z$, i.e., the classical
Poisson algebra closes in the algebraic span of the geodesic functions
$\{G_X,G_Y,G_Z\}$.

\subsubsection{Poisson geodesic algebras for higher genera}
In order to generalize the torus example, we must find a graph
on which graph simple geodesics constitute a convenient algebraic basis.
Such a graph is illustrated in Figure~3, where
$m$ edges pairwise connect two horizontal line
segments. Graph simple closed geodesics in this picture are those and
only those that pass through exactly two different
``vertical'' edges, and they are therefore enumerated by
ordered pairs of edges; we
denote the corresponding geodesic functions ${\cal G}_{ij}$ where $i<j$. The Poisson algebra
for the functions ${\cal G}_{ij}$ is described by
\be
\{{\cal G}_{ij},{\cal G}_{kl}\}=\left\{\begin{array}{l}
0,\quad j<k,\\
0,\quad k<i,\ j<l,\\
{\cal G}_{ik}{\cal G}_{jl}-{\cal G}_{kj}{\cal G}_{il},\quad i<k<j<l,\\
\frac12 {\cal G}_{ij}{\cal G}_{jl}-{\cal G}_{il},\quad j=k,\\
{\cal G}_{il}-\frac12 {\cal G}_{ij}{\cal G}_{il},\quad i=k,\ j<l\\
{\cal G}_{ik}-\frac12 {\cal G}_{ij}{\cal G}_{kj},\quad j=l,\ i<k.
\end{array}
\right. \label{P-geod} \ee

The graph in Figure~3 has genus $\frac{m}{2}-1$
and {two} faces (holes) if $m$ is even and genus $(m-1)/2$ and {
one} face (hole) if~$m$ is odd.  Such geodesic bases for $m$ even were considered in~\cite{NR}. The
Poisson algebras of geodesics obtained there {\it coincide
exactly\/} with (\ref{P-geod}). These are the so-called $so_q(m)$
algebras whose representations were constructed in \cite{Klim}.

In the mathematical literature, this
algebra has also appeared as the Poisson algebra
of the monodromy data (Stokes matrices) of certain matrix differential
equations~\cite{Ugaglia}
and on the symplectic
groupoid of upper-triangular matrices~${\cal G}$~\cite{Bondal}. These matrices
have entries given by unity on the main diagonal (i.e., we set ${\cal G}_{ii}\equiv1$)
and the entries ${\cal G}_{ij}$ above it.
For $m\times m$-matrices, there are $\left[\frac{m}{2}\right]$
central elements of this algebra generated by the polynomial
invariants $f_{\cal G}(\lambda)\equiv
\det({\cal G}+\lambda {\cal G}^{T})=\sum_{}^{}f_i({\cal G})\lambda^i$.
The total Poisson dimension $d$ of algebra
(\ref{P-geod}) is $\frac{m(m-1)}{2}-\left[\frac{m}{2}\right]$,
and for $m=3,4,5,6,\dots$ we have $d=2,4,8,12,\dots\,$. The
dimensions of the corresponding Teichm\"uller spaces are
$D=2,4,8,10,\dots\,$, so we see that the Teichm\"uller spaces are
embedded as the {Poisson leaves} in the algebra (\ref{P-geod}).

\be
\setlength{\unitlength}{.4mm}%
\begin{picture}(90,120)(15,-15)
\thicklines
\put(0,95){\line(1, 0){125}}
\qbezier(0,95)(-25,95)(-10,85)
\qbezier(125,95)(150,95)(135,85)
\put(10,85){\line(1, 0){10}} \put(30,85){\line(1, 0){15}}
\put(55,85){\line(1, 0){15}} \put(80,85){\line(1, 0){15}}
\put(105,85){\line(1, 0){10}} \put(0,0){\line(1, 0){125}}
\qbezier(0,0)(-25,0)(-10,10)
\qbezier(125,0)(150,0)(135,10)
\put(10,10){\line(1, 0){10}} \put(30,10){\line(1, 0){15}}
\put(55,10){\line(1, 0){15}} \put(80,10){\line(1, 0){15}}
\put(105,10){\line(1, 0){10}} \put(-10,10){\line(5, 3){125}}
\put(10,10){\line(5, 3){125}} \put(20,10){\line(1, 1){15}}
\put(30,10){\line(1, 1){31}} \put(95,85){\line(-1, -1){31}}
\put(105,85){\line(-1, -1){15}} \put(45,10){\line(1, 3){7.5}}
\put(55,10){\line(1, 3){11.4}} \put(70,85){\line(-1, -3){11.4}}
\put(80,85){\line(-1, -3){7.5}} \put(-10,85){\line(5, -3){62.5}}
\put(10,85){\line(5, -3){25}} \put(115,10){\line(-5, 3){25}}
\put(135,10){\line(-5, 3){62.5}} \put(20,85){\line(1, -1){34.6}}
\put(30,85){\line(1, -1){21.5}} \put(95,10){\line(-1, 1){21.5}}
\put(105,10){\line(-1, 1){34.6}} \put(45,85){\line(1, -3){11.7}}
\put(55,85){\line(1, -3){7.5}} \put(70,10){\line(-1, 3){7.5}}
\put(80,10){\line(-1, 3){11.7}}
\thinlines
\qbezier(25,90)(20,90)(25,85)
\qbezier(100,90)(105,90)(100,85)
\qbezier(25,5)(20,5)(25,10)
\qbezier(100,5)(105,5)(100,10)
\put(25,90){\line(1,0){75}}
\put(25,5){\line(1,0){75}}
\put(25,85){\line(1,-1){22}}
\put(100,10){\line(-1,1){22}}
\put(100,85){\line(-1,-1){14}}
\put(25,10){\line(1,1){14}}
\put(20,105){\makebox(0,0){$i$}}
\put(65,105){\makebox(0,0){$\cdots$}}
\put(105,105){\makebox(0,0){$j$}}
\put(20,-10){\makebox(0,0){$j$}}
\put(65,-10){\makebox(0,0){$\cdots$}}
\put(105,-10){\makebox(0,0){$i$}}
\end{picture}
\label{octopus}
\ee

\centerline{\bf Figure~3-the special spine for higher genera}

\hspace{10pt}

\newsection{Quantization}
A quantization of a Poisson manifold, which is equivariant under the action of a discrete
group $\cal D$,
is a family of $*$-algebras ${\cal A}^\hbar$ depending on a positive
real parameter $\hbar$ with
$\cal D$ acting by outer automorphisms and having the
following properties:

\leftskip .3in

\vskip .1in\noindent {\bf 1.} (Flatness.)  All algebras are isomorphic (noncanonically)
as linear spaces.

\vskip .1in\noindent {\bf 2.}
(Correspondence.) For $\hbar=0$, the algebra is isomorphic as a $\cal
D$-module to the $*$-algebra of complex-valued functions on the Poisson
manifold.

\vskip .1in\noindent {\bf 3.} (Classical Limit.) The Poisson bracket on ${\cal A}^0$ given
by $\{a_1, a_2\} = \lim_{\hbar \rightarrow 0}\frac{[a_1,a_2]}{\hbar}$ coincides with the
Poisson bracket given by the Poisson structure of the manifold.

\leftskip=0ex

\subsection{Quantizing Teichm\"uller spaces}

Here we construct a quantization ${\cal T}^\hbar (F)$ of the Teichm\"uller space ${\cal T}_H(F)$
that is equivariant with respect to the action of the mapping class group ${\cal D}=MC(F)$.

Fix a cubic fatgraph $\Gamma$ as spine of $F$, and let ${\cal T}^\hbar={\cal T}^\hbar(\Gamma)$
be the algebra generated by $Z_\alpha^\hbar$, one generator for each
unoriented edge $\alpha$ of $\Gamma$, with relations
\be
\label{comm}
[Z^\hbar_\alpha, Z^\hbar_\beta ] = 2\pi i\hbar\{z_\alpha, z_\beta\}
\ee
(cf.\ (\ref{WP-PB})) and the
$*$-structure
\be
(Z^\hbar_\alpha)^*=Z^\hbar_\alpha,
\ee
where $z_\alpha$  and
$\{\cdot,\cdot\}$ denotes the respective coordinate functions
and the Poisson bracket on the classical
Teichm\"uller space.
Because of (\ref{WP-PB}), the righthand side of (\ref{comm}) is a constant
taking only five values $0$, \ $\pm 2\pi i \hbar$, and $\pm 4 \pi i \hbar$
depending upon the coincidences of endpoints of edges labelled $\alpha$ and $\beta$.

\begin{lemma}\label{nondeg}
The center ${\cal Z}^h$ of the algebra ${\cal T}^h$ is generated by the sums
$\sum_{\alpha\in I}{Z^\hbar_\alpha}$ over all edges $\alpha\in I$ surrounding
a given boundary component, and the Poisson structure is non-degenerate on
the quotient ${\cal T}^h/{\cal Z}^h$.
\end{lemma}

A standard Darboux-type theorem for non-degenerate Poisson structures then gives
the following result.

\begin{cor}\label{corr11}
There is a basis for ${\cal T}^h/{\cal Z}^h$
given by operators $p_i,q_i$, for $i=1,\ldots ,6g-6+2s$
satisfying the
standard commutation relations $[p_i,q_j]=2\pi i
\hbar\delta_{ij}$, \ $i,j=1,\dots,3g-3+s$.
\end{cor}

\noindent Not only is the proof of Lemma \ref{nondeg} given in Appendix~B, but also an algorithm
for diagonalizing this Poisson structure is described there.

\vskip .2in

Now, define the Hilbert space ${\cal H}$ to be the set of all $L^2$ functions in the
$q$-variables and let each $q$-variable act by multiplication and each corresponding
$p$ variable act by differentiation, $p_i=2\pi i\hbar~{{\partial}\over{\partial p_i}}$. For
different choices of diagonalization of non-degenerate Poisson structures, these
Hilbert spaces are canonically isomorphic.

\vskip .1in

\noindent{\bf Torus Example}
In the case of the bordered torus, we have three generators, $X^\hbar$, $Y^\hbar$, and $Z^\hbar$,
the commutation relations (\ref{comm}) have the form
$[X^\hbar,Y^\hbar]=[Y^\hbar,Z^\hbar]=[Z^\hbar,X^\hbar]=4\pi i \hbar$, and
the single central element is $X^\hbar+Y^\hbar+Z^\hbar$. In the Darboux-type
representation, we can identify, e.g., $(X^\hbar+Y^\hbar)/2$ with $q$ and
$(-X^\hbar+Y^\hbar)/2$ with $(-2\pi i \hbar)\partial/\partial q$.

\vskip .1in

On the level of the modular groupoid as a category where all morphisms are invertible
and any two objects are related by a morphism, we have
constructed one $*$-algebra per object.
In order to describe the ${\cal D}$-equivariance we must associate a
homomorphism of the corresponding $*$-algebras to any morphism in the modular
groupoid.
For this, we associate a morphism of algebras to
any Whitehead move and must verify that the relations in Theorem~\ref{modgrpthm} are satisfied.

We now define the {\it quantum Whitehead move} or {\it flip} along an edge of $\Gamma$ by Eq. (\ref{abc}) using
the (quantum) function
\begin{equation} \label{phi}
\phi(z)\equiv \phi^\hbar(z) =
-\frac{\pi\hbar}{2}\int_{\Omega} \frac{e^{-ipz}}{\sinh(\pi p)\sinh(\pi \hbar
p)}dp,
\end{equation}
where the contour $\Omega$ goes along
the real axis bypassing the origin from above.
The function~(\ref{phi})
is Faddeev's generalization ~\cite{Faddeev} of the quantum dilogarithm.

\vskip .1in

\begin{proposition}
For each unbounded self-adjoint operator $Z$ on ${\cal H}$, $\phi (Z)$ is a well-defined
unbounded self-adjoint operator on
${\cal H}$.
\end{proposition}

\vskip .1in

\noindent{\sl Proof}~
The function $\phi^\hbar(Z)$ satisfies the relations (see~\cite{ChF})
$$
\phi^\hbar(Z)-\phi^\hbar(-Z)=Z,
$$
$$
\phi^\hbar(Z+i\pi\hbar)-\phi^\hbar(Z-i\pi\hbar)=\frac{2\pi i\hbar}{1+\e^{-Z}},
$$
$$
\phi^\hbar(Z+i\pi)-\phi^\hbar(Z-i\pi)=\frac{2\pi i}{1+\e^{-Z/\hbar}}
$$
and is meromorphic in the complex plane with the poles at the
points $\{\pi i(m+n\hbar),\ m,n\in {\Bbb Z}_+\}$ and
$\{-\pi i(m+n\hbar),\ m,n\in {\Bbb Z}_+\}$.

The function $\phi^\hbar(Z)$ is therefore holomorphic in the strip
$|\hbox{Im\,}Z|<\pi\hbox{\,min\,}(1,\hbox{Re\,}\hbar)-\epsilon$ for any $\epsilon>0$,
so we need only its asymptotic behavior as
$Z\in{\Bbb R}$ and $|Z|\to\infty$, for which we have (see, e.g., \cite{Kashaev3})
\be
\biggr.\phi^\hbar(Z)\biggl|_{|Z|\to\infty}=(Z+|Z|)/2+O(1/|Z|).
\label{phih-as}
\ee
Therefore, the function $\phi^\hbar(Z)$ increases as $|Z|$
at infinity and represents an operator in ${\cal H}$ by the functional calculus \cite{Lap}.~~~~~\hfill{\it
q.e.d.}.

\begin{theorem}\label{th-Q}

The family of algebras ${\cal T}^\hbar={\cal T}^\hbar(\Gamma )$ is a
quantization of ${\cal T}_H(F)$ for any cubic fatgraph spine
$\Gamma$ of $F$, that is:

\leftskip .4in

\vskip .1in\noindent {\bf 1}. In the limit $\hbar \mapsto 0$, morphism {\rm(\ref{abc})} using
(\ref{phi}) coincides
with the classical morphism {\rm(\ref{abc})} where $\phi(Z)=\log(1+\e^Z)$;

\vskip .1in\noindent {\bf 2.} Morphism {\rm(\ref{abc})} using {\rm(\ref{phi})} is indeed
a morphism of $*$-algebras;

\vskip .1in\noindent {\bf 3}. A flip $W_Z$ satisfies $W_Z^2=I$, cf.
{\rm(\ref{abc})}, and flips satisfy the commutativity relation;

\vskip .1in\noindent {\bf 4}. Flips satisfy the pentagon relation.\footnote{This
result was independently obtained by R.~M.~Kashaev~\cite{kashaev}.}

\leftskip=0ex

\vskip .1in

Furthermore,

\leftskip .4in

\vskip .1in

\noindent {\bf 5}. The morphisms ${\cal T}^\hbar(\Gamma)\rightarrow {\cal
T}^{1/\hbar}(\Gamma)$ given by $Z^\hbar_\alpha \mapsto Z^{1/\hbar}_\alpha$
commute with morphisms {\rm(\ref{abc})}.

\end{theorem}

\leftskip=0ex

\noindent {\em Sketch of the proof} ~Property 1 follows since $\lim_{\hbar \rightarrow
0}\phi^\hbar(z) =
\log(\e^z + 1)$, and Property 3 is obvious.

In order to prove Property 2, we must first verify that $[A+\phi^\hbar(Z),
B-\phi^\hbar(-Z)] = 0$ and  $[A+\phi^\hbar(Z), D - \phi^\hbar(-Z)]= -2\pi i
\hbar$ (since the other relations are obviously satisfied), which
follows from the identity $\phi^\hbar(z)-\phi^\hbar(-z)=z$.

For Property 5, we must verify that the morphism ${\cal
T}^\hbar(\Gamma)\rightarrow {\cal T}^{1/\hbar}(\Gamma)$ commutes with a flip, that is,
$(A+\phi^\hbar(Z))/\hbar = A/\hbar +\phi^\hbar(Z/\hbar)$,
$(B-\phi^{1/\hbar}(-Z))/\hbar = B/\hbar -\phi^\hbar(-Z/\hbar)$,~etc.,
which follows from $\phi^\hbar(z)/\hbar = \phi^{1/\hbar}(z/\hbar)$.

Turning finally to the most nontrivial Property 4, we may reformulate it as follows.
There are seven generators involved in the sequence of flips
depicted in Figure~2 for the dual cell decomposition,
which are denoted
$A,B,C,D,E,X,Y$ as in the figure. As a result of
a flip, the piece of graph shown in Figure~1 just gets
cyclically rotated. Denote by $A_i,B_i,C_i,D_i,E_i,X_i$,
and~$Y_i$ the algebra
elements associated to the edges of this piece of graph after $i$ flips are
performed. From (\ref{abc}), (\ref{phi}),
these elements evolve as follows:
\be
\left(
\begin{array}{l}X_{i+1}\\Y_{i+1}\\A_{i+1}\\ B_{i+1}\\C_{i+1}\\D_{i+1}\\E_{i+1}
\end{array} \right)
\to
\left(
\begin{array}{l}Y_i - \phi^\hbar(-X_i)\\-X_i\\D_i\\ E_i\\A_i +
\phi^\hbar(X_i)\\B_i-\phi^\hbar(-X_i)\\C_i +
\phi^\hbar(X_i) \end{array} \right)
\ee
We must prove that this operator is periodic with period five.

Assume for a moment that this five-periodicity of $X_i$ has been established.
Then
five-periodicity of $Y_i$ and other variables follow from simple
calculations.

It suffices to prove this five-periodicity for the $X_i$ since the five-periodicity of the
other operators such as $Y_i$ then follows from elementary calculations.
As to the five-periodicity of $X_i$,
let us ``take logarithms'' and introduce four new algebra elements
$$
U _i= \e^{X_i}; \quad V _i=\e^{Y_{i}}; \quad \tilde{U}_i=\e^{X_i/\hbar}; \quad
\tilde{V}_i=\e^{Y_{i}/\hbar},
$$
which satisfy the following commutation relations
\bea
&{}&U_iV_i=q^{-2}V_iU_i,\qquad
\tilde{U}_i\tilde{V}_i=\tilde{q}^{-2}\tilde{V_i}\tilde{U_i},
\nonumber\\
&{}&U_i\tilde{V}_i=\tilde{V}_iU_i,
\qquad V_i\tilde{U}_i=\tilde{U}_iV_i,
\label{aaa}
\eea
where
$$
q=\e^{-\pi i \hbar},\qquad \tilde{q} =\e^{-\pi i/ \hbar}.
$$
Under the flip,
these variables are transformed in an especially simple way,
\bea
&{}&U_{i-1}= V_i^{-1}\label{e1}\\
&{}&V_{i-1}= U_i(1+qV_i)\label{e2}\\
&{}&\tilde{U}_{i-1}=\tilde{V}_i^{-1}\label{e3}\\
&{}&\tilde{V}_{i-1}=V_i(1+\tilde{q}\tilde{V}_i).\label{e4}
\eea
As the first step of the proof, we consider the inverse
transformation laws for $X_i$ and $Y_i$:
$$
X_{i-1}=-Y_i;\qquad Y_{i-1} = X_i+\phi^\hbar(Y_i).
$$
Equations (\ref{e1}) and (\ref{e3}) are obvious.
Using the standard formula
$$
\e^{A+F(B)} =  \e^{{{1\over{[A,B]}}~\int_B^{B+[A,B]}F(z)dz}}~\e^A,
$$
we obtain
\bea
V_{i-1}&{}&=\e^{Y_{i-1}}=\e^{X_i+\phi^\hbar(Y_i)}=e^{X_i}\e^{\int_Y^{Y+2\pi i
\hbar}\phi^\hbar(z)dz}
\nonumber\\
&{}&=U_i\exp\left(-\frac{\pi \hbar}{2}\int_\Omega\frac{e^{-ipz}(e^{-2p\pi
\hbar}-1)}{(-ip)\sinh(\pi p)\sinh(\pi \hbar p)}dp\right) = U_i\exp\left(\frac{\pi
\hbar}{i}\int_\Omega\frac{e^{-ip(z-\pi i \hbar)}}{p\sinh(\pi p)}dp\right)
\nonumber\\
&{}&=U_i(1+q^{-2}V_i).
\nonumber
\eea
The proof of (\ref{e4}) is analogous.

Now, in order to finally prove that $X_i$ is five-periodic, it suffices to verify
that both $U_i$ and $\tilde{U_i}$ are five-periodic. Indeed,
if only the operator $U_i$ is five-periodic, it does not
suffice because the logarithm of an operator is ambiguously
defined. However, if we have two families of operators $U$ and $\tilde{U}$,
which depend continuously on $\hbar$, then, assuming the
existence of an operator $X$
(depending continuously on $\hbar$) such that $U=e^X$ and
$\tilde{U}=e^{X/\hbar}$, then this operator is evidently unique.
(It can be found as
$\lim_{(m+n/\hbar)\rightarrow 0} (U^m\tilde{U}^n)/(m+n/\hbar)$ for any
irrational value of $\hbar$.) The five-periodicity of sequence (\ref{e2})
(and (\ref{e4})) is a direct calculation using (\ref{aaa}).

The only subtlety remaining is the possibility that some of the edges
$A,B,C,D,E$ coincide. (Note, however, that $X$ and $Y$ must have
exactly one common vertex.)
 If, say, edges $A$ and $C$ coincide, then the value
of the commutator $[A,X]$ is doubled by definition, and we can then fictitiously
{\it split\/} the edge $A=C$ into two half-edges with the matrices $X_{A/2}$
assigned to each half-edge using the formula
\be
\label{N.N.}
X_A=X_{A/2}~~\left(\begin{array}{cc}
                0 & 1 \\
                -1 & 0 \\
                \end{array}\right) ~~X_{A/2}.
\ee
The quantities $A/2=C/2$ have the same commutation
relations with the rest of variables as well as
the same transformation laws as the quantities $A$ and $C$ before,
so the earlier formulas remain valid if we simply replace
there $A'$ by $A'/2$ and $A$ by $A/2$.  The net effect is that commutators
with $A=C$ are doubled.

If edges $A$ and $B$ coincide, then we must correct formulas (\ref{th-Q}) using a
splitting as above but demanding $A'/2=A/2+X/2$. Obviously, $[A,X]=0$ in this case
(in which formulas for the quantum ordering below must be also corrected, see
Section 3.5).

The formulas in Proposition~\ref{propcase} can thus be realized for exponentiated
quantities, although in the current quantum case there will be corrections. For
instance, formulas (\ref{e2}) and (\ref{e3}) will be different; indeed, letting $A=C=X$,
we calculate that $e^{X/2}\mapsto e^{X/2+\phi^\hbar(Z)}= :e^{X/2}(1+e^Z):\equiv
e^{X/2}+e^{X/2+Z}=e^{X/2}(1+e^{-i\pi\hbar}e^Z)$ since the commutator of $X$ and $Z$
is doubled, where the normal ordering (the Weyl ordering) $:\cdot :$ is explained in
the next section. The transformation for $e^X$ itself becomes more
complicated when $A=C$, namely, $e^X=
(e^{X/2})^2\mapsto e^X(1+q^3e^Z)(1+qe^Z)$.~~~~\hfill{\it q.e.d.}

\vskip0.1in

\begin{cor}\label{cor1}

{\bf 1.} Let $K$ be an operator acting in the Hilbert space
$L^2(\RR)$ and having the integral kernel
\be
\label{corr1}
K(x,z) = F^\hbar(z)e^{-\frac{zx}{2\pi i \hbar}},
\ee
where
\be
\label{dlc}
F^\hbar(z) = \exp\left(-\frac{1}{4}\int_\Omega \frac{e^{-ipz}}{p\sinh(\pi
p)\sinh(\pi \hbar p)}dp\right)
\ee
Then the operator $K$ is unitary up to a multiplicative constant and
satisfies the identity
\be
\label{pic}
K^5=\mbox{\em const}.
\ee

\noindent {\bf 2.} \ Let $\hbar=m/n$ be a rational number and assume that
both $m$ and $n$ are odd. Introduce a linear operator $L(u)$
acting in the space
${\Bbb C}^n$ and depending on one positive real parameter $u$
through its matrix
\be
\label{dld} L(u)^i_j = F^\hbar(j,u)q^{-4ij},
\ee
where
$$
F^\hbar(j,u) =
(1+u)^{j/n}\prod_{k=0}^{j-1}(1+q^{-4k+2}u^{1/n})^{-1}.
$$
Then the following
identity holds\/{\rm:}
\be
\label{pid}
L(u)L(v+uv)L(v + v u^{-1}+ u^{-1})L(u^{-1}v^{-1}+ u^{-1})L(v^{-1}) = 1.
\ee
\end{cor}
\vspace{3mm}

Using this construction, Kashaev~\cite{Kashaev3} constructed the
set of eigenfunctions of the quantum Dehn twist transformation.
Namely, Kashaev's dilogarithm function $e_b(z)$ is
$1/F^\hbar(2z)$ from (\ref{dlc}). We return to this
discussion when considering sets of quantum Dehn twists for
the torus in Section~5.

\subsection{Geodesic length operators}

We next embed the algebra of geodesics
(\ref{geod}) into a suitable
completion of the constructed algebra ${\cal T}^\hbar$.
For any $\gamma$,
the geodesic function $G_\gamma$ (\ref{geod}) can be expressed  in terms of
shear coordinates on ${\cal T}_H$:
\begin{equation}
G_\gamma \equiv \tr P_{Z_1\cdots Z_n}=
\sum_{j\in J}\exp\left\{{\frac{1}{2} \sum_{\alpha \in E(\Gamma)}
m_j(\gamma,\alpha) z_\alpha}\right\},  \label{clen}
\end{equation}
where $m_j(\gamma,\alpha)$ are
integers and $J$ is a finite set of indices.
In order to find the
quantum analogues of these functions, we
denote by $\widehat{\cal T}^\hbar$ a completion of the algebra ${\cal
T}^\hbar$ containing $e^{xZ_\alpha}$ for any real $x$.

For any closed path $\gamma$ on $F$, define the {\it quantum geodesic}
operator $G^\hbar_\gamma \in
\widehat{\cal T}^\hbar$ to be
\begin{equation} \label{qlen}
G^\hbar_\gamma \equiv
\ORD{\tr P_{Z_1\dots Z_n}}\equiv
\sum_{{j\in J\atop \kappa\in\{j\}}}
\exp\left\{{\frac{1}{2} \sum_{\alpha \in E(\Gamma)}
\bigl(m_j(\gamma,\alpha) Z^\hbar_\alpha
+2\pi i\hbar
c_j^\kappa(\gamma,\alpha)
\bigr)}\right\},
\end{equation}
where the {\it quantum ordering} $\ORD{\cdot}$ implies that we vary the classical
expression (\ref{clen}) by introducing additional
integer coefficients  $c_j^\kappa(\gamma,\alpha)$, which must
be determined from the conditions below.
Notice that the operators $\{G^\frac{1}{\hbar}_\gamma\}$ themselves can be considered as
belonging to the algebra $\widehat{\cal T}^\hbar$ insofar as
\begin{equation} \label{dlen}
G^\frac{1}{\hbar}_\gamma = \sum_{j\in J}
\exp\left\{{\frac{1}{2\hbar} \sum_{\alpha \in
E(\Gamma)} \bigl(m_j(\gamma,\alpha) Z^\hbar_\alpha
+2\pi i c_j^\kappa(\gamma,\alpha) \bigr)}\right\}.
\end{equation}

In what follows for the notational simplicity,
we shall sometimes omit the superscript $\hbar$ from $G^\hbar$
and write it merely $G$ assuming that $G$ is either an operator or a classical
geodesic function depending on the context. In this paper, we
shall concentrate on the case of quantum functions of the $G^\hbar$-type; the consideration
of the sector $G^{1/\hbar}$ is analogous and does not lead to new effects
(at least at the present stage of understanding),
so we omit it.  We also call a quantum geodesic function merely a
quantum geodesic (since implicitly quantum objects admit only functional, not
geometrical, descriptions).

We wish to associate an operatorial {\it quantum multicurve} QMC to a multiset
$\hat C$ of quantum geodesics corresponding to disjointly embedded families of nonnegative
integrally weighted geodesics.  One ansatz will be that operators corresponding to
disjoint underlying geodesics must commute; this implies that the ordering
in which the quantum geodesics enter the product QMC is immaterial, where the product
is defined as for GMCs.  We next formulate the defining properties of quantum
geodesics.

\leftskip .4in

\vskip .1in\noindent {\bf 1.} If closed paths $\gamma$ and $\gamma^\prime$ do not
intersect, then the operators $G^\hbar_\gamma$ and
$G^\hbar_{\gamma^\prime}$ commute.

\vskip .1in\noindent {\bf2.} {\it Naturality.} { The mapping class group}
$MC(F)$ (\ref{abc}) acts naturally, i.e., for any $\{G^\hbar_\gamma\}$, $\delta \in MC(F)$
and closed path $\gamma$ in a spine $\Gamma$ of $F$, we have
$\delta(G^\hbar_\gamma) = G^\hbar_{\delta \gamma}$.

\vskip .1in\noindent {\bf3.} {\it Geodesic algebra}. The product of two quantum geodesics
is a linear combination of QMC's governed by the (quantum) skein
relation \cite{Turaev}.

\vskip .1in\noindent {\bf4.} {\it Orientation invariance.} Quantum traces of
direct and inverse geodesic operators coincide.

\vskip .1in\noindent {\bf5.} {\it Exponents of geodesics.}
A quantum geodesic $G_{n\gamma}$ corresponding to the $n$-fold concatenation of $\gamma$
is expressed via $G_{\gamma}$ exactly as in the classical case, namely,
\be
\label{cheb}
G_{n\gamma}=2T_n\bigl(G_{\gamma}/2\bigr),
\ee
where $T_n(x)$ are Chebyshev's polynomials.

\vskip .1in\noindent {\bf6.} {\it Duality.} For any $\gamma$ and $\gamma^\prime$,
the operators $G^\hbar_\gamma$ and  $G^\frac{1}{\hbar}_{\gamma^\prime}$ commute.

\leftskip=0ex\vskip .1in

We shall let  the standard normal ordering symbol ${:}\e^{a_1}\e^{a_2}\cdots\e^{a_n}{:}$ denote the {\it Weyl
ordering} $\e^{a_1+\cdots+a_n}$,  i.e.,
\be
\label{normal}
{:}\e^{a_1}\e^{a_2}\cdots\e^{a_n}{:}=1+(a_1+\cdots +a_n)+{1\over {2!}}~(a_1+\cdots +a_n)(a_1+\cdots
+a_n)+\cdots
\ee
for any set of exponents with $a_i\neq -a_j$ for $i\neq j$,
In particular, the Weyl ordering implies total symmetrization in the susbscripts.

\vskip .1in

\begin{proposition}\label{lem31} 
For any graph simple geodesic $G$ with respect to any spine $\Gamma$,
the coefficients $c_j^\kappa(\gamma,\alpha)$ in {\rm(\ref{qlen})} are identically zero.
\end{proposition}

\vskip .2in

\noindent {\sl Proof}~Consider term-by-term the trace of the matrix product for the
quantum graph simple geodesic and expand it in Laurent monomials in
$\e^{Z_i/2}$. It is easy to see that each term $\e^{Z_i/2}$ comes either in power $+1$, or
$-1$ in the corresponding monomial and there are no equivalent monomials in
the sum.  This means that in order to have a {Hermitian} operator, we must apply the
Weyl ordering with no additional $q$-factors (by the correspondence principle, each
such factor must be again a Laurent monomial in $q$ standing by the corresponding
term, which breaks the self-adjointness unless all such monomials are unity).  Since quantum
Whitehead moves must preserve the property of being Hermitian, if a graph-simple geodesic transforms
to another graph-simple geodesic, then a Weyl-ordered expression transforms to a Weyl-ordered
expression, and only these expressions are self-adjoint.~~~~\hfill{\it q.e.d.}

\vskip .2in

\noindent {\bf Torus Example~}  For the torus
with one hole, there are three
graph simple quantum geodesics for any spine, which are exactly (\ref{torus-classic})
in the Weyl-ordered form.
The quantum geodesics $\wtd G_Z$ obtained
from $G_Z$ by the flip transformation is
\be
\label{tildeG} {\wtd
G_Z}=\e^{-X/2-Y/2-Z}+\e^{X/2-Y/2-Z}+\e^{X/2-Y/2}\cdot 2\cos(\pi\hbar)
+\e^{X/2-Y/2+Z}+\e^{X/2+Y/2+Z}.
\ee

The product of two graph simple quantum geodesics is
\be
\label{TurPR}
G_XG_Y=\e^{i\pi\hbar/2}{\wtd G}_Z+\e^{-i\pi\hbar/2}G_Z.
\ee

Denoting $q\equiv\e^{-i\pi\hbar}$, \
$[A,B]_q\equiv q^{1/2}AB-q^{-1/2}BA$, and~$\xi= q-q^{-1}$,
we obtain from (\ref{TurPR})
\be
[G_X,G_Y]_q=\xi G_Z,\quad
[G_Y,G_Z]_q=\xi G_X,\quad
[G_Z,G_X]_q=\xi G_Y.
\label{so3}
\ee
This algebra is exactly the~$so_q(3)$ quantum algebra
studied in~\cite{Klim}. There is a unique central element, the quantum Markov relation
\be
{\cal M}=G_XG_YG_Z-q^{1/2}(G_X^2+q^{-2}G_Y^2+G_Z^2).
\ee

\clearpage

\subsection{Algebra of quantum geodesics}

Let $G_1$ and $G_2$ correspond to
the respective graph simple geodesics with respect to the same spine having one nontrivial
intersection.
For $G_1$ and $G_2$, formula (\ref{qlen})
implies, by virtue of Proposition~\ref{lem31}, the mere Weyl ordering.

After some algebra, we obtain (cf. (\ref{Goldman}))
\be
\label{AG4}
G_1G_2=\e^{-i\pi\hbar/2}G_Z+\e^{i\pi\hbar/2}\wtd G_Z,
\ee
where $G_Z$
coincides with the Weyl-ordered $G_{\mbox{\tiny I}}$
in the classical case (cf.~(\ref{Goldman})) while $\wtd G_Z$
contains the quantum correction term
\bea
\label{AG6}
\wtd G_Z&=&\ORD{\tr^1\tr^2
\dots (e^1_{ij}\otimes e^2_{ji}) [X^1_Z\otimes X^2_Z]\dots}
\nonumber\\
&=&{:}\tr^1\tr^2
\dots (e^1_{ij}\otimes e^2_{ji}) \bigl[X^1_Z\otimes X^2_Z
+2(1-\cos\pi\hbar)e^1_{11}\otimes e^2_{22}\bigr]\dots\,{:}\ .\nonumber
\eea
Here $e^1_{ij}\otimes e^2_{ji}$ is the standard $r$-matrix that permutes
the spaces ``$1$'' and ``$2$,'' and as a result, the ``skein''
relation of form (\ref{Goldman}) appears. Locally, this relation
has exactly the form proposed by Turaev~\cite{Turaev}, i.e., for
two graph simple geodesics intersecting at a single point, we have the defining relation
\be
\setlength{\unitlength}{.8mm}%
\begin{picture}(90,40)(0,45)
\thicklines
\put(-10,50){\line(1, 1){9}}
\put( 1, 61){\line(1, 1){9}}
\put(-10,70){\line(1,-1){20}}
\put(18,60){\makebox(0,0){$=$}}
\put(30,61){\makebox(0,0){$\e^{-i\pi\hbar/2}$}}
\put(40,60){\oval(20,20)[r]}
\put(63,60){\oval(20,20)[l]}
\put(70,60){\makebox(0,0){$+$}}
\put(78,61){\makebox(0,0){$\e^{i\pi\hbar/2}$}}
\put(95,48){\oval(20,20)[t]}
\put(95,72){\oval(20,20)[b]}
\put(-11,68){\makebox(0,0)[rb]{\Large$G_1$}}
\put(-11,52){\makebox(0,0)[rt]{\Large$G_2$}}
\put(52,74){\makebox(0,0)[cb]{\Large$G_Z$}}
\put(95,74){\makebox(0,0)[cb]{\Large$\wtd G_Z$}}
\end{picture}
\label{skein}
\ee
(The order of crossing lines corresponding to
$G_1$ and $G_2$ depends on which quantum geodesic
occupies the first place in the product;
the rest of the graph remains unchanged for all items in (\ref{skein})).
Note, however, that if the quantum geodesics
$G_1$ and $G_2$ correspond to graph simple geodesics, we may turn the geodesic
$\wtd G_Z$ again into the graph simple geodesic $\wtd G'_Z$
by performing the quantum flip with respect to the edge $Z$.

If we now compare two unambiguously determined expressions:
$\wtd G'_Z$, which must be Weyl ordered, and $\wtd G_Z$ obtained from
the geodesic algebra, we find
that $\wtd G_Z=\wtd G'_Z$. This enables us to formulate the main assertion.

\begin{lemma}\label{lem34}{\rm \cite{ChF}}
There exists
a unique quantum ordering $\ORD{\dots}$ {\rm (\ref{qlen})}, which is
generated by the quantum geodesic
algebra (\ref{skein}) and consistent with the quantum mapping
class groupoid transformations {\rm(\ref{abc})}, i.e., so that the quantum
geodesic algebra is invariant under the action of the quantum mapping
class groupoid.
\end{lemma}

We can now relax the constraints of graph simplicity
of curves: as the quantum geodesic algebra is quantum
mapping class group invariant, having two {\it arbitrary\/} embedded geodesics
with a single intersection, we can transform them using quantum morphisms
to a canonical form of
graph simple geodesics and employ the Weyl order. Relation (\ref{skein})
remains valid in both cases.

Let us now address the problem of multiple intersections. Here, we have the
following lemma.

\begin{lemma}\label{lem35}{\rm \cite{ChF}}
If {\rm more\/} than one intersection of two QMCs
occurs, then the quantum skein
relations {\rm(\ref{skein})}
must be applied {\rm simultaneously\/} at {\rm all\/}
intersection points.
\end{lemma}

This lemma
implies the standard Reidemeister moves for curves on a graph, where
the empty loop gives rise to a factor $-\e^{-i\pi\hbar}-\e^{i\pi\hbar}$; that is,
for geodesics intersecting generically, apply (\ref{skein})
simultaneously at all intersection points to obtain the Reidemeister moves.

\begin{remark}\label{exam12}
{\rm
The quantum algebra ${\cal M}_1^2$ was studied in~\cite{ChF}, where the exact correspondence
with Kauffman bracket skein
quantization of the corresponding Poisson algebra of geodesics (see~\cite{Bull})
was observed.
In~\cite{Bondal}, this algebra arises as the
quantum deformation algebra of the classical groupoid Poisson relations
for the group $SL(4)$. This algebra is however only one
among many quantum Nelson--Regge algebras corresponding
to Riemann surfaces of higher genera which are described in the next section.
}
\end{remark}

\subsection{Quantizing the Nelson--Regge algebras}
The algebra (\ref{P-geod}) was quantized by the deformation
quantization method in~\cite{NR,NRZ}. We now explicitly implement
the quantization conditions (\ref{comm}). It is convenient
to represent the (classical or quantum) elements ${\cal G}_{ij}$ as chords connecting the
points of the cyclically ordered set of indices
$i,j$.  There is then a trichotomy: if two chords do not intersect, then the
corresponding geodesics do not intersect either, and the quantum
geodesics commute (Figure~4a); if two chords have a common endpoint,
then the corresponding geodesics intersect at one point, and the
three quantum geodesics ${\cal G}^\hbar_{ij}$, ${\cal G}^\hbar_{jk}$,
${\cal G}^\hbar_{ki}$ (as
depicted in Figure~4b) constitute the quantum subalgebra $so_q(3)$;
if two chords intersect at an interior point (as depicted in Figure~4c), then the
corresponding geodesics intersect in two points, and the
corresponding quantum geodesics ${\cal G}^\hbar_{ij}$ and ${\cal G}^\hbar_{kl}$, $i<k<j<l$,
satisfy the commutation relation
\be
[{\cal G}^\hbar_{ij},{\cal G}^\hbar_{kl}]=\xi({\cal G}^\hbar_{ik}{\cal G}^\hbar_{jl}
-{\cal G}^\hbar_{il}{\cal G}^\hbar_{jk})
\ee
with the usual commutator (not the $q$-commutator) and where again $\xi= q-q^{-1}$.

\be
\setlength{\unitlength}{.6mm}%
\begin{picture}(100,50)(10,30)
\thicklines
\put(0,60){\oval(30,30)} \put(60,60){\oval(30,30)}
\put(120,60){\oval(30,30)} \put(-15,55){\circle*{3}}
\put(-5,45){\circle*{3}} \put(5,75){\circle*{3}}
\put(15,65){\circle*{3}} \put(-15,55){\line(1, 1){20}}
\put(-5,45){\line(1, 1){20}} \put(-20,55){\makebox(0,0){$i$}}
\put(-7,40){\makebox(0,0){$l$}} \put(6,80){\makebox(0,0){$j$}}
\put(19,68){\makebox(0,0){$k$}} \put(0,37){\makebox(0,0){a}}
\put(47.5,50){\circle*{3}} \put(72.5,50){\circle*{3}}
\put(60,75){\circle*{3}} \put(47.5,50){\line(1, 2){12.5}}
\put(72.5,50){\line(-1, 2){12.5}} \put(47.5,50){\line(1,0){5}}
\put(72.5,50){\line(-1,0){5}} \put(55,50){\line(1,0){10}}
\put(45,47){\makebox(0,0){$i$}} \put(77,47){\makebox(0,0){$k$}}
\put(60,80){\makebox(0,0){$j$}} \put(60,37){\makebox(0,0){b}}
\put(108,48){\circle*{3}} \put(132,48){\circle*{3}}
\put(108,72){\circle*{3}} \put(132,72){\circle*{3}}
\put(108,48){\line(1, 1){24}} \put(108,48){\line(1, 1){24}}
\put(132,48){\line(-1, 1){11}}
\put(119,61){\line(-1,1){11}}\put(108,48){\line(1,0){6}}
\put(117,48){\line(1,0){6}} \put(126,48){\line(1,0){6}}
\put(108,72){\line(1,0){6}} \put(117,72){\line(1,0){6}}
\put(126,72){\line(1,0){6}} \put(108,48){\line(0,1){6}}
\put(108,57){\line(0,1){6}} \put(108,66){\line(0,1){6}}
\put(132,48){\line(0,1){6}} \put(132,57){\line(0,1){6}}
\put(132,66){\line(0,1){6}} \put(108,55.5){\circle*{1}}
\put(108,64.5){\circle*{1}} \put(132,55.5){\circle*{1}}
\put(132,64.5){\circle*{1}} \put(106,41){\makebox(0,0){$l$}}
\put(106,79){\makebox(0,0){$i$}} \put(137,43){\makebox(0,0){$j$}}
\put(135,79){\makebox(0,0){$k$}} \put(120,37){\makebox(0,0){c}}
\end{picture}
\label{circles} \ee

\centerline{\bf Figure~4-picture of Nelson--Regge quantum relations}

\hspace{10pt}

\subsection{Improving the quantum ordering}

We now extend the construction of QMCs by taking
the products of operatorial matrices $X_Z$, $L$, and $R$ along oriented geodesics as before,
but we do not apply the trace operation as in the definition of the geodesic operators.

For the three cases of oriented curves depicted in Figure~5 below, we may apply the
indicated quantum Whitehead move and calculate as follows:

\vskip .1in

\noindent
For curve 1, we obtain
\bea
X_{B'}LX_{-Z}LX_{A'}&=&
\left(\begin{array}{cc} \e^{-B'/2}\e^{-Z/2}\e^{A'/2}+\e^{-B'/2}\e^{Z/2}\e^{A'/2}
                         & -\e^{-B'/2}\e^{Z/2}\e^{-A'/2}\\
                \e^{B'/2}\e^{-Z/2}\e^{A'/2} & 0\end{array}\right)
\nonumber\\
&=&\left(\begin{array}{cc} \e^{A/2-B/2}
                         & -\e^{-A/2-B/2-\frac14 [A,B]}\\
                \e^{A/2+B/2-\frac14 [A,B]} & 0\end{array}\right)\ \ \hbox{for\ $B'\ne A'$}
\nonumber\\
&=&\e^{-\frac18[A,B]}\left(\begin{array}{cc} \e^{-B/2}\e^{A/2}
                         & -\e^{-B/2}\e^{-A/2}\\
                \e^{B/2}\e^{A/2} & 0\end{array}\right)\ \ \hbox{for\ $B'\ne A'$}
                \nonumber\\
                &=&\left\{\begin{array}{ll}
                  \e^{-\frac18[A,B]}X_BLX_A&\hbox{for\ $A\ne B$},
                 \\
                  X_BLX_A&\hbox{for\ $A=B$}. \\
                \end{array}
               \right.
\label{1-1}
\eea

\noindent For curve 2, we obtain
\bea
X_{C'}RX_{-Z}LX_{A'}&=&
\left(\begin{array}{cc} \e^{-C'/2}\e^{Z/2}\e^{A'/2}
                         & -\e^{-C'/2}\e^{Z/2}\e^{-A'/2}\\
                \e^{C'/2}\e^{Z/2}\e^{A'/2}+\e^{C'/2}\e^{-Z/2}\e^{A'/2}
                 & -\e^{C'/2}\e^{Z/2}\e^{-A'/2}\end{array}\right)
\nonumber\\
&=&\left(\begin{array}{cc} \e^{A/2-C/2+Z/2+\frac14 [A,B]}
                         & -\e^{-A/2-C/2+Z/2}-\e^{-A/2-C/2-Z/2}\\
                \e^{A/2+C/2+Z/2} & -\e^{-A/2+C/2+Z/2-\frac14 [A,B]}\end{array}\right)
\nonumber\\
&=&\left(\begin{array}{cc} \e^{-C/2}\e^{Z/2}\e^{A/2}
                         & -\e^{-C/2}\e^{Z/2}\e^{-A/2}-\e^{-C/2}\e^{-Z/2}\e^{-A/2}\\
                \e^{C/2}\e^{Z/2}\e^{A/2}
                 & -\e^{C/2}\e^{Z/2}\e^{-A/2}\end{array}\right)
\nonumber\\
                 &=&X_CLX_ZRX_A,
\label{1-2}
\eea
and if edges $A$ and $C$ coincide, we merely use the same splitting as in (\ref{N.N.});
no additional factors arise.

\noindent For curve 3, we obtain
\bea
X_{C'}RX_{-Z}RX_{D'}&=&
\left(\begin{array}{cc} 0 & -\e^{-C'/2}\e^{Z/2}\e^{-D'/2}\\
                \e^{C'/2}\e^{-Z/2}\e^{D'/2} &
                -\e^{C'/2}\e^{Z/2}\e^{-D'/2}-\e^{-C'/2}\e^{-Z/2}\e^{-D'/2}\end{array}\right)
\nonumber\\
&=&\left(\begin{array}{cc} 0
                         & -\e^{-C/2-D/2+\frac14 [A,B]}\\
                \e^{C/2+D/2+\frac14 [A,B]} & -\e^{C/2-D/2}
                \end{array}\right)\ \ \hbox{for\ $C'\ne D'$}
\nonumber\\
&=&\e^{\frac18[A,B]}\left(\begin{array}{cc} 0
                         & -\e^{-C/2}\e^{-D/2}\\
                \e^{C/2}\e^{D/2} & -\e^{C/2}\e^{-D/2}\end{array}
                \right)\ \ \hbox{for\ $C'\ne D'$}
                \nonumber\\
                &=&\left\{\begin{array}{ll}
                \e^{\frac18[A,B]}X_CRX_D.&\hbox{for\ $C\ne D$},
                 \\
                  X_CRX_D&\hbox{for\ $C=D$}. \\
                \end{array}
               \right.
\label{1-3}
\eea

In these formulas, we have used the identities
\bea
\e^{A'/2+C'/2+Z/2}+\e^{A'/2+C'/2-Z/2}&=&:\e^{A/2+C/2}\frac{1}{1+\e^{-Z}}(\e^{Z/2}+\e^{-Z/2}):~~\nonumber\\
&=&\e^{A/2+C/2+Z/2}
\label{2-1}
\eea
and
\bea
\e^{-A'/2-C'/2+Z/2}&=&:\e^{-A/2-C/2}(1+\e^{-Z})\e^{Z/2}:\nonumber\\
&=&\e^{-A/2-C/2+Z/2}+\e^{-A/2-C/2-Z/2},
\label{2-2}
\eea
where vertical dots denote the Weyl ordering
(\ref{normal}) as before implying total symmetrization
with respect to all the variables $\{A,B,C,D,Z\}$.

As we have just observed and amazingly enough, the only thing that changes under a
quantum Whitehead move is the overall factor standing by the product of matrices,
and even this factor can easily be taken into account if included with each
left-turn matrix $L$ is an overall factor $\e^{-\frac18[A,B]}$, and included with
each right-turn matrix $R$ is an overall factor $\e^{\frac18[A,B]}$.

The subtlety of potentially multiple intersections between curves corresponds to the possibility that edges of
the graph may coincide.  As in the previously useful convenient fiction of splitting such an edge into two
half-edges, the formulas (\ref{1-1}), (\ref{1-2}), and (\ref{1-3}) remain valid if we replace there $A'/2$ by
$A'/4$ and $A/2$ by $A/4$ in case $A=C$ for instance.

We thus have the quantum analogue of Lemma~\ref{lem1}.

\begin{lemma}\label{qlem1}
For any oriented non boundary-parallel geodesic $\gamma$, taking the finite (periodic) sequence of
matrices with quantum entries in $X_{Z_i}$ as before but making the replacement
\be
\label{LRtilde}
\tilde L= q^{-1/4}L, \qquad \tilde R= q^{1/4}R,
\ee
where $q=\e^{-i\pi\hbar}$,
the resulting product of matrices
\be
\label{Pcyclic}
P_\gamma =\bigl(X_{Z_n}{\tilde L}X_{Z_{n-1}}{\tilde R}\dots
X_{Z_3}{\tilde L}X_{Z_2}{\tilde R}X_{Z_1}{\tilde R}\bigr)
\ee
is invariant under quantum
Whitehead moves.
\end{lemma}

\begin{remark}\label{simple}
{\rm In the case of a boundary-parallel curve, where a sub-word of the form $X_{D'}{\tilde
R}X_{-Z}{\tilde R}$ transforms to a sub-word of the form
$X_D{\tilde {\tilde R}}$,
we must set the resulting turn matrices to be
$$
{\tilde {\tilde L}}=q^{-1/2}L, \qquad {\tilde {\tilde R}}=q^{1/2}R,
$$
i.e., they have quantum factors doubled in comparison with ${\tilde L}$ and ${\tilde R}${\rm}.
}
\end{remark}

We now address the
question of the proper quantum ordering.  Using MCG
transformations, we can reduce {any} simple curve either to the form of a
graph simple curve with exactly one left and one right turn (c.f. formulas
(\ref{torus-classic}) and Figure~5) if this curve is {not} boundary-parallel, or to the form
$X_D{\tilde {\tilde R}}$ or $X_D{\tilde {\tilde L}}$ for a boundary-parallel curve.
In both cases, each term in the corresponding Laurent polynomial must be
self-adjoint, which immediately results in the Weyl ordering by
Proposition~\ref{lem31}.

\centerline{
\epsffile{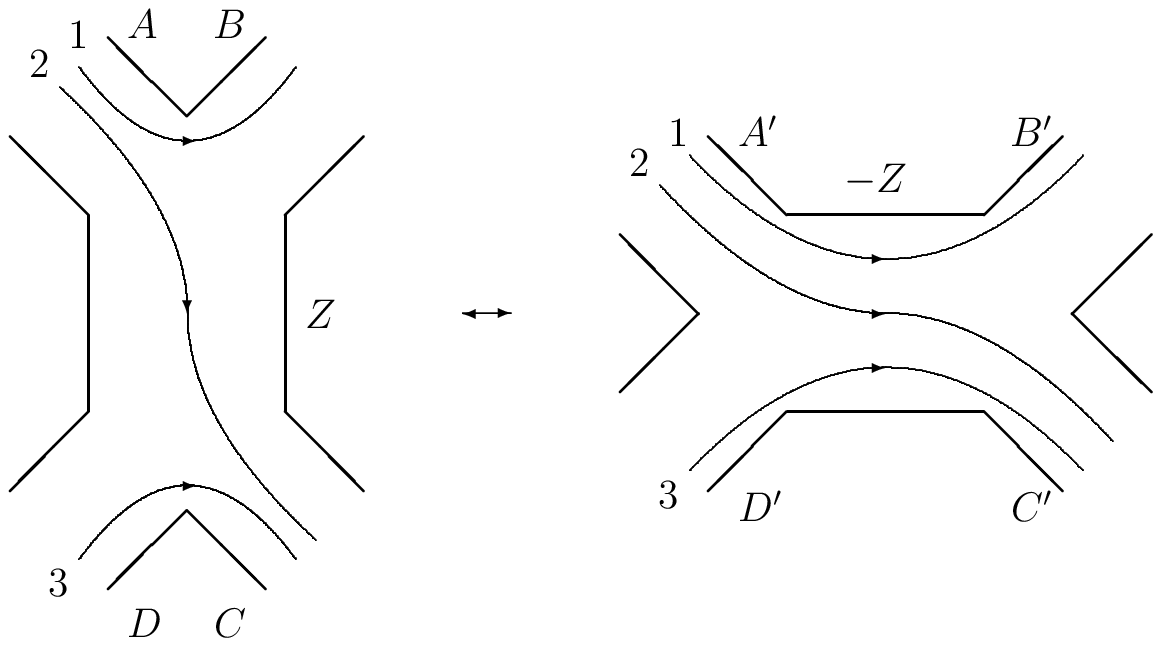}}


\centerline{\bf Figure~5-three cases of flips for geodesics}

\vspace{10pt}

This observation
does not suffice to derive the proper quantum ordering in our quantization of Thurston theory since we cannot consider in
a consistent way an infinite product of matrices corresponding to infinite leaves of a measured foliation.  Our tools for
analyzing the Thurston theory in this case will devolve to naturality of the MCG action on the QMC
algebra and an operatorial version of infinite continued fractions in Section~\ref{cfe2sec}.  On the other hand, the
improved quantum ordering is used in our analysis of closed geodesics on the torus, i.e., of operatorial
finite continued fraction expansions in Section~\ref{cfesec}.

\clearpage

\vskip .2in

\newsection  {Classical Thurston Theory of Surfaces}

\vskip .1in

Let $F_g^s$ denote an oriented smooth surface with $s\geq 0$ punctures
(so $F_g^s$ may be closed without boundary in this section), with genus $g\geq 0$, and
with negative Euler characteristic $2-2g-s< 0$.

\vskip .1in

\subsection{Measured foliations and Thurston's boundary}

\vskip .1in

Define a {\it measured foliation} on $F=F_g^s$ to be a
one-dimensional topological foliation ${\cal F}$ of $F$, where in a
neighborhood of any $p\in F$, ${\cal F}$ must restrict (in an appropriate chart) to
the horizontal foliation as in Figure 6a or it must restrict
to a foliation with one $n$-pronged singularity at $p$, for $n\geq 3$, as
illustrated in Figure 6b.

\centerline{
\epsffile{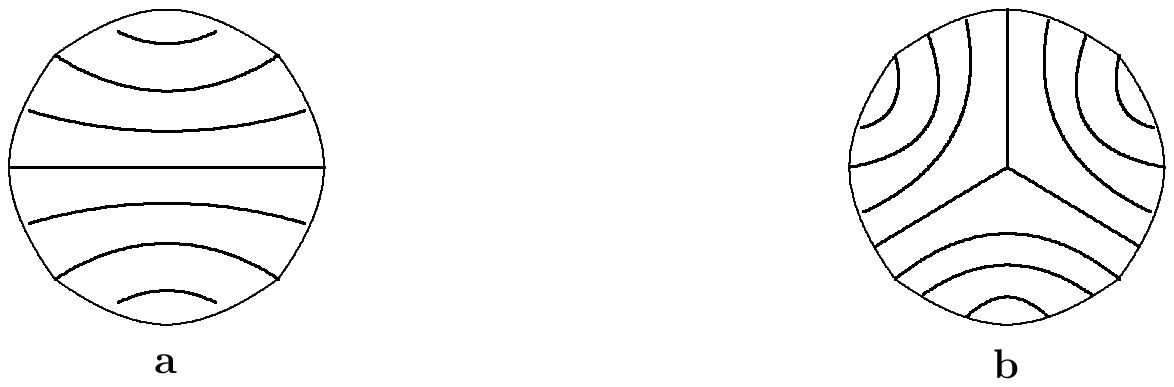}}

\centerline{\bf Figure 6-pictures of foliations}

\vskip .2in
\noindent Furthermore, ${\cal F}$ comes equipped with a
{\it transverse measure} $\mu$, which assigns to any arc $a$ in $F$
that is transverse to ${\cal F}$ a
real number $\mu (a)\in{\Bbb R}_{\geq 0}$, where $\mu$ is required to
satisfy:

\vskip .1in

\leftskip .3in

\noindent {\it No Holonomy}:~~If $a_0,a_1$ are homotopic through
arcs $a_t$ transverse to ${\cal F}$, for $0\leq t\leq 1$, keeping the
endpoints of $a_t$ on the same leaf for all $t$, then $\mu (a_0)=\mu
(a_1)$;

\vskip .1in

\noindent {\it $\sigma$-Additivity}:~If $a$ is the serial concatenation
of transverse arcs $a_1,a_2,\ldots $, then $\mu (a)=\sum _{i\geq 1}\mu (a_i)$.

\vskip .1in

\leftskip=0ex

In other words, in the neighborhood of a non-singular point, there is a local chart $\phi:U\to {\Bbb R}^2=\{
(x,y):x,y\in{\Bbb R}\}$ so that $\phi ^{-1}(\{ y=~{\rm constant}\})$ are the leaves of the foliation
${\cal F}\cap U$.  If two charts $U_i$ and $U_j$ intersect, then the transition function $\phi _{ij}$ are
of the form $\phi _{ij}(x,y)=(h_{ij}(x,y), c_{ij}\pm y)$, where $c_{ij}$ is constant.  In these coordinates,
the transverse measure is $|dy|$.  In case the transition functions can be chosen with constant sign
$\phi_{ij}(x,y)=(h_{ij}(x,y),c+y)$, i.e., if the foliation is ``transversely orientable'', then (away
from the singular points) $y$ is the primitive of a closed one-form on $F$.

Another canonical construction of a measured foliation on a
Riemann surface is given by taking the leaves of the foliation to be the level sets of
a harmonic function, where the transverse measure is given by integrating
the conjugate differential along transverse arcs.  Still another example is given by
the homotopy class of a finite collection of disjointly embedded (weighted) curves,
as we shall see.

\centerline{
\epsffile{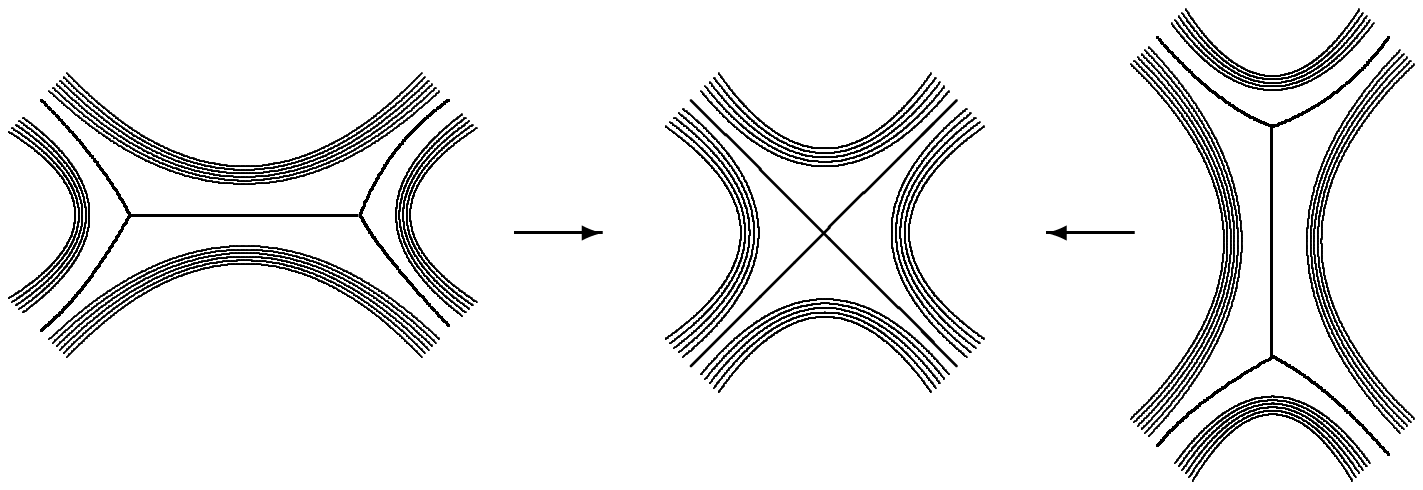}}

\centerline{\bf Figure 7-Whitehead collapses on foliations}

\vskip .2in

There is an equivalence relation on measured foliations generated by
isotopy and {\it Whitehead collapse} as illustrated in Figure~7, and the set of all
equivalence classes (including the empty measured foliation $\varnothing$) is denoted ${\cal
MF}={\cal MF}(F)$, where the class of $({\cal F},\mu )$ is denoted $[{\cal F},\mu ]$.
To naturally topologize
${\cal MF}$, we introduce the discrete set ${\cal S}={\cal S}(F)$
consisting of all free homotopy classes
$[c]$ of simple closed curves
$c$ in
$F$ which are neither null homotopic nor puncture-parallel.
We shall also require the set ${\cal S}'={\cal S}'(F)$ consisting of homotopy
classes of all disjointly embedded families of curves in $F$, where each component of
the family lies in ${\cal S}$ and no two components are homotopic.

One can
show \cite{FLP} that for any
$[c]\in {\cal S}$, there is a representative $c_{\cal F}\in[c]$ which
minimizes the $\mu$-transverse measure in its homotopy class; furthermore,
given  equivalent measured foliations $({\cal F}_i,\mu _i)$, for $i=1,2$,
we have $\mu _1(c_{{\cal F}_1})=\mu _2(c_{{\cal F}_2})$ \cite{FLP}, and hence there is a
well-defined mapping
\bea
\nonumber
J:{\cal MF}&{}&\to~~~~{\Bbb R}_{\geq 0}^{\cal S}\\
\nonumber({\cal F},\mu )&{}&\mapsto ~~~\biggl (i_{({\cal F},\mu)}:[c]\mapsto \mu (c_{\cal F})
\biggr ),
\eea
where the empty measured foliation $\varnothing$ is identified with $\vec 0=\{ 0\}
^{\cal S}$. This mapping $J$ is an injection \cite{FLP} and induces a topology on
${\cal MF}$
(where a neighborhood of $\varnothing$ is homeomorphic to a cone from $\varnothing$ over ${\cal MF}/{\Bbb R}_{>0}$,
with the natural action by homothety of ${\Bbb R}_{>0}$ on measures).  The
function
$i_{({\cal F},\mu)}$ is called the {\it (geometric) intersection function} of the measured foliation ${({\cal F},\mu)}$.

Given $[c]\in{\cal S}$, there is a corresponding measured foliation defined as
follows.
Choose a representative $c$ of $[c]$ and collapse $F-c$ onto a spine (making further choices)
to build a foliation
${\cal F}_c$ of $F$, whose leaves either lie in the spine or are homotopic to $c$, and choose a transverse measure
$\mu _c$ on
${\cal F}_c$ that pulls back under the collapsing map the counting
measure $\delta _c$ on $c$.  It is a classical fact due to Whitehead \cite{FLP}
that the resulting Whitehead
equivalence class
$[{\cal F}_c,\mu _c]$ is well-defined independent of any choices.
Taking the projective class of this foliation, we may thus
regard
$$
{\cal S}\subseteq {\cal PF}_0.
$$
Furthermore by construction, $i_{[{\cal F}_c,\mu _c]}([d])$, for $[d]\in {\cal S}$, is just the geometric
intersection number of $[c]$ and $[d]$, i.e., the total number of intersections of representatives $c$ with $d$,
where $c$ and $d$ intersect minimally.  More generally, given a family of curves $c_1,\ldots ,c_n$ representing
a point of ${\cal S}'$, together with a collection $w_1,\ldots ,w_n\in{\Bbb R}_{>0}$ of ``weights'', we may again
collapse to a spine of $F-\cup\{ c_i\} _1^n$ to produce a foliation ${\cal F}$ and choose a measure $\mu$ on
${\cal F}$ that pulls back the weighted sum $\sum _1^n w_i\delta_{c_i}$ to get a well-defined Whitehead
equivalence class $[{\cal F},\mu ]$.  As mentioned before, one may thus associate a measured foliation to a
weighted curve family.

Of special interest is the subspace ${\cal MF}_0={\cal MF}_0(F)$
consisting of all measured foliations $({\cal F},\mu )$ of compact
support, i.e., any leaf of ${\cal F}$  with a transverse arc
$a$ so that $\mu (a)>0$ must be disjoint from a neighborhood of the
punctures and no such leaf is puncture parallel.  There is a natural ${\Bbb
R}_{>0}$-action on
${\cal MF}-\{\varnothing\}$ and
${\cal MF}_0-\{\varnothing\}$ given by scaling the transverse measure, and the corresponding quotients
$$
{\cal PF}={\cal PF}(F)=({\cal MF}(F)-\{\varnothing\})/{\Bbb R}_{>0},
$$
$$
{\cal PF}_0={\cal PF}_0(F)=({\cal MF}_0(F)-\{\varnothing\} )/{\Bbb R}_{>0}
$$
are the spaces of central interest in the sequel.  The projective class of
$({\cal F},\mu)$ and $[{\cal F},\mu ]$, respectively,  will be denoted
$({\cal F},\bar\mu)$ and $[{\cal F},\bar\mu ]$.

A point of ${\cal T}={\cal T}(F_g^s)={\cal T}_g^s$ may be regarded as the class of a
hyperbolic metric on $F$, i.e., a complete finite-area Riemannian metric on $F$ of
constant Gauss curvature -1.  We shall also require the ``Yamabe space'' ${\cal
Y}={\cal Y}(F)$ of all complete finite-area Riemannian metrics on $F$ of constant
Gauss curvature
$-x^2$, for some $x\in{\Bbb R}_{>0}$.  ${\cal Y}$ is canonically homeomorphic to
${\cal T}_g^s\times{\Bbb R}_{>0}$, where $(\rho ,x)$ corresponds to the class of the
metric
$x\rho$, and we let $\pi:{\cal Y}\to{\cal T}$ denote the projection onto the
first factor.  Define the map
\bea\nonumber
I:{\cal Y}{}\to{\Bbb R}_{\geq 0}^{\cal S}\\
\nonumber\rho {}\mapsto\ell_\rho (\cdot ),
\eea
where $\ell _\rho ([c])$ is the $\rho$-length of the unique $\rho$-geodesic in the
homotopy class $[c]$.  Thus, if $\rho\in{\cal Y}$ corresponds to $(\pi (\rho
),x)$, then $I(\rho )=I(\pi (\rho ),x)=x~I(\pi (\rho ),1)$.

The basic facts \cite{FLP} are that $I:{\cal Y}\to{\Bbb R}_{\geq
0}^{\cal S}$ and $J:{\cal MF}_0(F)\to{\Bbb R}_{\geq
0}^{\cal S}$ are embeddings with disjoint images, and we may define a completion
$\overline{\cal Y}$ of ${\cal Y}$ in ${\Bbb R}_{\geq 0}^{\cal S}$ by setting
$$\overline{\cal Y}=I({\cal Y})\cup J({\cal MF}_0)$$ and identifying ${\cal Y}$ with
$I({\cal Y})$.  Passing to quotients under the homothetic actions of ${\Bbb R}_{>0}$
on ${\Bbb R}_{\geq 0}^{\cal S}-\{\vec 0\}$, on ${\cal Y}-\{\vec 0\}$, and on ${\cal
MF}_0-\{\varnothing\}$, we obtain Thurston's compactification
$$
\overline{\cal T}_g^s~~=~~\bigl ( {\cal T}_g^s\cup{\cal PF}_0\bigr )~~\approx ~~\bigl (\overline{\cal Y}-\{\vec
0\}\bigr )/{\Bbb R}_{> 0}
$$
of ${\cal T}\approx ({\cal Y}-\{\vec 0\} )/{\Bbb R}_{> 0}$ by ${\cal
PF}_0\approx ({\cal MF}_0-\{\varnothing\} )/{\Bbb R}_{> 0}$.

\begin{theorem}~\noindent {\bf 1.} {\rm \cite{Thurs},\cite{FLP},\cite{HarPen}} ${\cal
PF}_0(F_g^s)$ is naturally a piecewise linear sphere of dimension
$6g-7+2s$ which compactifies ${\cal T}_g^s$ to produce a closed ball
$\overline{\cal T}_g^s$.

\vskip .1in

\noindent{\bf 2.} {\rm \cite{Thurs},\cite{Penn0}} ~The action of the mapping
class group $MC_g^s$ on
${\cal T}_g^s$ extends continuously to an action on $\overline{\cal
T}_g^s$, where the action on ${\cal PF}_0(F_g^s)$ is the natural one, and there
are explicit piecewise linear formulas for the action of Dehn twist generators.

\vskip .1in

\noindent {\bf 3.} {\rm \cite{Thurs},\cite{FLP}} Suppose that a sequence
of hyperbolic metrics
$\rho _i$ on $F$ tends to a point $[{\cal F},\bar\mu ]\in{\cal PF}_0$.  In the
projectivization of ${\Bbb R}_{\geq 0}^{\cal S}-\{\vec 0\}$, the projectivized
length functions $\bar\ell_{\rho _i}$ of $\rho _i$ converge to the projectivized
intersection function $\bar i_{({\cal F},\bar\mu )}$ of $({\cal F},\bar\mu )$.

\vskip .1in

\noindent{\bf 4.} {\rm \cite{Thurs},\cite{Bon1}} The function $i_{({\cal F},\mu)}:{\cal
S}\to{\Bbb R}_{\geq 0}$ extends continuously to the ``geometric
intersection pairing'' ${\cal MF}_0\times{\cal MF}_0\to{\Bbb R}_{\geq 0}$, that vanishes on the diagonal and which
is also invariant under $MC_g^s$.

\vskip .1in

\noindent {\bf 5.} {\rm \cite{Thurs},\cite{PP2},\cite{Bon3}} The Weil--Petersson K\"ahler
two-form on
${\cal T}_g^s$ continuously extends (in the appropriate sense on Yamabe space) to a
non-degenerate symplectic form,
called ``Thurston's symplectic form'', on
${\cal MF}_0(F_g^s)$ ($\approx \overline{\cal Y}-{\cal Y}$), which is invariant under
$MC_g^s$.
\end{theorem}

Though there is all this beautiful natural structure on Thurston's boundary, the quotient ${\cal
PF}_0(F_g^s)/MC_g^s$ is maximally non-Hausdorff (i.e., its largest Hausdorff quotient is a singleton) as we shall
see, so there is no correspondingly nice Thurston compactification on the level of Riemann's moduli space.

\vskip .2in

\noindent{\bf Torus Example}~Recall that for the once-punctured torus $F=F_1^1$, the Teichm\"uller space
is ${\cal T}_1^1\approx\{ z\in{\Bbb C}:|z|<1\}$, and the
mapping class group is
$MC_g^s\approx PSL_2({\Bbb Z})$.  Indeed, the right Dehn twists $M$ and $L$ on the meridian and longitude,
respectively, generate $MC_1^1$, and a complete list of relations between them is given by $\iota = MLM=LML$ and
$\iota ^2=1$.

A point of
${\cal PF}_0(F)$ is uniquely determined by its ``slope'', defined as follows.  Fix two disjointly embedded ideal
arcs $x,y$ asymptotic to the puncture $p$ which decompose $F$ into an ideal quadrilateral, where $x\cup\{ p\}$
is homotopic to the meridian, and $y\cup\{ p\}$ is homotopic to the longitude.
Given
$[{\cal F},\mu ]\in {\cal MF}_0(F_1^1)-\{\varnothing\}$, compact support guarantees that $i_{[{\cal
F},\mu ]}(x)$ and $i_{[{\cal
F},\mu ]}(y)$ are well-defined and finite, and the ratio $|\theta |= i_{[{\cal
F},\mu ]}(y)/i_{[{\cal
F},\mu ]}(x)\in [0,\infty ]$ is therefore projectively well-defined; we further imbue $\theta$ with a sign (when
it is finite and non-zero) in the natural way, where the sign is positive if one (in fact, any) leaf of ${\cal F}$
immediately after meeting
$x$ then meets the copy of $y$ in the frontier of the ideal quadrilateral which lies to the right
(where the orientation of the ideal quadrilateral is inherited from that of $F_1^1$).  It is easy to see that
the {\it slope} $\theta$ is a well-defined and complete invariant of $[{\cal F},\bar \mu ]\in {\cal PF}_0(F_1^1)$,
where we regard $\theta\in S^1$ in the natural way.  Thus, ${\cal PF}_0(F_1^1)\approx S^1$, and $\overline{\cal
T}_1^1={\cal T}_1^1\cup {\cal PF}_0(F_1^1)$ is a closed ball (which you should {\sl not} identify with the Poincar\'e
disk together with its circle at infinity). The slope
$\theta =p/q$ is rational if and only if the measured foliation corresponds to the simple closed curve wrapping $p$
times around the meridian and $q$ times around the longitude, where $p$ and $q$ are relatively prime integers.
Geometrically, deforming hyperbolic structure to pinch this curve, it is clear from elementary considerations of
hyperbolic geometry that the corresponding geodesic length functions converge projectively to the geometric intersection
number with this curve.  Furthermore, the geometric intersection number of curves with slopes $p/q, r/s$ written
in least terms is given by $|ps-qr|$.

\subsection{Train tracks}

\vskip .2in

A {\it train track} $\tau\subseteq F$ is a graph (where vertices are
called ``switches'' and edges are called ``branches'') together with the
following extra structure:

\vskip .1in

\leftskip .3in

\noindent {\it Smoothness}:~$\tau$ is $C^1$ away from its switches.
Furthermore, for each switch $v$ of $\tau$, there is a tangent line $\ell$
to $\tau$ at $v$ in the tangent plane to $F$ at $v$ so that for each
half-branch whose closure contains $v$, the one-sided tangent at $v$ lies
in $\ell$;

\vskip .1in

\noindent {\it Non-degeneracy}:~Vertices of $\tau$ are at least trivalent,
and for any switch $v$ of $\tau$, there is an embedding $(0,1)\to F$ with
$f({1\over 2})=v$ which is $C^1$ as a map into $F$;

\vskip .1in

\noindent {\it Geometry}:~Suppose that $C$ is a component of $F-\tau$, and
let $D(C)$ denote the double of $C$ along the $C^1$ frontier edges of $C$
so the non-smooth points in the frontier of $C$ give rise to punctures
of $D(C)$.  We demand that the Euler characteristic of $D(C)$ be negative.

\vskip .1in

\leftskip=0ex

The smoothness condition is synonymously called the structure of a
``branched one-submanifold'' and leads to the fundamental notion of a
graph smoothly supporting a curve or another train track as we shall see. According to
the non-degeneracy condition, the half-branches incident on a fixed switch
decompose canonically into two non-empty sets of ``incoming'' and
``outgoing'' branches.  The geometric condition rules out the following
complementary regions: smooth disks (i.e., nullgons), monogons, bigons, smooth annuli,
and once-punctured nullgons, and will be further explained below.

Let $B(\tau )$ denote the set of branches of $\tau$.  A function $\mu:B(\tau )\to{\Bbb
R}_{\geq 0}$ induces $\mu:\{ {\rm half-branches~of}~\tau\}\to{\Bbb R}_{\geq 0}$ in
the natural way (where $\mu (b_{1\over 2})=\mu (b)$ if $b_{1\over 2}\subseteq b$) and satisfies
the {\it switch conditions} provided that for each switch
$v$ of
$\tau$, we have
$$
\sum_{{{{\rm outgoing}\atop~{\rm half-branches}~b}}}\mu (b)~~~=~~~
\sum_{{{{\rm incoming}\atop{\rm half-branches}~b}}}\mu (b).
$$

Such a function satisfying the switch conditions is called a {\it
{\rm(}transverse\/{\rm)} measure} on $\tau$, and $\tau$ itself is said to be
{\it recurrent} if it supports a {\it positive} measure $\mu$
with $\mu (b)>0$ for each branch of $\tau$.  In the sequel, train
tracks will tacitly be assumed to be recurrent.

\vskip .1in

\noindent{\bf Torus Example}~There is a unique combinatorial type of recurrent
trivalent train track $\tau$ in the surface $F=F_1^1$, and two embeddings of it as
spine are illustrated in Figure~8. $F-\tau$ consists of a single once-punctured
bigon. There are two branches of $\tau$ so that $\mu$ is uniquely determined by its
values on these branches, and the weight on the remaining branch of $\tau$ is given
by their sum according to the switch condition. Notice that this train track is
``orientable'' in the sense that the graph  underlying $\tau$ admits an orientation
where incoming points toward outgoing at each vertex. Thus, fixing an orientation on
$\tau$, a measure $\mu$ on $\tau$ uniquely determines a homology class in
$H_1(F^1_1,{\Bbb R})$.  Equivalently, every measured foliation of compact support on
$F_1^1$ is transversely orientable.

\vskip .1in

\begin{constr}\label{constr}
{\rm
Given a train track $\tau$
and a positive measure $\mu$ on it, we may construct a measured foliation of a
neighborhood of $\tau$ in the following way.  For each branch $b$ of $\tau$, take a
rectangle of width $\mu (b)$ and length unity foliated by horizontal leaves.  For
each switch, place the rectangles of the incoming branches next to one another and
likewise for the outgoing branches, and then finally glue the vertical edges of all
the incoming to all the outgoing rectangles at each switch in the natural way
preserving the transverse measure along the widths of the rectangles by the switch
conditions. This produces from $\mu$ a measured foliation of a {\it tie
neighborhood} of $\tau$, where a vertical leaf in any rectangle is called a {\it
tie}, and the {\it singular ties} arise from the vertical sides of the rectangles.
As before by Whitehead's result, the Whitehead equivalence class of the resulting
measured foliation is well-defined.
}
\end{constr}

If a measured foliation arises in this way from a measure on a train track, then we say that
the train track {\it carries} the measured foliation.

\vskip .1in

Let $U(\tau )$ denote the cone of all measures on
$\tau$, i.e., the subspace of ${\Bbb R}_{\geq 0}^{B(\tau )}
$ determined by the switch conditions.  There is again the natural ${\Bbb
R}_{>0}$-action on $U(\tau )-\{\vec 0\}$ by homothety, and the quotient
$V(\tau )=(U(\tau )-\{\vec 0\} )/{\Bbb R}_{>0}$ is the {\it polygon of projective
measures} on $\tau$.  Construction~\ref{constr} thus gives well-defined maps
$U(\tau )\to {\cal MF}_0(F)~~{\rm and}~~V(\tau )\to {\cal PF}_0(F)$.

A recurrent train track is {\it maximal} if it is not a proper sub track
of any recurrent train track.  For general $F=F_g^s$, complementary
regions to a maximal train track are either trigons or once-punctured
monogons, but in the special case of the once-punctured torus, a maximal
train track has a single complementary once-punctured bigon.

\vskip .1in

\begin{theorem}~{\rm \cite{Thurs},\cite{HarPen}}~\it For any maximal recurrent train track $\tau$
in
$F$,  Construction \ref{constr} determines continuous embeddings
$$U(\tau )\to{\cal MF}_0(F)~~and~~V(\tau )\to{\cal PF}_0(F)$$
onto open sets.

\end{theorem}

\vskip .1in

In fact, the geometric condition in the definition of train track
precisely guarantees the injectivity in this theorem.  In light of
this result, one may regard a maximal train track in $F$ as indexing
a chart on the manifold ${\cal PF}_0$, and we next study the transition
functions of this putative manifold structure.

\vskip .1in

\noindent{\bf Torus Example}~For $F=F_1^1$, two embeddings of train tracks as spine
are illustrated in Figure~8, and in fact, every foliation is carried by one of
these two train tracks.  The corresponding charts on the circle are also illustrated
as well as the two points of intersection in the closures of these charts.

\centerline{
\epsffile{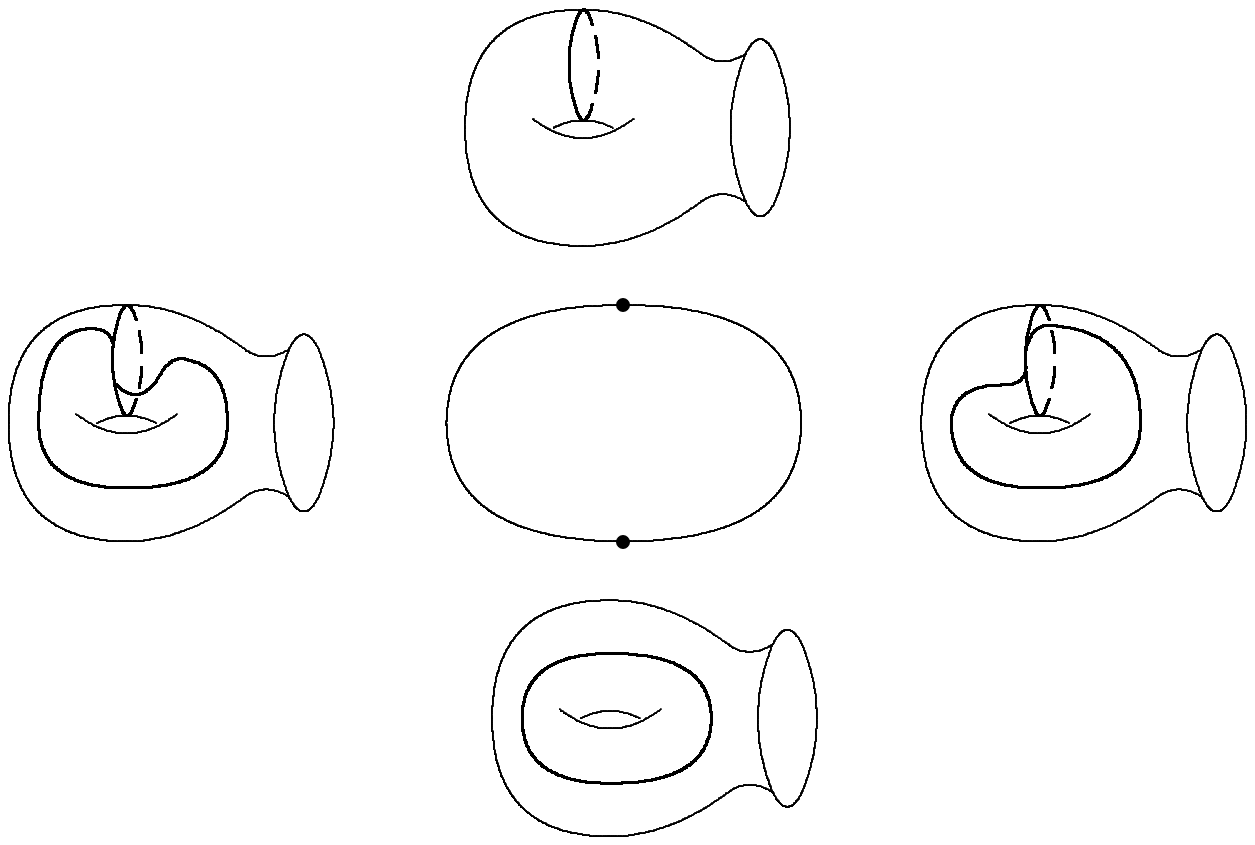}}

\centerline{\bf Figure 8-charts for torus}

\vskip .1in

It is most convenient now to restrict
to the ``generic'' case, where all switches of
$\tau$ are trivalent.  For each switch of $\tau$, the decomposition of
incident half-branches into incoming/outgoing thus consists of one
singleton and one doubleton, and we say a branch of $\tau$ is {\it large}
if it is a singleton at both its endpoints (which are then necessarily
distinct) as illustrated with the branch labeled $e$ in  the left-hand side of Figure~9.
Likewise, if a branch is a doubleton at both its endpoints, then it is called {\it small}, while a
branch which is neither small nor large is called {\it half-large}.

Define the combinatorial
{\it splitting} of a measured train track $(\tau ,\mu )$ along a large
branch $e$ as illustrated in Figure~9, where we identify an edge with its $\mu$ transverse
measure
$\mu (e)$ for convenience.
One imagines
separating bands of horizontal leaves in the rectangle associated to $e$ by
excavating along the two ``singular leaves'' beginning at the endpoints of
the large branch, i.e., beginning at the singular ties.  If the measure $\mu$ is so that
(either of) the singular leaves starting at an endpoint of $e$ turn left or right,
then the respective split is called a
{\it left\/} (case 1) or {\it right\/} (case 3) split, while if the
two singular leaves coincide for $e$, then the
split is called a {\it collision} (case 2).

\centerline{
\epsffile{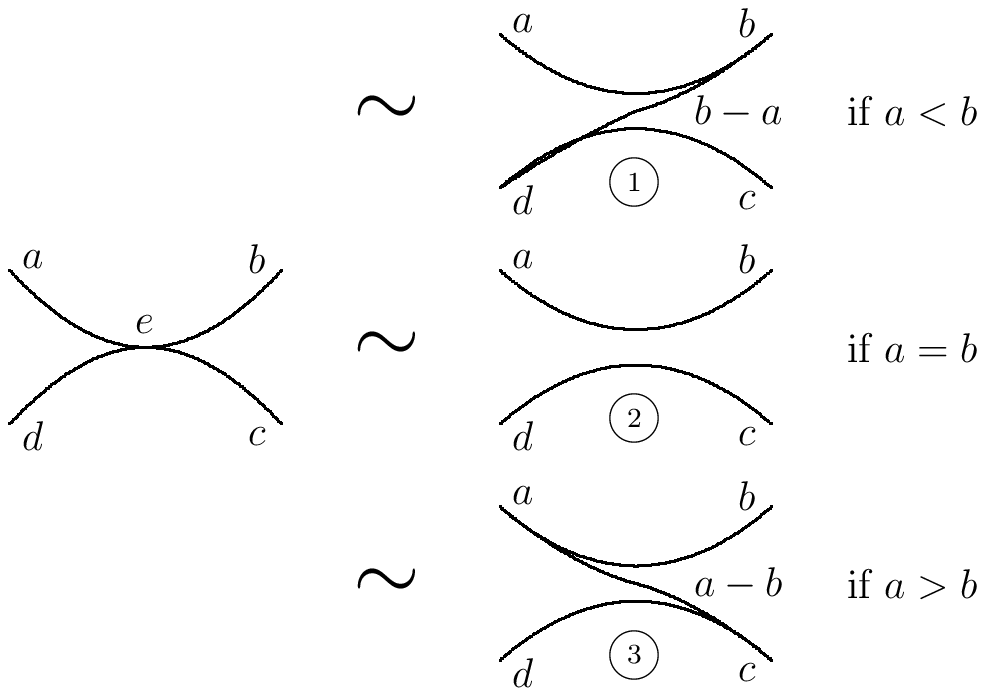}}

\centerline{\bf Figure 9-splitting}

Splitting  and smooth isotopy of measured train tracks generates an equivalence
relation on the set of all measured train tracks in $F$, and we shall let $[\tau ,\mu
]$ denote the equivalence class of the measured train track $(\tau ,\mu )$
and $[\tau ,\bar\mu ]$ denote the equivalence class of the projectively
measured train track $(\tau ,\bar\mu )$.

\begin{theorem}~{\rm \cite{HarPen}}\it
~~If $(\tau _i,\mu _i)$ are positively measured
train tracks giving rise to corresponding measured foliations $({\cal
F}_i,\mu _i)$ via Construction~\ref{constr}, for $i=1,2$, then
$$
[{\cal F}_1,\mu _1]=[{\cal F}_2,\mu _2] ~{ if~and~only~if}~
[\tau _1,\mu _1]=[\tau _2,\mu _2].
$$

\end{theorem}

\noindent Thus, the space of all Whitehead equivalence classes of
(projectivized) measured foliations is identified with the space of all
splitting equivalence classes of (projectivized) measured train tracks up to isotopy.

There is another aspect to the splitting equivalence relation on the set of all
measured train tracks.  In addition to splitting, one considers also {\it shifting}
along a half-large branch by pushing two confluent branches of a train track past
one another as illustrated in Figure~10.  Shifting plays a role in the later
discussion, and a basic result in train track theory \cite{HarPen} is
that if two train tracks are related by shifting, splitting, and smooth isotopy,
then they are also related by splitting and smooth isotopy alone.

\centerline{
\epsffile{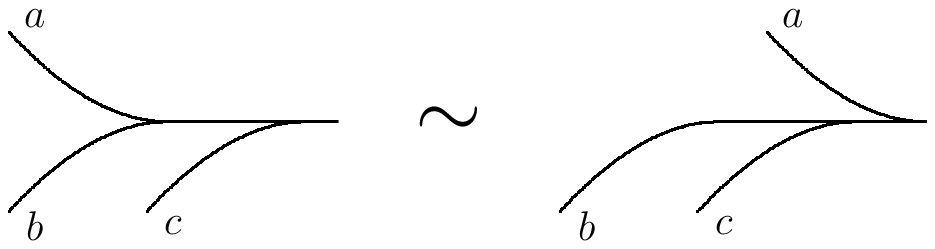}}

\centerline{\bf Figure 10-shifting}

\vskip .1in

Let us finally give the idea of the proof that ${\cal S}$ is dense in ${\cal
MF}_0$, as was mentioned before, by explaining density of ${\cal S}$ in each
chart $V({\cal T})\subseteq {\cal PF}_0$, for some maximal train track
$\tau\subseteq F$.  We may approximate any
$\mu\in V(\tau )\subseteq ({\Bbb R}^{B(\tau )}_{\geq 0}-\{\vec 0\})/{\Bbb R}_{>0}$
by a rational measure $\mu '\in ({\Bbb Q}^{B(\tau )}_{\geq 0}-\{\vec 0\})/{\Bbb R}_{>0}$
(satisfying the switch conditions).  Furthermore clearing denominators in $\mu '$, there are
$N\in{\Bbb Z}_{>0}$  and
$\nu\in({\Bbb Z}^{B(\tau )}_{\geq 0}-\{\vec 0\})/{\Bbb R}_{>0}$
so that $N\mu ' = \nu$.  We may construct an embedded family of
curves in $F$ from $\nu$ by arranging $\nu(b)\geq 0$
tie-transverse strands parallel to $b$ in a tie neighborhood of $\tau$. By
the switch conditions, there are at each vertex exactly as many incoming strands as
outgoing, and there is a unique way to combine strands near vertices to produce a
disjointly embedded family of curves.  Let us give each component
curve a weight $1/N$ and combine any parallel curves while adding their weights to
produce our desired weighted family of disjointly embedded curves.  Let $[{\cal
F},\mu _1]$ denote the corresponding measured foliation (discussed before); tracing
through the constructions, one finds that $[\tau ,\mu ']=1/N [\tau , \nu ]$ gives rise
to $[{\cal F},\mu _1]$.  Letting $\mu '\to\mu$ and projectivizing, it follows easily that families of
disjointly embedded curves are dense in $V(\tau )$.  With a little more work \cite{HarPen},
one can approximate (in the topology of $V(\tau )$) such disjointly embedded families
with a single curve, and this gives the asserted density of ${\cal S}$ itself.

Thus, Thurston's boundary ${\cal PF}_0(F)$ is a completion of the set ${\cal S}$.
In fact, one can approximate with a single non-separating
curve (provided $g\neq 0$); since any two such curves are equivalent under the action of $MC_g^s$,
it follows that the action of $MC_g^s$ on ${\cal PF}_0(F)$ has a dense orbit, and the
maximal non-Hausdorffness of the quotient, which was mentioned before, is thereby
established.

\subsection{Laminations}

\vskip.2in

Each basic formulation of the objects presented so far, namely,
measured foliations and measured train tracks, requires passage to the
quotient under an appropriate equivalence relation.  Thurston has given a more
ethereal, elemental, and elegant description of these objects as ``measured geodesic
laminations'', where no passage to equivalence classes is
necessary.  Here we simply give the definition and a few basic
properties referring the reader to \cite{HarPen} for instance for further details.  A ``lamination''
${\cal L}$ in $F$ is a foliation of a closed subset of $F$, and a ``(transverse)
measure'' to ${\cal L}$ is defined much as before as a $\sigma$-additive measure on
arcs transverse to ${\cal L}$ with the analogous condition of no-holonomy (where the
homotopy is through arcs transverse to ${\cal L}$ with endpoints disjoint from ${\cal
L}$).
${\cal L}$ is a ``geodesic lamination'' if its leaves are geodesic for some specified
hyperbolic metric.  The simplest case of a geodesic lamination is the geodesic representative
of an element of ${\cal S}'(F)$. (In fact, for different
choices of metric, the spaces of measured geodesic laminations are naturally identified via the circle at
infinity in their universal covers, so we may speak of a geodesic lamination without the {\it a
priori} specification of a metric.)  A measured geodesic lamination has zero measure
in $F$, and the intersection of a measured
geodesic lamination
${\cal L}$ with a transverse arc $a$ in
$F$ is a Cantor set together with isolated points corresponding to intersections with
simple geodesic curve components or arc components of ${\cal L}$, if any.
There is a natural topology on the set of all measured geodesic laminations in
$F$, which is induced by the weak topology on the set of all
$\pi _1(F)$-invariant measures supported on the M\"obius band past infinity.
${\cal ML}(F)$ (and ${\cal ML}_0(F)$) is the corresponding space of measured geodesic
laminations (and with compact support) and corresponding projectivization ${\cal
PL}(F)$ (and ${\cal PL}_0(F)$).  A basic result in Thurston theory is ${\cal
ML}(F)\approx {\cal MF}(F)$, ${\cal
ML}_0(F)\approx {\cal MF}_0(F)$,${\cal
PL}(F)\approx {\cal PF}(F)$, ${\cal
PL}_0(F)\approx {\cal PF}_0(F)$, where train tracks give suitable charts on any of these piecewise
linear manifolds. Furthermore, the deformation theory due to Thurston, called
``earthquaking'' cf. \cite{Ker2}, which we do not further discuss here, is most
conveniently expressed in the context of laminations.

\vskip .2in

\subsection{Dynamics on train tracks}

\vskip .2in

In this section, we simply recall Thurston's classification of surface automorphisms
as well as recall several basic facts about ``pseudo-Anosov'' mappings.  In the
process, we develop further basic techniques which will be required in quantization.
Since Thurston's compactification produces a closed ball upon which the mapping
class group acts continuously, one immediately is led to consider fixed points of
this action.

\begin{theorem} {\bf Thurston's Classification}~{\rm ~~\cite{Thurs},\cite{FLP}
}~~\it Any orientation-preserving homeomorphism $f:F\to F$ is homotopic to
a diffeomorphism $f':F\to F$ which satisfies one of the following conditions
(and the only overlap is between 1. and 2.):

\vskip .1in

\leftskip .3in

\noindent {\bf 1.}~$f'$ fixes a unique point of ${\cal T}$ and is of finite order;

\vskip .1in

\noindent {\bf 2.}~$f'$ is ``reducible'' in the sense that $f'$ fixes an element of
${\cal S}'$;

\vskip .1in

\noindent {\bf 3.} $f'$ is ``pseudo-Anosov'' in the sense that there is some $\lambda
>1$ together with two measured foliations
$({\cal F}_\pm ,\mu _\pm)$ , which share singular points and
are otherwise transverse, so that $f'({\cal F}_\pm,\mu _\pm)=\lambda ^{\pm 1}~({\cal
F}_\pm ,\mu _\pm )$.  The projective classes $[{\cal F}_\pm ,\bar \mu _\pm ]$ are the
unique fixed points of $f$ on $\overline{\cal T}$.  The invariant $\lambda$ is called the
{\it dilatation} of $f$ or $f'$.

\leftskip=0ex

\end{theorem}

Notice the similarity with the trichotomy elliptic/parabolic/hyperbolic for fractional
linear transformations.   In the reducible case, one simplifies the dynamics by
cutting
$F$ along a representative of the invariant element of ${\cal S}'$.  A pseudo-Anosov
mapping is the analogue of an Anosov map of the torus in the current context of
surfaces with negative Euler characteristic.

In fact, train tracks provide a powerful tool for analyzing the dynamics of surface
automorphisms owing to the fact that since a train track has a well-defined tangent
line at each point, there is a coherent notion of a train track smoothly ``carrying''
a curve, another train track, or a lamination.

Suppose that $\kappa$ is a smooth curve, a train track, or a measured geodesic
lamination.  We say that the train track
$\tau$ {\it carries} $\kappa$ and write $\kappa <\tau$, if there is a $C^1$
map $\phi :F\to F$ homotopic
to the identity, called the {\it supporting map}, so that $\phi (\kappa )\subseteq \tau$,
where the restriction of the differential $d\phi _p$ to the tangent line to $\kappa $
at $p$ is non-zero for every $p\in \kappa $.

We think of $\phi$ as squashing together nearly parallel strands of $\kappa$.  For
instance, any curve arising as before from an integral measure on $\tau$ is carried
by $\tau$, and more generally, any curve, train track or lamination $\kappa$ which
lies in a tie-neighborhood of $\tau$ and is transverse to the ties satisfies
$\kappa <\tau$, where the supporting map collapses ties.

For example, if a train track $\sigma$ arises from $\tau $ by splitting and shifting
(but no collapsing), then
$\sigma <\tau$, and we say that $\sigma$ arises from $\tau$ by {\it unzipping}.

Suppose that the train track $\sigma$ is transverse to a tie-neighborhood of
$\tau$, say with supporting map $\phi :F\to F$.  Let us enumerate the branches $b_j$
of $\tau$, for $j=1,\ldots ,n$, and
$a_i$ of $\sigma$, for $i=1,\ldots ,m$, and choose $x_j\in b_j$ for each $j$.
There is then an $m\times n$-matrix
$A=(A_{ij})$ called the {\it incidence matrix}, where $A_{ij}$ is the cardinality
of $\phi ^{-1}(x_j)\cap a_i$.
It is clear that the incidence matrix $A:{\Bbb R}_{\geq 0}^m\to{\Bbb R}_{\geq 0}^n$
describes the inclusion $U(\sigma )\to U(\tau )$ in the train track coordinates.

We close this section with several basic results about pseudo-Anosov mappings.

\begin{theorem}~{\rm \cite{Thurs},\cite{PP1}}\it  ~~A homeomorphism $f:F\to F$
is a pseudo-Anosov map if and only if there is a train track $\tau$ in $F$, with each
component of $F-\tau$ an at most once punctured polygon, so that
$\tau$ unzips to $f( \tau)$ with no collisions.   Furthermore, the incidence matrix
$A$ of the carrying
$f(\tau)<\tau$ is Perron--Frobenius, the eigenvector of $A$ corresponding to the
spectral radius $\lambda$ gives the projective measure
$\bar\mu _+$, and likewise the extreme eigenvector of the transpose of $A$
gives $\bar \mu _-$.
\end{theorem}

Given a measured train track $(\tau
,\mu )$, consider the foliated neighborhood of $\tau$ determined by $\mu$ via
Construction~\ref{constr}.  Choose some enumeration of the switches of $\tau$ and serially
follow the singular leaves from the switches until the first splitting (ignoring shifting),
for the first
switch, second switch, $\ldots$ , last switch, and then begin anew from the first
switch.  Suppose there are no collisions, and record the resulting sequence of right
or left splits, so as to produce a semi-infinite word of rights and lefts.

\begin{theorem}{\rm \cite{PP1}}   The right-left
sequence is eventually periodic if and only if the corresponding measured foliation
is fixed by some pseudo-Anosov mapping.
\end{theorem}

An explicit and simple construction of pseudo-Anosov maps is given by the following result.

\begin{theorem} {\rm \cite{Penn4}}~ Suppose that ${\cal C},{\cal D}\in {\cal S}'$
admit representative arc families $C,D$ intersecting minimally which satisfy the
condition that each component of
$F-\cup({ C}\cup{ D})$ is an at most once-punctured polygon.  Take any composition
$w$ of Dehn twists to the right along elements of ${\cal C}$ and to the left along
elements of
${\cal D}$ so that the Dehn twist along each element of ${\cal C}$ or ${\cal D}$
occurs at least once in
$w$.  Then
$w$ is pseudo-Anosov.

\end{theorem}

\vskip .1in

\noindent{\bf Torus Example}~Consider a generic train track $\tau$ in $F_1^1$, so
$\tau$ has one large branch $e$ and two small branches.  The two small branches are
canonically linearly ordered by first taking the branch $a$ to the right and then
the branch $b$ to the left at either endpoint of $e$, and furthermore $U(\tau
)\approx{\Bbb R}_{\geq 0}^{\{ a,b\}}$, i.e., the measures of the small branches $a,b$
are unconstrained and uniquely determine the measure on $e$ as well.  Given a
measure $\mu\in U({\tau })-\{ \vec 0\}$, start unzipping $(\tau ,\mu )$ along either
singular leaf, i.e., split along $e$, to produce another measured train track $(\tau
_1,\mu _1)$; of course, $\tau _1$ is combinatorially equivalent to $\tau$.  For
definiteness, suppose that $B=\mu (b)>\mu (a)=A$, so the split is a left split.  The
edge corresponding to $b$ is the large edge of $\tau _1$, and the two small edges,
in right/left order, have measures $(A,B-A)$.  Continue unzipping, i.e., next split
$(\tau _1,\mu _1)$ along its large edge to produce $(\tau _2,\mu _2)$. Again suppose
that $B>2A$ for definiteness, so the second split is a left split as well, and the
small branches of $\tau _2$ have measures $(A,B-2A)$.  Continue unzipping (under the
assumption that there are no collisions) until the first right split, say there are
$a_1=[{B\over A}]>1$ left splits before the first right split.  Perform the right
split along  the large branch of $(\tau _{a_1},\mu _{a_1})$, where the measures on
the small branches are $(A,B_1)=(A,B-m_1A)$, to produce  a train track whose small
edges have measures $(A-B_1,B_1)$.  Continue unzipping and suppose there are no
collisions, i.e., suppose $A$ and $B$ are not rationally related, to produce a
semi-infinite sequence of symbols $L$ (for left splits) and $R$ (for right splits).
Let $a_1$ denote the number of $L$'s that begin this sequence, $a_2$ denote the
length of the next consecutive sequence of $R$'s, $a_3$ the length of the next
consecutive sequence of $L$'s, and so on.  It follows from the discussion above that
the continued fraction expansion of $B/A$ is given by

\bea\nonumber
B/A={a_1+{\displaystyle\frac{1}{
  \renewcommand{\arraystretch}{0}
  \arraycolsep=0.05em
  \begin{array}{ccc}
   \strut a_2+{\displaystyle\frac{1}{
  \renewcommand{\arraystretch}{0}
  \arraycolsep=0.05em
  \begin{array}{ccc}
   \strut a_3+{} & & \\
     & \ddots & \\
       \end{array}
     }} & & \\
       \end{array}
     }}}.
\eea

Continued fractions occur in another related guise as well.  The isomorphism
$MC_1^1\approx PSL_2({\Bbb Z})$ is induced by
$M\mapsto\pmatrix{1&1\cr 0&1\cr}$ and
$L\mapsto\pmatrix{1&0\cr -1&1\cr}$, where $M$ and $L$ are the right Dehn twists on
the meridian and longitude respectively.
A mapping class in  $MC_1^1$ is pseudo-Anosov, reducible, periodic if and only if the
corresponding fractional linear transformation is hyperbolic, parabolic, elliptic respectively.
Every hyperbolic element of $PSL_2({\Bbb Z})$ is conjugate to a product
$$
\pmatrix{1&m_1\cr 0&1\cr}~\pmatrix{1&0\cr n_1&1\cr}~\cdots ~\pmatrix{1&m_k\cr 0&1\cr}
~\pmatrix{1&0\cr n_k&1},
$$
where $m_i,n_i>0$ are unique up to cyclic permutation.  Furthermore,
$m_1,n_1,\ldots ,m_k,n_k$ are the partial quotients of
the periodic continued fraction expansion of
the dilatation of the corresponding pseudo-Anosov map. It follows from this
discussion that in $MC_1^1$, all pseudo-Anosov mappings arise from the previous
theorem. Indeed, the theorem gives the construction of two semi-groups corresponding
to right/left or left/right twisting on meridian/longitude.  For each semi-group, it
is easy to construct a train track $\tau$ in $F_1^1$ so that the matrix
representation above precisely describes the action of the corresponding semi-group
on the measures of the linearly ordered small branches of $\tau$; indeed, these two
train tracks are illustrated in Figure~8.

\vskip .2in

\subsection {Decorated measured foliations and freeways}

\vskip .2in

In this section, we recall material from \cite{PP2} which is required for quantization. If
$\Gamma\subseteq F$ is a cubic fatgraph spine of $F$, then we may blow-up each vertex of
$\Gamma$ into a little trigon as illustrated in Figure~11.  The resulting object
$\tau=\tau _\Gamma$ has both a natural branched one-submanifold structure and a fattening, and
furthermore, components of $F-\tau$ are either little trigons or once-punctured
nullgons.  Thus, $\tau$ is not a train track, but it is almost a train track, and is
called the {\it freeway} associated to $\Gamma$.  Notice that each edge of $\Gamma$ gives rise
to a corresponding large branch of $\tau$, and each vertex gives rise to three small
branches.  It is easy to see that every measured lamination of compact
support in $F$ is carried by the freeway $\tau$.    The frontier of a once-punctured
nullgon component of $F-\tau$ is a puncture-parallel curve called a {\it collar
curve} of $F$. A small branch is contained in exactly one collar curve, while a
large branch may be contained in either one or two collar curves.

\vskip .1in

\setlength{\unitlength}{0.3mm}%
\begin{picture}(50,120)(-250,-40)
\thicklines
\put(-130,19){\line(-2,-1){50}}
\put(-130,19){\line(2,-1){50}}
\put(-130,19){\line(0,1){47}}
\put(-60,25){\makebox(0,0)[cc]{\huge $\to$}}
\put(-15,2.5){\line(-2,-1){25}}
\put(35,2.5){\line(2,-1){25}}
\put(10,51){\line(0,1){25}}
\qbezier(-15,2.5)(-2.5,24.25)(10,51)
\qbezier(35,2.5)(22.5,24.25)(10,51)
\qbezier(-15,2.5)(10,2.5)(35,2.5)
\put(80,25){\makebox(0,0)[cc]{\huge $\sim$}}
\qbezier(120,0)(150,17.3)(180,0)
\qbezier(120,0)(150,17.3)(150,56)
\qbezier(150,56)(150,17.3)(180,0)
\put(120,0){\line(-2,-1){20}}
\put(180,0){\line(2,-1){20}}
\put(150,56){\line(0,1){20}}
\end{picture}

\centerline{\bf Figure 11--freeway from fatgraph}

\vskip .1in

A {\it measure} on a freeway $\tau$ is a function $\mu \in {\Bbb R}^{B(\tau )}$
satisfying the switch conditions, where we wish to emphasize that the measure is not
necessarily nonnegative (as it is for train tracks).  Let $U(\tau )$ denote the
vector space of all measures on $\tau$.  Notice that $\mu\in U(\tau )$ is uniquely
determined by its values on the small branches alone, and the switch conditions are
equivalent to the following ``coupling equations''
$$
\mu (a_1)+\mu (b_1)=\mu (e)=\mu (a_2)+\mu (b_2),
$$
for any large branch $e$ whose closure contains the switches $v_1\neq v_2$, where
$a_i,b_i$ are the small branches incident on $v_i$ for $i=1,2$. On the other hand,
the values on the large branches alone also uniquely determine $\mu$, and in fact,
these values are unconstrained by the switch conditions.  Indeed, letting $a_i$
denote the large branches incident on a little trigon with opposite small branches
$\alpha _i$, for $i=1,2,3$, we may uniquely solve for a measure $\mu$ on $\tau$,
where
\bea
\nonumber
\mu (\alpha _i )={1\over 2}\{ \mu (a_1)+\mu (a_2)+\mu (a_3) -2\mu (a_i)\},~{\rm for}~i=1,2,3,
\eea
and so we identify $U(\tau )\approx{\Bbb R}^{LB(\tau )}$, where $LB(\tau )$
denotes the set of large branches of $\tau$.

In particular, if $\mu$ is a nonnegative measure on $\tau$, then the analogue of
Construction~\ref{constr} in the current context produces a well defined equivalence
class of measured foliations in $F$, where this measured foliation will typically
contain a collection of puncture-parallel annuli foliated by curves homotopic to
collar curves.  Deleting these foliated annuli produces a well-defined (but possibly
empty) class in ${\cal MF}_0(F)$. Thus, a nonnegative measure on $\tau$ canonically
determines a point of ${\cal MF}_0(F)$ together with a nonnegative ``collar
weight'', i.e., the transverse measure of a transverse arc connecting the boundary
components of the corresponding foliated annulus.

In the general case that $\mu$ is not necessarily nonnegative, suppose that $C$ is a
collar curve of $\tau$.  The switches of $\tau$ decompose $C$ into a collection of
arcs, each of which inherits a corresponding real-valued weight from $\mu$.  Let
$\{\gamma _i\} _1^n$ denote the collection of real numbers associated to the small
branches of $\tau$ that occur in $C$.  Define the {\it collar weight} of $C$ for
$\mu$ to be $\mu _C=\min\{\gamma _i\} _1^n$.  Define a {\it collar weight} on $F$
itself to be the assignment of such a weight to each puncture.

We may modify the original measure $\mu\in U(\tau )$ by defining
$\mu '(b)=\mu (b)-\mu _C$ if $b$ is contained in the collar curve $C$ for any
small branch $b$ of
$\tau$.  Thus, $\mu '$ is a nonnegative measure on $\mu$ with identically vanishing
collar weights that determines a corresponding element of ${\cal MF}_0(F)$.

We are led to define the space
$\widetilde{\cal MF}_0(F)={\cal MF}_0(F)\times{\Bbb R}^s$ of {\it decorated measured
foliations} and summarize the previous discussion:

\begin{theorem} ~{\rm \cite{PP2}}~
The space $U(\tau)\approx {\Bbb R}^{LB(\tau )}$ gives global coordinates on
$\widetilde{\cal MF}_0$, and there is a canonical fiber bundle
$\Pi :\widetilde{\cal MF}_0\to{\cal MF}_0$, where the fiber over a point is the set
${\Bbb R}^s$ of all collar weights on $F$.
\end{theorem}

\begin{remark}
{\rm
The natural action of $MC_g^s$ is by bundle isomorphisms of $\Pi$.  Furthermore,
$\Pi$ admits a natural $MC_g^s$-invariant section
$\sigma :{\cal MF}_0\to\widetilde{\cal MF}_0$
which is determined by the condition  of identically
vanishing collar weights.  The restriction of $\sigma$ to
${\cal MF}_0\subseteq\widetilde{\cal MF}_0$ gives a piecewise-linear embedding of the
piecewise-linear manifold ${\cal MF}_0$ into the linear manifold (vector space)
$\widetilde{\cal MF}_0\approx U(\tau)\approx {\Bbb R}^{LB(\tau )}$.
}
\end{remark}

\subsection{Shear coordinates for measured foliations}\label{scmf}

We now give an equivalent parametrization of measured foliations in terms of
``Thurston's shear coordinates'' that are close
analogues of
Thurston's shear coordinates $Z_\alpha$ on ${\cal T}_H(F)$. In fact, we have already
encountered these quantities when describing the splitting procedure train tracks (see Figure~9).
There, excavating along two different singular leaves, we have obtained the ``new'' edge, which can
turn either left or right (for splittings) or be absent (for collisions).

We assign a corresponding signed quantity (positive for right, negative for left) as follows.
Given a measure $\mu$ on the long branches of the freeway $\tau$ associated to the fatgraph spine
$\Gamma\subseteq F$, define the {\it (Thurston's foliation-)shear coordinate} of the edge indexed by
$\alpha$ to to be $$\zeta_\alpha ={1\over 2}(\mu (A)-\mu (B)+\mu (C)-\mu (D)),$$
in the notation of Figure~1 for nearby branches.
>From the very definition, $\zeta _\alpha$ is independent of collar weights.
Again, Thurston's foliation-shear coordinates are alternatively defined in terms of the
signed transverse length of the arc between the singular leaves along $Z_\alpha$, in analogy
to the geometric interpretation given before of the
shear coordinates on ${\cal T}_H$.

Note that the shear coordinates $\zeta _\alpha$ are not independent. They are subject to the
restrictions that
\be
\label{facecond}
\sum_{\alpha\in I}\zeta _\alpha=0
\ee
for the sum over edges $\alpha\in I$ surrounding any given boundary component,
and we shall refer to these conditions as the {\it face conditions} for shear coordinates.
Thus, the space of foliation-shear coordinates is of dimension
${LB(\tau)-n}$, where we let
$n$ denote the number of boundary components.
One sees
directly that for any assignment of shear coordinates, there is a well-defined point of ${\cal MF}_0$
realizing them, thereby
establishing a homeomorphism between ${\cal MF}_0$ and this sub-vector space ${\Bbb R}^{LB(\tau
)-n}\subseteq {\Bbb R}^{LB(\tau )}$ of shear coordinates on the long branches of $\tau$.

To describe the action of the mapping class group on foliation-shear coordinates,
we shall give the transformation under Whitehead moves, i.e., derive the analogue of
formula (\ref{abc}) for measured foliations, which is an elementary calculation using the formulas for
splitting as follows.

\begin{lemma}\label{lem4-1}
Under the Whitehead move in Figure~2, the corresponding foliation-shear coordinates
of the edges $A$, $B$, $C$, $D$, and $Z$ situated as in Figure~2 are transformed
according to formula (\ref{abc}) $$M_Z:
(\zeta _A,\zeta _B,\zeta _C,\zeta _D,\zeta _Z)\quad\mapsto\quad\quad\quad\quad\quad\quad\quad$$
$$\quad\quad\quad(\zeta _A+\phi_H(\zeta _Z),
\zeta _B-\phi_H(-\zeta _Z), \zeta _C+\phi_H(\zeta _Z), \zeta _D-\phi_H(-\zeta _Z),-\zeta _Z)$$ with
\begin{equation}\label{h-trans}
\phi_H(\zeta_Z)=(\zeta_Z+|\zeta_Z|)/2,
\end{equation}
i.e., $\phi_H(x)=x$, for $x>0$, and zero otherwise.  All other
shear coordinates on the graph remain unchanged.
\end{lemma}

\centerline{
\epsffile{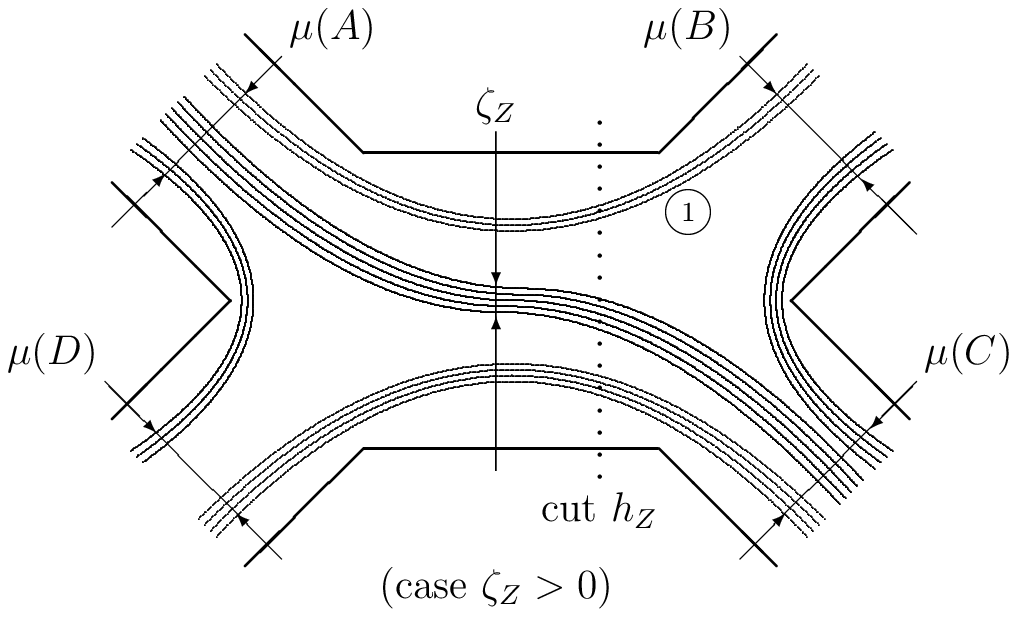}}

\centerline{\bf Figure~12-foliation-shear coordinates}

\vspace{10pt}

\begin{remark}\label{rem4-1}
{\rm
Comparing expressions for
the classical function $\phi(x)=\log(1+\e^x)$ and (\ref{h-trans}), one finds that the
latter is a {\em projective limit\/} of the former:
\begin{equation}\label{h-tr2}
\phi_H(x)=\lim_{\lambda\to+\infty}\frac1\lambda \phi(\lambda x)
=\lim_{\lambda\to+\infty}\frac1\lambda \phi^\hbar(\lambda x),
\end{equation}
that is, {all three} transformations coincide {asymptotically}
in the domain of large absolute values {\rm(}or large eigenvalues for the
corresponding operators\/{\rm)}
of Teichm\"uller space coordinates $\{Z_\alpha\}$.
We shall actively use this property in Section~\ref{mcg-sect}
when proving the existence of the quantization of Thurston's boundary for the punctured torus.
}
\end{remark}

\clearpage

\newsection{On quantizing Thurston theory}

\subsection{Proper length of geodesics}

\begin{defin}\label{6-1}
{\rm
The {\em proper length\/} ${\rm p.l.}(\gamma)$ of a closed curve
$\gamma$ in the classical or quantum case is
constructed from the
quantum ordered operator $P_\gamma$
associated to a closed oriented edge-path with basepoint
(to begin the linearly ordered word $P_\gamma$) as
\be
\label{proj-l}
{\rm p.l.}(\gamma)=\lim_{n\to\infty}\frac1n
~{\rm tr}~\log 2T_n(P_\gamma/2),
\ee
where we take the principal branch of the logarithm and $T_n$ are Chebyshev's
polynomials (cf. {\rm (\ref{cheb})}). Since $T_n(\cosh {t\over 2})=\cosh {{nt}\over 2}$,
it follows that ${\rm p.l.}(\gamma )$ agrees with half the
hyperbolic length of $\gamma$ in the Poincar\'e metric in the classical case.

More explicitly in the operatorial case, we can determine ${\rm p.l.}(\gamma )$ explicitly in terms of the
spectral expansion of the operator $G_\gamma$, which is known exactly.  Namely, the basis of
eigenfunctions of $G_\gamma$ is ``doubly reduced'' in the sense that each eigenvalue (except 2, which
is singular) with corresponding eigenfunction $\alpha _S$, has the form $\rm e ^{S/2}+\e ^{-S/2}$,
where $S$ ranges over the entire real axis, and $\alpha _S$ has the same eigenvalue as $\alpha _{-S}$.
In fact, these functions coincide, so there is actually a representation on the positive real axis,
which is nevertheless complete, and is singular at infinity and zero.  We may define the proper length
operator to the be one with the same eigenfunctions $\alpha _S$ for $S$ positive (which constitute a basis
in the function space) and with eigenvalues to be $|S/2|$.  This operator ${\rm p.l.}(\gamma )$ is then a
well-defined operator on any compactum in function space.

The {\em proper length\/} of a QMC or GMC $\hat C$,
again denoted ${\rm p.l.}(\hat C)$, is the sum of the proper lengths
of the constituent geodesic length operators
(or the sum of half geodesic lengths calculated in the Poincar\'e metric
in the classical case) weighted by the number of appearances in the multiset.
}
\end{defin}

\subsection{Approximating laminations and the main theorem}

Fix once and for all a spine $\Gamma$ of $F$ with corresponding freeway $\tau$.  A measure $\mu$ on
$\tau$ gives rise to a (possibly empty) measured foliation in $F$ together with a
collar weight on the boundary components of $F$.  Erasing collars yields an underlying measure $\mu
_1\geq 0$ on $\tau$ whose support is a sub-train track $\tau _1\subseteq\tau$, and
the measured train track $(\tau _1,\mu _1)$ determines a (possibly empty) measured foliation.   Via
the canonical embedding of
${\cal MF}_0$ into
$\widetilde{\cal MF}_0$ with vanishing collars, we may thus uniquely determine a point of ${\cal
MF}_0$ by specifying foliation-shear coordinates on the
long branches of $\tau$ satisfying the face conditions \ref{facecond}, i.e., ${\cal MF}_0$ is
naturally identified with a codimension $n$ subspace of ${\Bbb R}^{LB(\tau )}$.  Passing to projective
foliations, a point of
${\cal PF}_0$ is given by the projectivization
$P\vec\zeta$ of a vector of foliation-shear coordinates $\vec\zeta=(\zeta _i)$, where $\vec\zeta\in
{\Bbb R}^{LB(\tau )}-\{\vec 0\}$, and $i$ indexes the long branches of $\tau$, i.e., the edges of
$\Gamma$.

\begin{defin}\label{lam1}
{\rm
A sequence $\vec n^\beta=(n_i^\beta\}$, for $\beta \geq 1$, of integer-valued
$n_i$, for $i =1,\ldots ,LB(\tau )$, on
$\tau$ is an {\em approximating sequence\/} for the projectivized measured
foliation $P\vec \zeta$
if the face conditions \ref{facecond} hold on $\vec n^\beta$ and if
$\lim_{\beta\to\infty}n_i^\beta/n_j^\beta=\zeta
_i/\zeta _j$ for all
$i,j$ with $\zeta _j\neq 0$.
}
\end{defin}

Constructed from $\vec n$ as an integral measure on $\tau$ is a GMC $\hat C$ with integral collar
weights.  Just as with decorated measured foliations, components of $\hat C$ which are puncture- or
boundary-parallel can be erased to produce a corresponding multicurve to be denoted $\hat C _{\vec
n}$.  $\hat C _{\vec n}$ is carried by a sub-train track of $\tau$, and it traverses the long
branch of $\tau$ indexed by $i$ some number, say, $m_i\geq 0$ of times, so $m_i$ is the standard train
track coordinate of integral transverse measure.  In the usual notation as in Figure~12, one sees
directly that
$n_Z=\frac12
\bigl(m_A-m_B+m_C-m_D\bigr)$.  We may also sometimes write $\hat C_{\vec m}$ for $\hat C_{\vec n}$

\begin{defin}\label{lam4}
{\rm
A {\em graph length\/} function with respect to the spine $\Gamma$
is any linear function
\be
{\rm g.l.}_\Gamma^{\vec a}(\hat C_{\vec n})={\rm g.l.}_\Gamma^{\vec a}(\hat C_{\vec m})=\sum_i a_im_i.
\ee
In particular, when all $a_i$ are unity, the
graph length is just the combinatorial length of
$\hat C_{\vec m}$, i.e., the total
number of edges of $\Gamma$
traversed (with multiplicities) by all the
component curves of $\hat C_{\vec m}$.
When the spine $\Gamma$ and $\vec a$ are fixed or unimportant, then we shall write
simply ${\rm g.l.}(\hat C_{\vec n})$ or ${\rm g.l.}(\hat C_{\vec m})$ for the
graph length. }
\end{defin}

Any graph length function is evidently additive over disjoint unions of
multicurves, and more generally, is a linear function of $\vec m$.

We next describe the bordered or punctured torus case
in detail.
Each multicurve on the torus is uniquely determined by
three nonnegative integers $(m_X,m_Y,m_Z)$ that satisfy one of the three
triangle {equalities} $m_i=m_j+m_k$,
where $\{ i,j,k\} =\{ X,Y,Z\}$.
Projectivization allows us to re-scale so that
$m_j$ and $m_k$ are relatively prime using this degree of freedom, so that the
corresponding multicurve has just one component.

As in Figure~8, the space ${\cal
PF}_0(F_1^1)$ of projectivized measured foliations with compact support
is a piecewise-linear circle for $r+s=1$, and an alternative family of charts on this circle
is given in Figure~13.  The relation $\sim$ in Figure~13 denotes the equivalence
between different boundary cases between two different charts, and arrows
represent one-simplices in ${\cal PF}(F_1^1)$.

\centerline{
\epsffile{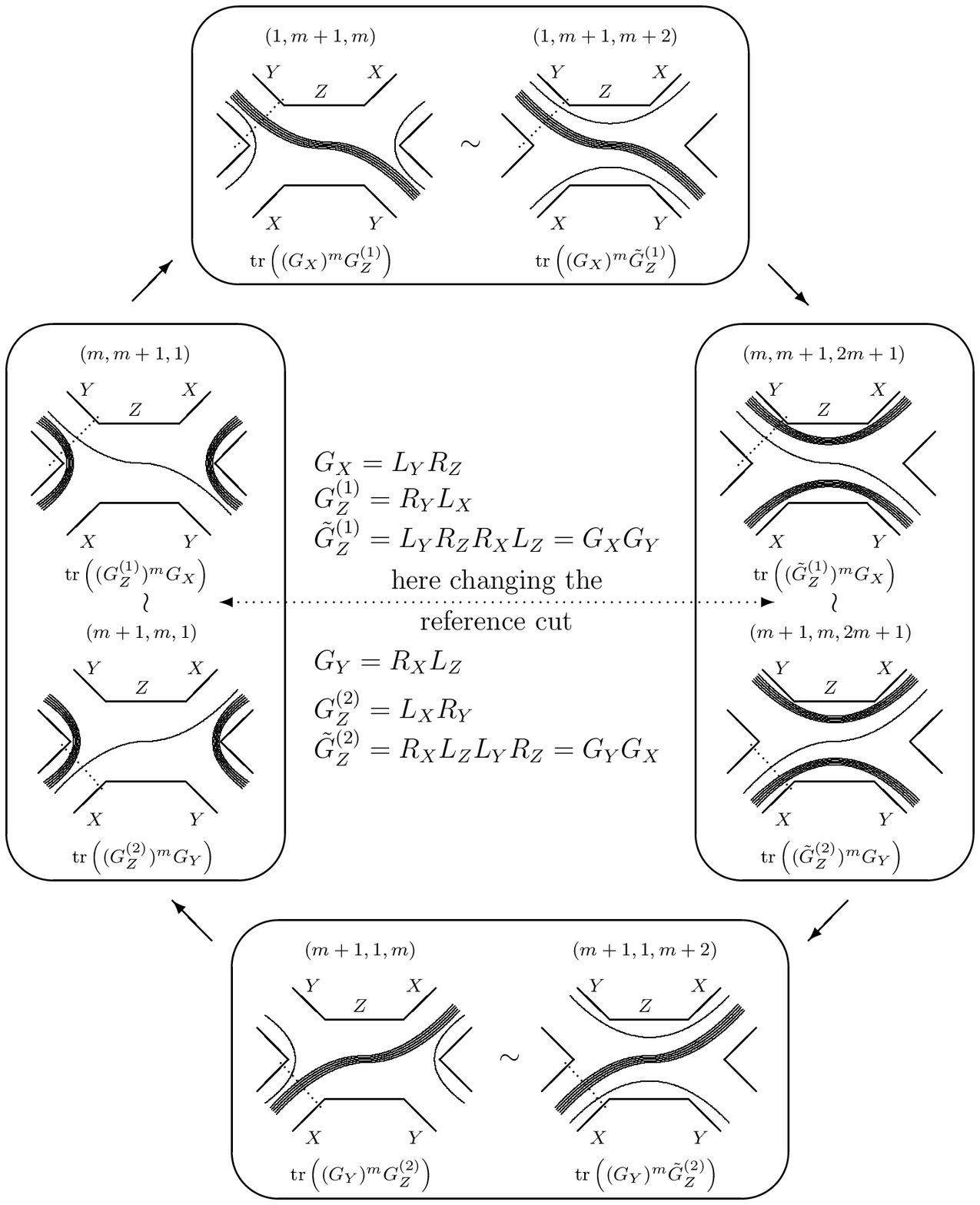}}

\vspace{10pt}

\centerline{\bf Figure 13--the circle ${\cal PL}_0(F_1^1)$}

\vspace{10pt}

In Figure~13, we use the previous notation $G_X$, etc. (see (\ref{torus-classic}))
but in a slightly different sense.
Now, these quantities are $(2\times 2)$-{\it matrices\/}, not just geodesic
functions, i.e., we do not evaluate traces in the corresponding formulas.
There is thus an ambiguity in choosing the place in the graph where the
matrix products begin. We indicate this place by drawing the reference cut
(the dotted line). Changing the reference cut when moving along the circle
corresponds to passing from one chart to another in the chart
covering of the circle.  Of course, choosing the reference cut does not affect the
quantum trace operation.

In order to have a good transition
in the boundary cases, for instance, in the upper case in Figure~13 (and the other cases are similar and
omitted), we must ensure that the corresponding functions for the quantities $G_X^mG_Z^{(1)}$ and
$G_X^m{\tilde G}_Z^{(1)}$ must coincide in the limit $m\to\infty$ with each
other and with the corresponding quantity calculated merely for the ``short'' geodesic
function $G_X$. To prove this, given two
$(2\times 2)$-matrices $G_X$ and $G_Z$ corresponding to geodesic curves, we
can conjugate them by respective unitary transformations
$U_X$ and $U_Z$ to diagonal form with
real eigenvalues $\e^{\pm l_X/2}$ and $\e^{\pm l_Z/2}$ since $G_X,G_Z$
are hyperbolic. We then have
\be
\label{7-1}
\tr G_X^mG_Z =\tr\left(%
\begin{array}{cc}
  \e^{ml_X/2} & 0 \\
  0 & \e^{-ml_X/2} \end{array}%
\right)
V\left(%
\begin{array}{cc}
  \e^{l_Z/2} & 0 \\
  0 & \e^{-l_Z/2} \\
\end{array}%
\right)V^{-1}, \quad V=U_X^{-1}U_Z
\ee
and the proper length (\ref{proj-l}) is
$ml_X/2+O(1)$.\footnote{Unless the matrix $U_X^{-1}G_ZU_X$
has the form
$\left(%
\begin{array}{cc}
  0 & -r^{-1} \\
  r & \phi \\
\end{array}%
\right)$.  For this matrix to
determine a hyperbolic element, the quantity $\phi$ must be
real greater than two.
Multiplying this matrix by the diagonal matrix above, we obtain
$\left(%
\begin{array}{cc}
  0 & -r^{-1}\e^{ml_X/2} \\
  r\e^{-ml_X/2} & \phi\e^{-ml_X/2} \\
\end{array}%
\right),$ so for sufficiently large
$m$, the resulting product ceases to be hyperbolic,
which is absurd.}
In order to have a well-defined projective limit, we shall ``kill''
the factor $m$ in a consistent way, and this can be achieved by dividing the result by
any graph length function of the curve since
${\rm g.l.}(G_X^mG_Z)=m\cdot {\rm g.l.}(G_X)+{\rm g.l.}(G_Z)$.

In the quantum case, however, the situation is much more involved. Indeed, let us consider
an example of the product, which is of form $U^mV$, as in (\ref{7-1}), where $U=\e^{X/2}$,
$V=\e^{Y/2}$, and $[X,Y]=4\pi i\hbar$. Thus, $(U^m)V=\e^{mX/2+Y/2+m\pi i\hbar/2}$
while $V(U^m)=\e^{mX/2+Y/2-m\pi i\hbar/2}$, and the corresponding logarithms do not
coincide as $m\to\infty$. Moreover, even the Hermiticity condition does not often
suffice to determine the proper length. For instance, given operators $U=\e^{\alpha X^2}$
and $V=\e^{i\beta \partial_X}$, we may calculate that $V(U^m)V=\e^{m\alpha X^2 +2i\beta \partial_X
+m\alpha\beta/3}$, i.e., the proper length in this case is $\alpha X^2+\alpha\beta/3$
and depends on the parameter~$\beta$ (of course, this correction is
purely quantum). This illustrates that proving continuity for the boundary transitions in the quantum
case requires more subtle estimates, which we perform in the next section after deriving recurrence
relations for the operators of quantum approximating multicurves.

We may now formulate our main result on quantizing Thurston theory:

\begin{theorem}\label{main}
Fix a spine $\Gamma$ of $F_1^1$ with corresponding
freeway
$\tau$.  Fix any projectivized vector $P{\vec \zeta}$
of foliation-shear coordinates on $\tau$ and any graph length function ${\rm g.l.}$.  For
any approximating sequence ${\vec n}^\beta$ to $P\vec \zeta$,
the limit
\be
\label{mainlimit}
\lim_{\beta\to\infty}\frac{{\rm p.l.}(\hat C_{\vec n^\beta})}{{\rm g.l.}(\hat C_{\vec n^\beta})}
\ee
exists both in the classical case as a real number and in the quantum case as a
weak operatorial limit.
\end{theorem}

Because both the numerator and
denominator in the limit are additive, this limit is projectively invariant and defines
a continuous function (in the classical case) or a weakly continuous family of operators
(in the quantum case) on the circle ${\cal PL}_0(F_{1,r}^s)$, for $r+s=1$.

The proof of the previous theorem occupies the remainder of this section.  The
continued fraction structure intrinsic to the torus case is used extensively, and various analogous
operatorial recursions are derived and studied.
There is a second essentially combinatorial proof of this result, however, only in the classical case since
we have no means to control the quantum ordering of the procedure.  Nevertheless, the structures discovered
are interesting, and we present this second proof in Appendix~A, which depends upon the
recursion (\ref{cf7}) derived later in Lemma~\ref{qlem2}. Indeed, this basic recursion arises from
``Rauzy--Veech--Zorich
induction'' \cite{Rauzy} in the special case of the torus, which is derived from first principles in the next section.

\subsection{Elements of the proof}
\subsubsection{Continued fraction expansion}\label{cfesec}

In each one-simplex in ${\cal PL}_0(F)$, illustrated as arrows in Figure~13, the approximating
multicurve is determined by two nonnegative integers, $m_1$ and $m_2$,
where we assume that $m_1>m_2$ with $m_1$ and $m_2$ relatively prime.  It is
convenient to represent the ratio $m_2/m_1$ as a simple continued fraction:
\be
\label{cf1}
m_2/m_1=\frac{1}{a_1+{\displaystyle\frac{1}{
  \renewcommand{\arraystretch}{0}
  \arraycolsep=0.05em
  \begin{array}{ccc}
   \strut a_2+{} & & \\
     & \ddots & \\
     & &  {}+{\displaystyle\frac{1}{a_{n-1}+{\displaystyle\frac1{a_n}}}}
     \end{array}
     }}}.
\ee

We concentrate on the case $m_X=m_2$, $m_Y=m_1$, $m_Z=m_X+m_Y$,
with the other cases following by symmetry,
and describe the recurrence procedure for constructing the corresponding approximating multicurves
(i.e., approximating geodesics, since any multicurve in the torus is just a multiple of a single
geodesic).

Referring to Figure~13, it is convenient to visualize this case by drawing a line
\be
\label{cf2}
\begin{array}{ccccccccc}
{\setlength{\unitlength}{1.5mm}%
\begin{picture}(6,4)(-3,-2)
\put(0,0){\makebox(0,0)[cc]{$\scriptstyle2$}}
\put(0,0){\circle{3}}
\end{picture}
}
&&&&&
{\setlength{\unitlength}{1.5mm}%
\begin{picture}(6,4)(-3,-2)
\put(0,0){\makebox(0,0)[cc]{$\scriptstyle1$}}
\put(0,0){\circle{3}}
\end{picture}
}
&&&
\\
\quad\bullet\quad&\quad\bullet\quad&\cdots&\bullet&\ \bullet\
&\circ&\cdots&\circ&\circ\\
1&2&&m_1-1&m_1&m_1+1&&m_1+m_2-1&m_1+m_2
\end{array},
\ee
which represents (the right side of) the cut over the edge $Z$
(homotopic to the cut $h_Z$ in Figure~12).
The geodesic function is constructed as follows.
Considering the multicurve in a neighborhood of the edge $Z$,
we start from the lowest thread that goes from the edge $Z$ to
the edge $X$ (labelled by the circled number one in Figure~12 and in (\ref{cf2}))
and moves to the right (this corresponds to the leftmost $\circ$
in (\ref{cf2})). We have the matrix $G_Y$ corresponding to passing
consecutively through edges $Z$ and $X$ and come back to the line
(\ref{cf2}), which must be periodically continued, at the point
that is situated to the right by $m_2$ points (in this case, just
the leftmost $\bullet$ in (\ref{cf2})
marked there by a circled number two). Each time passing the
bullet sign, we must set the matrix $G_X$ and passing the
$\circ$ sign we must set the matrix $G_Y$.

We thus move along the periodically continued line of circles and
bullet signs in (\ref{cf2}), at each step jumping $m_2$ points to the right
and setting the corresponding matrices $G_X$ or $G_Y$.  We
describe a recurrence procedure that produces from the continued fraction expansion
(\ref{cf1}) the correct sequence of operators, i.e., the correct
ordered sequence of branches traversed by the corresponding curve.  We
shall split the construction procedure into stages. Let us consider
the trajectory of the starting thread (the point $m_1+1$ in
(\ref{cf2})).  Each stage
terminates as soon as we come {closer} to the starting
thread~1 than the distance to it at the beginning of the stage
(approaching it from the opposite side compared to the beginning of the stage, see (\ref{cf3}) below).
We shall let $L_i$ denote the string of matrices for the first $i$
stages, where by convention we arrange matrices from right to left.
We terminate the first stage at the thread that is closer
than $m_2$ to the left of the starting point, i.e.,
$L_1=(G_X)^{a_1}G_Y$.

Define $\tilde
L_0\equiv G_X$ and introduce $\tilde L_i$,
for $i\ge1$, which is the matrix $L_i$ (composed from the
elementary matrices $G_X$ and $G_Y$)
in which the first two symbols of elementary
matrices must be interchanged: if the first symbol of $L_i$ is
$G_Y$ and the second is $G_X$, then the first symbol of $\tilde L_i$
is $G_X$ and the second is $G_Y$; all other elementary matrices
retain their forms. We illustrate this procedure in
Figure~14.\footnote{Note that we terminate a stage whenever we come closer
to the starting thread; this occurs at each stage on the opposite side from the previous stage, just as
for continued fractions.  The
appearance of tilded quantities is explained as follow: each time during the recursion when
we start a string ``parallel'' to some $L_k$ from the right to the starting thread
(in the circled domain), we must follow the same string of branches because, by
definition, there are no threads in the string $L_k$ except the very
first thread that appears at a distance closer than $|L_k|$ to
the starting thread. On the other hand, if we start from the left of the starting thread,
then we must interchange exactly  the first two appearances of matrices in the
resulting string.}
\begin{equation}\label{cf3}
\setlength{\unitlength}{.8mm}%
\begin{picture}(50,75)(30,5)
\thicklines
\dottedline[$\bullet$]{2.5}(60,62)(10,62)
\dottedline[$\circ$]{2.5}(60,62)(90,62)
\put(59,67){\vector(-1, 0){ 48}}
\put(35,72){\makebox(0,0)[cb]{($L_{i-2}$)}}
\put(61,67){\vector(1, 0){ 28}}
\put(75,72){\makebox(0,0)[cb]{$L_{i-1}$}}
\thinlines
\put(59,52){\vector(-1, 0){ 28}}
\put(59,52){\vector(1,0){0}}
\put(45,54){\makebox(0,0)[cb]{$\scriptstyle |L_{i-1}|$}}
\thinlines
\put(105,47){\vector(-1,1){12}}
\put(105,47){\line(1,0){10}}
\put(117,47){\makebox(0,0)[lc]{\scriptsize end\ of\ $(i-1)$th\ stage}}
\thicklines
\dottedline[$\bullet$]{2.5}(60,32)(35,32)
\dottedline[$\circ$]{2.5}(60,32)(90,32)
\dottedline{1.5}(30,32)(30,62)
\put(59,37){\vector(-1, 0){ 23}}
\put(48,42){\makebox(0,0)[cb]{$L_{i}$}}
\thinlines
\put(61,22){\vector(1, 0){ 23}}
\put(61,22){\vector(-1,0){0}}
\put(72,24){\makebox(0,0)[cb]{$\scriptstyle |L_{i}|$}}
\thinlines
\put(20,17){\vector(1,1){12}}
\put(20,17){\line(-1,0){10}}
\put(8,17){\makebox(0,0)[rc]{\scriptsize end\ of\ $i$th\ stage}}
\thicklines
\dottedline[$\bullet$]{2.5}(60,2)(35,2)
\dottedline[$\circ$]{2.5}(60,2)(75,2)
\dottedline{1.5}(85,32)(85,2)
\put(61,7){\vector(1, 0){ 13}}
\thinlines
\put(90,15){\vector(-1,-1){12}}
\put(90,15){\line(1,0){10}}
\put(102,15){\makebox(0,0)[lc]{\scriptsize end\ of\ $(i+1)$th\ stage}}
\put(68,12){\makebox(0,0)[cb]{$L_{i+1}$}}
\thinlines
\end{picture}
\end{equation}

\vspace{10pt}

\centerline{\bf Figure 14--threads a curve}

One thus sees directly that
$$
  \begin{array}{c}
L_2=(L_1)^{a_2-1}{\tilde L_0}L_1\\
    L_3=(\tilde L_2)^{a_3-1}{L_1}L_2\\
    L_4=(L_3)^{a_4-1}{\tilde L_2}L_3\\
    L_5=(\tilde L_4)^{a_5-1}{L_3}L_4\\
    \vdots\\
  \end{array}
$$
and this leads to the following recurrence relation.

\begin{lemma}\label{qlem2}Given the simple continued fraction expansion (\ref{cf1}) of $m_2/m_1$,
the sequence of matrices $L_n$ associated to the corresponding geodesic is given by the following
recursion
\begin{equation}\label{cf7}
\begin{array}{l}
 L_{2i}=({L_{2i-1}})^{a_{2i}-1}{\tilde L_{2i-2}}L_{2i-1}\\
 L_{2i+1}=(\tilde L_{2i})^{a_{2i+1}-1}{L_{2i-1}}L_{2i}
 \end{array}
 \quad\hbox{for}\quad i\ge1,\quad \tilde L_0=G_X,\quad
 L_1=(G_X)^{a_1}G_Y.
\end{equation}
\end{lemma}

Turning to the quantum case, we first show that the proper length operator {\rm(\ref{proj-l})}
{\rm p.l.}($\ORD{ L_{2i}}$)/{\rm g.l.}$(L_{2i})$ must agree with
{\rm p.l.}($\ORD{ L_{2i-1}}$)/{\rm g.l.}$(L_{2i-1})$ for $a_{2i}\to\infty$, that is,
the operators corresponding to continued fractions of form (\ref{cf1}) with large coefficient
$a_{n}$ must converge to the operator corresponding to the continued fraction terminated at
the $(n-1)$st step. To this end, we must analyze the structure of matrix products and
corresponding operators.  Indeed, the operators $L_i$ and $\tilde L_i$ enjoy elegant
commutation relations as we shall next see.

Notice that for every stage $i$
we have a geodesic corresponding to the matrix $L_i$
(because we can close the corresponding geodesic line without
self-intersections), and we can therefore define
the corresponding QMC
$$
{\cal L}_i\equiv \ORD{\tr L_i}.
$$
The first
observation pertains to $L_i$ and $\tilde L_i$: One of the corresponding curves can be
obtained from the other by a parallel shift along the cut $h_Z$ illustrated in Figure~12.
Thus, the curves are disjoint and hence are homotopic on the torus, and so
$\tr L_i=\tr \tilde L_i$, i.e.,
\begin{equation}\label{cf8}
{\cal L}_i=\tilde{\cal L}_i,~{\rm for}~i\geq 1.
\end{equation}
Furthermore, the curves $L_i$ and $L_{i+1}$ can be perturbed to have exactly one
intersection as one sees by considering how the corresponding
geodesic curves pass through the cut $h_Z$.\footnote{One of these
curves necessarily has odd subscript $2k+1$, and we can make a small
shift of all threads of this curve to the right from the threads
of the second curve $L_{2k}$ (or $L_{2k+2}$, depending on the
situation). We find that
there are no intersections of threads outside the region between
terminating points of $L_{2k+1}$ and $L_{2k}$ (or $L_{2k+2})$, and
in this domain, when closing the curves, there is produced exactly one
intersection point.} In the case where $i=2k$,
we obtain (in the notation of (\ref{so3}))
\bea\label{cf9}
[{\cal L}_{2k},{\cal L}_{2k+1}]_q=\xi \ORD{\tr L_{2k}L_{2k+1}},~{\rm and}~
[{\cal L}_{2k+1},{\cal L}_{2k}]_q=\xi \ORD{\tr L_{2k-1}(\tilde L_{2k})^{a_{2k+1}-1}},
\eea
where the proper quantum ordering is assumed for the terms in the right-hand
sides.

Formulas (\ref{cf9}) are crucial when proving the continuity. Letting
$$
I_m\equiv \ORD{\tr (L_{2i-1})^{m-1}\tilde L_{2i-2}L_{2i-1}},
$$
$${\cal L}_{2i-1}\equiv \e^{\ell_X/2}+\e^{-\ell_X/2},$$ we find
\bea
\label{cont1}
I_{m-1}{\cal L}_{2i-1}&=&q^{1/2}I_m+q^{-1/2}I_{m-2},\\
\label{cont2}
{\cal L}_{2i-1}I_{m-1}&=&q^{-1/2}I_m+q^{1/2}I_{m-2},
\eea
where the basis of the recursion is given by
$I_0\equiv \e^{l_Y/2}+\e^{-l_Y/2}$ and
$I_{-1}\equiv \e^{l_Z/2}+\e^{-l_Z/2}$.
>From (\ref{cont1}), we have the exact equalities
\be
\label{cont3}
I_m=q^{-m/2}I_0\bigl(\e^{m\ell_X/2}+\e^{-m\ell_X/2}\bigr)-q^{-m/2-1/2}I_{-1}
\bigl(\e^{(m-1)\ell_X/2}+\e^{-(m-1)\ell_X/2}\bigr),
\ee
or, equivalently,
\be
\label{cont4}
I_m=q^{m/2}\bigl(\e^{m\ell_X/2}+\e^{-m\ell_X/2}\bigr)I_0-q^{m/2+1/2}
\bigl(\e^{(m-1)\ell_X/2}+\e^{-(m-1)\ell_X/2}\bigr)I_{-1}.
\ee

We wish to present the expression (\ref{cont3}) or (\ref{cont4}) in the form
$\e^{mH_1+H_0}$, where $H_1$, $H_0$ are Hermitian operators independent of $m$, so the proper
length is then just $H_1$ (while $H_0$ introduces quantum corrections that
do not affect the proper length limit but are
important for ensuring the proper commutation relations).

Let us choose a compact domain ${\cal F}$ in the function space $L^2({\Bbb R})$
such that norms of all the operators in play are bounded for functions from
this domain. Notice that neither term on the right-hand side in (\ref{cont3}) is self-adjoint, but
if one of these terms prevails (and it can be only the first term since we have an expression with
coefficients that are all {positive} in the classical limit of the right-hand side, and hence
we must have a positive left-hand side as well), then we can replace the total sum by this
prevailing term and obtain an approximate equality
$$
q^{-m/2}I_0\bigl(\e^{m\ell_X/2}+\e^{-m\ell_X/2}\bigr)\sim
q^{m/2}\bigl(\e^{m\ell_X/2}+\e^{-m\ell_X/2}\bigr)I_0
$$
in this limit.
By considering the spectral expansion with respect to the eigenfunctions of the operator $\ell_X$,
we immediately conclude that $I_0\sim \e^{- 2\pi i\hbar
{\partial}/{\partial |\ell_X|}}$ in this limit, and then
\be
\label{cont6}
H_1=|\ell_X|/2,
\ee
where the modulus has to be understood in terms of the spectral expansion: having a QMC
operator ${\cal L}_X$ which admits a spectral decomposition (see formulas
(\ref{cf12})--(\ref{cf14}) below) in functions $|\alpha_S\rangle$, we define the
operator $|\ell_X|$ by its action on these functions:
$|\ell_X|\,|\alpha_S\rangle=|S|\,|\alpha_S\rangle$.

A potentially problematic situation is when neither of the terms prevails and their
difference remains finite as $m\to\infty$. This would correspond, as in the classical
case discussed before, to a situation where the corresponding element fails to be
hyperbolic.  As we next show, we must obtain ``long'' curves when $m$ goes to infinity, so this is
also impossible in the quantum case; that is, for $G_X$, $G_Y$, $G_Z$, and ${\wtd G}_Z$ from
(\ref{torus-classic}), (\ref{tildeG}) and representing $G_X=\e^{\ell_X/2}+\e^{-\ell_X/2}$,
we must prove the operatorial inequality
\be
\label{cont7}
G_Y\e^{|\ell_X|/2}>q^{1/2}{\wtd G}_Z.
\ee

To prove this, first express $\e^{\ell_X/2}$ through $G_X$. Taking the positive branch of the square
root, we find from the left-hand side of (\ref{cont7}) the expression
$$
G_Y\left(\frac{G_X}{2}+\sqrt{\frac{G_X^2}{4}-1}\right),
$$
and since $G_YG_X=q^{1/2}{\wtd G}_Z+q^{-1/2}G_Z$, we must compare two expressions
$G_Y\sqrt{\frac{G_X^2}{4}-1}$ and $\bigl(q^{1/2}{\wtd G}_Z-q^{-1/2}G_Z\bigr)/2$.
Multiplying by the Hermitian conjugate on the right in both expressions, we eliminate
the square root and arrive at Laurent polynomial expressions. After some simple
algebra, we come to the inequality to be proved:
$$
\frac12\bigl(q{\wtd G}_ZG_Z+q^{-1}G_Z{\wtd G}_Z\bigr)\quad>\quad G_Y^2.
$$
To see this, we expand the left-hand side
\bea
&&\bigl(\e^{-X-Z}+\e^{X-Z}+\e^{X+Z}+(q+q^{-1})\e^X\bigr)\nonumber\\
&+&\frac{q+q^{-1}}{2}\bigl(\e^{-X-Y-Z}+\e^{X+Y+Z}+\e^Z+2\e^{-Z}+(q+q^{-1})\bigr)
\nonumber\\
&+&\frac{q^2+q^{-2}}{2}\bigl(\e^{Y+Z}+\e^{-Y-Z}+\e^{-Y+Z}+(q+q^{-1})\e^{-Y}\bigr)
+\frac{q^3+q^{-3}}{2}\e^Z,
\nonumber
\eea
while the right-hand side is expressed as
$$
G_Y^2=\e^{-X-Z}+\e^{X-Z}+\e^{X+Z}+2+(q+q^{-1})(\e^X+\e^{-Z}).
$$
Subtracting this expression from the previous one, we obtain that this difference is
\bea
\frac{q^2+q^{-2}}{2}G_X^2+\frac{q+q^{-1}}{2}\bigl[\e^{X+Y+Z}+\e^{-X-Y-Z}-q-q^{-1}\bigr],
\eea
and {both} these terms are positive definite for $|q|=1$.

The operatorial inequality (\ref{cont7})
has therefore been established.  This proves that the limit (\ref{mainlimit}) exists and is
well defined at rational points of the continued fraction expansion, and we next turn to the
case of infinite continued fraction expansions, i.e., infinite sequences of elementary operators.

\subsubsection{Mapping class group transformations and the unzipping procedure}\label{mcg-sect}

We consider now an infinite continued fraction expansion
$a_1,a_2,\ldots ,a_n,a_{n+1},\ldots$ extending the notation of (\ref{cf1}).
As we shall see, there is a corresponding sequence of unzippings of the freeway $\tau$ associated to
a spine $\Gamma$ of $F$ as in Section 4.5 as well as an associated sequence of mapping class group
elements, expressed as Dehn twists, which reduce an approximating multicurve to one of two possible
graph simple curves.

Given the recursive representation (\ref{cf7}) for the operator of a geodesic curve determined by
a continued fraction expansion (\ref{cf1}) and applying two (unitary) operators $D_X$ and $D_Y$
of the modular transformations of the form (\ref{corr1}) that correspond to the
respective Dehn twists along the corresponding closed curves $\gamma_X$ and $\gamma_{Y}$
(with the respective geodesic functions $G_X$ and $G_Y$),
we shall construct the sequence
of zipping or unzipping transformations.

\begin{defin}\label{def5-1}
{\rm
An approximating multicurve is determined by two nonnegative integers $m_1$ and $m_2$.
As in Figure~13, enumerate such a pair as a triple $(m_1,m_2,m_1+m_2)$. If $m_1>m_2$,
the action of the Dehn twist $D_{Y}^{-1}$
along $\gamma_{Y}$ (an {\em unzipping\/}
transformation---the Dehn twist in the opposite direction)
is $$D_{Y}^{-1} :(m_1,m_2,m_1+m_2)\mapsto (m_1-m_2,m_2,m_1)$$
while if $m_1<m_2$, we apply
the unzipping transformation along the curve $\gamma_{X}$, which gives
$$D_{X}^{-1} :(m_1,m_2,m_1+m_2)\mapsto (m_1,m_2-m_1,m_2).
\footnote{In terms of the symbolic dynamics
of elementary operators $G_X$ and $G_Y$, we can present this action as follows:
$D_{X}^{-1}(G_X)=G_X$, $D_{X}^{-1}(G_XG_Y)=G_Y$ and $D_{Y}^{-1}(G_XG_Y)=G_X$, $D_{Y}^{-1}(G_Y)=G_Y$.}$$
}
\end{defin}

Given a continued fraction expansion
(\ref{cf1}), we construct the sequence of unzipping transformations
\begin{equation}\label{cf10}
D_{(X\ \hbox{\scriptsize or}\ Y)}^{-a_n}
D_{(Y\ \hbox{\scriptsize or}\ X)}^{-a_{n-1}}
\cdots D_Y^{-a_3}D_X^{-a_2}D_Y^{-a_1},
\end{equation}
which, when applied to the approximating multicurve $(m_1,m_2,m_1+m_2)$, reduces
it either to $(1,0,1)\equiv \gamma_Y$ for $n$ even or to $(0,1,1)\equiv \gamma_X$
for $n$ odd.

\begin{defin}\label{def5-2}
{\rm
Equivalently, we can consider the {\em zipping\/} procedure, that is, given
a sequence of
Dehn twists $D_Y^{a_1}D_X^{a_2}D_Y^{a_3}\cdots D_Y^{a_{n-1}}D_X^{a_n}$ applied
to the curve $\gamma_Y$, we obtain the curve $(m_1,m_2,m_1+m_2)$.
}
\end{defin}

Considering the sequence (\ref{cf10}) of quantum Dehn twist operators (\ref{corr1})
and exploiting the quantum invariance from Lemma~\ref{qlem1}, we come to the { main
observation} that having an involved expression
for the proper limit (\ref{proj-l}) of a
QMC operator constructed by the rules described in Lemmas~\ref{qlem1}
and~\ref{qlem2} in terms of the elementary operators $X,Y,Z$, we may
perform the sequence (\ref{cf10}) of unzipping quantum modular transformations,
which reduces this operator to a standard form of the quantum operator $G_Y$ or $G_X$
expressed through the new operators $X^{(n)},Y^{(n)},Z^{(n)}$ related to the
initial operators by this sequence of quantum modular transformations.  This is the operatorial
statement of naturality of lengths under the mapping class group action.

It is intuitively natural to imagine that as the geodesic lengths must diverge as
$m_1$, $m_2$ tend to infinity, we must eventually come to an asymptotic regime
where all quantities $X^{(n)}$, $Y^{(n)}$ are large in the literal or operatorial sense
for all sufficiently large $n$.  The inexorability of the approach to this asymptotic regime is
not obvious and is described in the next section.

\subsubsection{Asymptotic regime}

Let us recall the modular transformations for $X$, $Y$, and $Z$ variables:
\begin{equation}\label{DX1}
    D^{-1}_X\:(X,Y,Z)\mapsto (X+2\phi^\hbar(Z),-Z,Y-2\phi^\hbar(-Z))
\end{equation}
and
\begin{equation}\label{DY1}
    D^{-1}_Y\:(X,Y,Z)\mapsto (-Z,Y-2\phi^\hbar(-Z),X+2\phi^\hbar(Z)).
\end{equation}
In terms of the quantities $U\equiv \e^{X/2}$ and $V\equiv\e^{-Y/2}$
in the case where $X+Y+Z=0$, we have
\begin{eqnarray}
  D_X^{-1}\left(
\begin{array}{c}
  U \\
  V \\
\end{array}
  \right)D_X &=& \left(\begin{array}{c}
    \e^{X/2}+\e^{-X/2-Y}\\
    \e^{-Y/2+X/2}        \\
  \end{array}\right)\equiv
  \left(\begin{array}{c}
    U+VU^{-1}V \\
    q^{1/2}U^{-1}V \\
  \end{array}\right),
\label{UVas1}
  \\
  D_Y^{-1}\left(
\begin{array}{c}
  U \\
  V \\
\end{array}
  \right)D_Y &=&
  \left(\begin{array}{c}
    \e^{X/2+Y/2} \\
     \e^{-Y/2}+\e^{Y/2+X}\\
  \end{array}\right)\equiv
  \left(\begin{array}{c}
    q^{1/2}UV^{-1} \\
    V+UV^{-1}U \\
  \end{array}\right).
\label{UVas2}
\end{eqnarray}
Worth mentioning is that since each operator $Z_\alpha$ is Hermitian,
each exponential is positive definite, so we can always write, for instance,
that $U+VU^{-1}V>U$ in the sense of spectral expansion:
$\langle f|U+VU^{-1}V|f\rangle>\langle f|U|f\rangle$ for any function $f\in L^2({\Bbb R})$.

Using now an alternating sequence of transformations (\ref{UVas1}),
(\ref{UVas2}), we shall subsequently show that we attain the asymptotic regime of
large positive $X$ (large $U$) and large in absolute value
negative $Y$ (large $V$) starting
from every pair of $X$ and $Y$ lying in a compact domain of the
$(X,Y)$-plane in the classical case or acting within a compactum
of test functions with bounded derivatives in the function space in
the quantum operatorial case.

In this section, we verify that the asymptotic regime is attained
for the distinguished sequence of modular
transformations corresponding to the Fibonacci number sequence (golden mean),
namely, for alternating $D_X^{-1}$ and $D_Y^{-1}$.
The proof in the
general case is analogous although the structure is more involved, as described
in the next section.

Given the sequence of
transformations $D_X^{-1}D_Y^{-1}\cdots
D_X^{-1}D_Y^{-1}=\bigr(D_X^{-1}D_Y^{-1}\bigl)^n$, we obtain
\begin{eqnarray}
&&D_X^{-1}D_Y^{-1}\left(
\begin{array}{c}
  U \\
  V \\
\end{array}
  \right) =
    \label{cf11}
  \\
  &&=
  \left(
\begin{array}{c}
  UV^{-1}U+V \\
  V^{1/2}(V^{-1}Uq^{-1/2}+U^{-1}Vq^{1/2})U(V^{-1}Uq^{-1/2}+U^{-1}Vq^{1/2})V^{1/2}+V^{1/2}U^{-1}V^{1/2} \\
\end{array}
  \right).
  \nonumber
\end{eqnarray}
In the classical case, the asymptotics is already clear from this formula; for the first
entry in (\ref{cf11}), we have $UV^{-1}U+V=U(V^{-1}U+U^{-1}V)>2U$ as the expression
in the parentheses has the form $\e^S+\e^{-S}\geq 2$ for any real $S$. The same logic applies to the
second entry in (\ref{cf11}), and we deduce that the classical part $V(U+U^{-1})$
has the same property.

The proof given below in the quantum case is more subtle as it needs a
thorough operatorial analysis. Nevertheless,
the estimates turn out to be close to those in the classical case, which we briefly discuss here: we
must prove that a lower bound on the operatorial spectrum on a compactum in the function space
diverges with $n$.  This is a routine procedure,
which uses that the action of operators $X$ and $Y$ in the basis of, say,
normalized Hermitian functions $h_n$ has the form
\bea
X|h_n\rangle&=&\sqrt{4\pi\hbar}\bigl(\sqrt{n+1}|h_{n+1}\rangle+\sqrt{n}|h_{n-1}\rangle\bigr)
\nonumber
\\
Y|h_n\rangle&=&\sqrt{4\pi\hbar}\frac{i}{2}\bigl(\sqrt{n+1}|h_{n+1}\rangle-\sqrt{n}|h_{n-1}\rangle\bigr).
\nonumber
\eea
These operators ``almost'' commute in the domain of large $n$, which allows
the combinatorics to be analyzed semiclassically.

In the quantum case, we recall the construction of quantum
Dehn twists and their eigenfunctions from
\cite{Kashaev3}. The generator of the Dehn twist $D_X$ has the
form
\begin{equation}\label{cf12}
  D_X=\e^{q_1^2/2\pi i\hbar}F^\hbar(q_1+p_1),\quad q_1=X/2,\
  p_1=2\pi i\hbar\partial_X
\end{equation}
and because it commutes with the geodesic length operator $G_X$,
they share the common set of eigenfunctions
\begin{equation}\label{cf13}
  |\alpha_S\rangle=\e^{-X^2/16\pi i \hbar}F^\hbar(S+X)F^\hbar(-S+X)
\end{equation}
with the eigenvalues
\begin{equation}\label{cf14}
  G_X|\alpha_S\rangle=2\cosh(S/2)|\alpha_S\rangle;\qquad
  D_X|\alpha_S\rangle=\e^{S^2/2\pi i\hbar}|\alpha_S\rangle.
\end{equation}
The functions $|\alpha_S\rangle$ constitute a complete set of
functions in the sense that
$$
\langle\alpha_T|\alpha_S\rangle=\delta(S-T)\nu^{-1}(S),\qquad
\nu(S)=4\sinh(\pi S)\sinh(\pi\hbar S)
$$
and
$$
\int_0^\infty\nu(S)dS\,|\alpha_S\rangle\langle\alpha_S|=\hbox{Id}
$$

We now split the plane of the variables $(X,Y)=(X^{(0)},Y^{(0)})$ into four sub-domains and consider
the action of the Dehn twists $D_X^{-1}$ and $D_Y^{-1}$ in each sub-domain.

\vskip .2in

\leftskip .4in

\noindent {\sl Domain I.}\
$\{X^{(0)}>0,\ Y^{(0)}>0\}\cup\{X^{(0)}>0,\ Y^{(0)}<0\ \hbox{and}\ |X^{(0)}|>|Y^{(0)}|\}$.
\vskip -.3in
\bea
&&D_X^{-1}\:(X^{(0)},Y^{(0)})=(X^{(1)},Y^{(1)})=(X^{(0)},X^{(0)}+Y^{(0)})\in \hbox{\,Domain\
I};
\nonumber
\\
&&D_Y^{-1}\:(X^{(0)},Y^{(0)})=(X^{(1)},Y^{(1)})=(X^{(0)}+Y^{(0)},-2X^{(0)}-Y^{(0)})\in
\hbox{\,Domain\ II};
\nonumber
\eea
\noindent {\sl Domain II.}\ $\{X^{(0)}>0,\ Y^{(0)}<0\ \hbox{and}\ |X^{(0)}|<|Y^{(0)}|\}\cup
\{X^{(0)}<0,\ Y^{(0)}<0\}$.
\vskip -.3in
\bea
&&D_X^{-1}\:(X^{(0)},Y^{(0)})=(X^{(1)},Y^{(1)})=(-X^{(0)}-2Y^{(0)},X^{(0)}+Y^{(0)})\in
\hbox{\,Domain\ I};
\nonumber
\\
&&D_Y^{-1}\:(X^{(0)},Y^{(0)})=(X^{(1)},Y^{(1)})=(X^{(0)}+Y^{(0)},Y^{(0)})\in \hbox{\,Domain\
II};
\nonumber
\eea
\noindent {\sl Domain IIIa.}\ $\{X^{(0)}<0,\ Y^{(0)}>0\ \hbox{and}\ |X^{(0)}|<|Y^{(0)}|\}$.
\vskip -.3in
\bea
&&D_X^{-1}\:(X^{(0)},Y^{(0)})=(X^{(1)},Y^{(1)})=(X^{(0)},X^{(0)}+Y^{(0)})\in \hbox{\,Domain\
IIIa or IIIb};
\nonumber
\\
&&D_Y^{-1}\:(X^{(0)},Y^{(0)})=(X^{(1)},Y^{(1)})=(X^{(0)}+Y^{(0)},-2X^{(0)}-Y^{(0)})\in
\hbox{\,Domain\ I or II};
\nonumber
\eea
\noindent {\sl Domain IIIb.}\ $\{X^{(0)}<0,\ Y^{(0)}>0\ \hbox{and}\ |X^{(0)}|>|Y^{(0)}|)$.
\vskip -.3in
\bea
&&D_X^{-1}\:(X^{(0)},Y^{(0)})=(X^{(1)},Y^{(1)})=(-X^{(0)}-2Y^{(0)},X^{(0)}+Y^{(0)})\in
\hbox{\,Domain\ I or II};
\nonumber
\\
&&D_Y^{-1}\:(X^{(0)},Y^{(0)})=(X^{(1)},Y^{(1)})=(X^{(0)}+Y^{(0)},Y^{(0)})\in \hbox{\,Domain\
IIIa or IIIb}.
\nonumber
\eea

\vskip .2in

\leftskip=0ex

We see that only Domain III is potentially problematic. This regime is however unstable:
absolute values of $X$ and $Y$ variables decrease in this regime
and they eventually leave the asymptotic regime as
soon as we {\em remain\/} in Domain III; immediately upon leaving this domain, we come to domains I and II and
will never leave this three quarters of the $(X,Y)$-plane. The above considerations of $U$, $V$ just
demonstrate that even if we were initially in Domain III, we come to the nonasymptotic domain of
bounded $X$ and $Y$ and then will leave this compactum moving toward asymptotic expansions in
domains I and II.

The asymptotic dynamics {always} takes place in the first three quarters of
the $(X,Y)$-plane.
Nevertheless, even this dynamics is rather involved. The stable regime corresponds to the case where
we are in Domain II {before} applying one or several operators $D_X^{-1}$. The
application of the first of these operators brings us to Domain I,
and upon subsequent applications of the operators $D_X^{-1}$ we remain in Domain I.
Next, if we were in Domain I, then the very first application of the operator $D_Y^{-1}$ brings us
to Domain II, and we then remain in Domain II upon subsequent applications of $D_Y^{-1}$.

We turn now to actual geodesic lengths of curves or proper lengths of operators.
If a sequence of unzipping transformations terminates,
this means that we have a graph simple geodesic, which is either $G_X$ if the
{ last} transformation was $D_Y^{-1}$ or $G_Y$ if the last transformation
was $D_X^{-1}$. Considering the corresponding geodesic or proper lengths,
we find that up to exponentially small corrections, the leading contributions
in the above domains are
\begin{eqnarray}
  {\rm p.l.}(\gamma_Y) &=& X+Y/2 \quad \hbox{in domain I},
  \label{cf16a}
  \\
  {\rm p.l.}(\gamma_X) &=& -Y-X/2 \quad \hbox{in domain II}
  \label{cf16b}
\end{eqnarray}
(see expressions (\ref{torus-classic})).

Thus, although the transformation laws for the
variables $X,Y$ themselves do not possess the property of linearity with respect to the parameters
$a_i$,
$a_j$, when applying sequences of transformations $(D_X^{-1})^a_i\equiv D_X^{-a_i}$ or
$(D_Y^{-1})^{a_j}\equiv D_Y^{-a_j}$, the proper lengths do possess this property!
Namely, starting with variables $(X^{(0)}, Y^{(0)})$ lying in the corresponding domains and
applying the sequences of transformations $D_X^{-a_i}$ or $D_Y^{-a_j}$, we obtain for the
resulting proper lengths the following expressions:
\begin{eqnarray}
\label{cf15a}
 {\rm p.l.}(\gamma_{Y^{(i)}}) &=& -\frac{Y^{(0)}}2+a_i\left(-Y^{(0)}-\frac{X^{(0)}}{2}\right)
 \quad\hbox{for $X^{(0)},Y^{(0)}\in$ domain II}  \\
\label{cf15b}
 {\rm p.l.}(\gamma_{X^{(j)}}) &=& \frac{X^{(0)}}2+a_j\left(X^{(0)}+\frac{Y^{(0)}}{2}\right)
 \quad\hbox{for $X^{(0)},Y^{(0)}\in$ domain I}
\end{eqnarray}

Let us now explore the asymptotic formulas (\ref{cf16a}), (\ref{cf16b}) and (\ref{cf15a}), (\ref{cf15b})
first in the classical case to close this section, relegating the discussion of the quantum case
to the next section.

Assume that we start from the variables $(X^{(0)},Y^{(0)})$ in domain I
and have the corresponding
initial length ${\rm p.l.}(\gamma_{Y^{(0)}})$
from (\ref{cf16a}). Applying the transformation $D_Y^{-a_j}$, we obtain new variables
$(X^{(j)},Y^{(j)})$ and the new proper length ${\rm p.l.}(\gamma_{X^{(j)}})$
(\ref{cf15b}) having form (\ref{cf16b}) in these new variables, which must now lie in domain II.
Note that explicitly
$$
(X^{(j)},Y^{(j)})=(X^{(0)}+Y^{(0)}+(a_j-1)(-2X^{(0)}-Y^{(0)}),-2X^{(0)}-Y^{(0)}).
$$
We then apply the transformation $D_X^{-a_i}$
to obtain variables $(X^{(j,i)},Y^{(j,i)})$, and
the proper length ${\rm p.l.}(\gamma_{Y^{(j,i)}})$ is expressed as in (\ref{cf15a}), where the term
multiplied by $a_i$ is none other than ${\rm p.l.}(\gamma_{X^{(j)}})$ and the term $-\frac{Y^{(j)}}2$ is
exactly ${\rm p.l.}(\gamma_{Y^{(0)}})$. We thus find in the asymptotic regime that the corresponding
lengths are related by {\em exactly the same\/} recurrence relation as for a
graph length (the latter of which follows immediately from (\ref{cf7})):
\begin{eqnarray}
   {\rm p.l.}(\gamma_{Y^{(j,i)}})&=&a_i{\rm p.l.}(\gamma_{X^{(j)}})+{\rm p.l.}(\gamma_{Y^{(0)}}),
\label{cf17a}
   \\
   {\rm g.l.}(\gamma_{Y^{(j,i)}})&=&a_i{\rm g.l.}(\gamma_{X^{(j)}})+{\rm g.l.}(\gamma_{Y^{(0)}}).
\label{cf17b}
\end{eqnarray}

It is then easy to conclude that the {ratio} of these two quantities has a
definite limit as $i\to\infty$ for {any} sequence of numbers $a_i$.  It is a standard
estimate: given two numerical sequences
(\ref{cf17a}) and (\ref{cf17b}) and denoting the relative error of their ratio as
$\eps_i$, i.e., at the $i$th step, the ratio is $S(1+\eps_i)$, where $S$ is constant,
for $a_{i+1}>1$ at the ($i+1$)th step, we obtain $\eps_{i+1}<\eps_i/(a_{i+1}-1/2)$,
or if we have two coefficients $a_{i+1}=a_i=1$, then $\eps_{i+1}<\eps_{i-1}/1{.}5$. In general, for $\eps_i$
small enough, we always have $\eps_{i+1}<\eps_i$. This shows that the
relative error decreases exponentially with the index $i$.

\subsubsection{Quantum continued fraction expansion}\label{cfe2sec}

Let us turn again to the sequence (\ref{cf10}) of unzipping transformations.
In order to obtain operatorial expressions, we consider the unitary operators
$D_X$, $D_Y$ and explicitly indicate the
variables in which these operators are expressed, i.e., we write
$D_{X^{(j)}}\equiv D_X(X^{(j)},Y^{(j)})$ for the Dehn twist along $\gamma _X$ at the ($j+1$)th step. Thus,
\begin{equation}\label{qf1}
(X^{(n)},Y^{(n)})=
D_{X^{(n-1)}}^{-a_n}
D_{Y^{(n-2)}}^{-a_{n-1}}
\cdots D_{X^{(1)}}^{-a_2}D_{Y^{(0)}}^{-a_1}
(X^{(0)},Y^{(0)})
D_{Y^{(0)}}^{a_1}D_{X^{(1)}}^{a_2}\cdots D_{Y^{(n-2)}}^{a_{n-1}}D_{X^{(n-1)}}^{a_n},
\end{equation}
In order to represent such long strings of operators in terms of  the original
operators $(X^{(0)},Y^{(0)})$, we invert the dependence, i.e., we remember that, for instance,
$$
D_{X^{(1)}}^{-a_2} = D_{Y^{(0)}}^{-a_1}D_{X^{(0)}}^{-a_2}D_{Y^{(0)}}^{a_1},
$$
etc., which gives
\begin{equation}\label{qf2}
(X^{(n)},Y^{(n)})=
D_{Y^{(0)}}^{-a_1}D_{X^{(0)}}^{-a_2}\cdots D_{Y^{(0)}}^{-a_{n-1}}D_{X^{(0)}}^{-a_n},
(X^{(0)},Y^{(0)})
D_{X^{(0)}}^{a_n}D_{Y^{(0)}}^{a_{n-1}}\cdots D_{X^{(0)}}^{a_2}D_{Y^{(0)}}^{a_1}
\end{equation}

We shall compute with bases of functions that are convenient in the asymptotic regime. Let
\be
\label{qf*}
|f_{\mu,s}\rangle\equiv \e^{i\mu(x-s)^2/2}\quad \mu,s\in {\Bbb R}.
\ee
These functions constitute a basis at each $\mu$:
\be
\langle f_{\mu,t}|f_{\mu,s}\rangle=\frac{2\pi}\mu \delta(s-t),
\qquad
\int_{-\infty}^{\infty}ds|f_{\mu,s}\rangle\,\langle f_{\mu,s}|=\frac{2\pi}\mu\, \hbox{Id}.
\ee
For two arbitrary real numbers $w$ and $\gamma$, we have
\bea
\e^{iwx^2/2}|f_{\mu,s}\rangle&=&\e^{i\mu ws^2/2}\left|f_{\mu+w,\frac{s}{1+w/\mu}}\right\rangle,
\label{qf3a}
\\
\e^{i\gamma\partial_x^2/2}|f_{\mu,s}\rangle&=&\frac1{\sqrt{1+\gamma\mu}}\left|f_{\frac{1}{\gamma+1/\mu},s}\right\rangle,
\label{qf3b}
\eea
and
\bea
\label{qf4a}
\langle f_{\mu,s}|x|f_{\mu,t}\rangle&=&\frac{2\pi}{i\mu^2}\delta'(s-t)+\frac{2\pi}{\mu}s\delta(s-t),
\\
\label{qf4b}
\langle f_{\mu,s}|\frac1i\partial_x|f_{\mu,t}\rangle&=&\frac{2\pi}{i\mu}\delta'(s-t).
\eea
We now define the dimensionless variable~$x$ and set
\be
\label{qf5}
X|f_{\mu,s}\rangle=\sqrt{4\pi\hbar}x\cdot|f_{\mu,s}\rangle;\qquad Y|f_{\mu,s}\rangle=
\sqrt{4\pi\hbar}\frac1i \frac\partial{\partial x}\cdot|f_{\mu,s}\rangle.
\ee
The explicit formulas for the operators $D_X$ and $D_Y$ acting on $|f_{\mu,s}\rangle$
in the asymptotic regime are
\bea
\label{qf6a}
D_X^{a_i}|f_{\mu,s}\rangle=\e^{i(a_i-1)x^2/2}\,\e^{-i\partial_x^2}\,\e^{ix^2/2}|f_{\mu,s}\rangle,
\\
\label{qf6b}
D_Y^{a_j}|f_{\mu,s}\rangle=\e^{i(a_j-1)\partial_x^2/2}\,
\e^{-ix^2}\,\e^{i\partial_x^2/2}|f_{\mu,s}\rangle.
\eea
In order to establish the required recurrence relation, we must compare
matrix elements of the three consecutive length operators
in the corresponding operatorial decompositions:
\bea
A^{(0)}_{st}&=&\langle f_{\mu,s}|X+\frac{Y}2|f_{\mu,t}\rangle,
\nonumber
\\
A^{(j)}_{st}&=&\langle f_{\mu,s}D_Y^{-a_j}|-Y-\frac{X}2|D_Y^{a_j}f_{\mu,t}\rangle,
\nonumber
\\
A^{(j,i)}_{st}&=&\langle f_{\mu,s}D_Y^{-a_j}D_X^{-a_i}|X+\frac{Y}2|D_X^{a_i}D_Y^{a_j}f_{\mu,t}\rangle.
\nonumber
\eea
Now, using formulas (\ref{qf3a})--(\ref{qf6b}), it is straightforward to show that
\be
\label{qf7}
A^{(j,i)}_{st}=A^{(0)}_{st}+a_iA^{(j)}_{st}\quad ,
\ee
for all $s,t$, i.e., we again attain the recurrence relation (\ref{cf17a}) in the asymptotic regime,
but now for the matrix elements of the operators of the
quantum proper lengths. Estimates show that the corrections due to both
the (operatorial) deviations from the asymptotic regime and the error parameters
$\eps_i$ (as for (\ref{cf17a}), (\ref{cf17b})) decrease exponentially with the index~$i$, so the limit
(\ref{mainlimit}) exists in a weak operatorial sense.
We conclude that ratios (\ref{mainlimit})
define a weakly continuous family of operators parameterized by projective transverse
measures on the freeway associated to a spine of the
once-punctured torus. This completes the proof of Theorem~{\ref{main}}.

\newsection{Conclusion}

We hope to have added to the mathematical foundation and general understanding of the
quantization of Teichm\"uller space and its geometric underpinnings in the first several sections of
this paper.  We also hope that the survey given here of train tracks and their extensions
might be useful.

The quantization of Thurston's boundary in general seems to be a substantial project, which we have only just begun here
with the quantization of continued fractions.
First of all, one would like a better understanding of the operators we have constructed, for instance, an intrinsic
characterization or an explicit calculational framework for them.  At the same time, our current constructions depend
upon a choice of spine, and there would seem to be a more invariant version of the theory, where the
choice of spine is dictated by the combinatorics of the cell decomposition of Teichm\"uller space; the calculations in
this paper apply to each such spine (since there is a combinatorially unique cubic one) for the once-punctured torus.

Second of all, the quantization of Thurston's boundary for higher-genus or multiply-punctured surfaces
may be approachable using the improved quantum ordering.  Namely, in any fixed spine of the surface,
there is a fixed finite family of ``edge-simple'' closed edge-paths which by definition never twice traverse the same
oriented edge.  It is elementary to see that any closed edge-path on $\Gamma$ may be written
non-uniquely as a concatenation of edge-simple paths, where the particular concatenation depends
upon a starting point.  (Edge-simple paths were studied as ``canonical curves'' on train tracks in \cite{Penn0}; they
contain the extreme points of the polyhedron of projective measures on the track.)
It follows that an arbitrary leaf of a measured foliation carried by a freeway can be written as a concatenation of
paths from this finite collection of edge-simple paths.  The corresponding quantum operatorial statement
results from the improved quantum ordering described here.  Thus, whereas the quantization of the once-punctured torus
devolved, in effect, to an analysis of two-letter words, the quantization of Thurston's boundary sphere in general may depend
upon an analysis of words comprised of letters which are edge-simple paths.

One appealing long-term goal would be to
discover the Thurston classification already on the operatorial level, for instance, with the dilatation in the
pseudo-Anosov case explicitly computable from the MCG operator or from the invariant projective
foliation operator.

Another intriguing aspect involves generalizations of graph
length
functions insofar as the proof of Theorem~\ref{main} holds taking as
graph length any continuous positive definite function which is homogeneous of
degree one.  A natural choice of such a function is
induced by the geodesic
length of the corresponding geodesic curve taken for a fixed
basepoint in Teichm\"uller space on its fixed
spine, for instance,
vanishing shear coordinates on
the usual spine in the once-punctured torus.
What sort of regularity (e.g., piecewise smoothness) is achieved in the operators
corresponding to points of Thurston's boundary under such ``gauge
fixing''?

Also worth mentioning are very recent advances in the description of quantum $sl(n,{\Bbb R})$ connections
\cite{FockG}, where one finds an improved
quantum ordering in a more complicated higher-dimensional setting.

\setcounter{section}{0}
\appendix{Combinatorial proof of Theorem~\ref{main}}
In this appendix, we give a complementary, combinatorial proof of the classical Theorem~\ref{main}
using  the recurrence relation (\ref{cf7}).  At the present state of
understanding, the proof applies only to the classical case as we cannot control the quantum ordering.

Let us recall the structure of the matrix product (\ref{Pz}). It is a sequence of matrices
$L_{Z}$, $R_{Z}$ with different $Z$. It can be always segregated into clusters of
matrices
$$
L_{\vec Z}\equiv L_{Z_{i+s}}L_{Z_{i+s-1}}\cdots L_{Z_i}
$$
and
$$
R_{\vec Z}\equiv R_{Z_{j+k}}L_{Z_{j+k-1}}\cdots L_{Z_j}.
$$
The periodic extension of
expression (\ref{cf7}) is always an alternating sequence of matrices $L_{\vec Z}$ and
$R_{\vec Z}$:
$$
P_{Z_1,\dots Z_n}=\cdots (L_{\vec Z_s}R_{\vec Z_{s-1}})\cdots (L_{\vec Z_\cdot}R_{\vec Z_\cdot})\cdots\,.
$$
First note that it is impossible to have arbitrarily long sequences of only left or right matrices for a given
graph: the maximum length is restricted to be less of equal the maximum graph length of geodesics around holes.
This means that the length of a single cluster for a given graph is always bounded once the topology is fixed.

One can directly
calculate the product $(L_{\vec Z_1}R_{\vec Z_2})$ for $L_{\vec Z_1}=L_{Z_1}\cdots L_{Z_m}$ and
$R_{\vec Z_2}=R_{Z_{m+1}}\cdots R_{Z_{m+k}}$:
\bea
(L_{\vec Z_1}R_{\vec Z_2})&=&As^+_1s^+_2+B(s^-_1s^-_2+S_1s^-_2+S_1S_2)
\nonumber
\\
&&\quad+D(s^+_1s^-_2+s^+_1S_2)
+PS_1s^+_2,
\label{AA1}
\eea
where $S_j,s^\pm_j$, for $j=1,2$, are the following coefficient functions:
\bea
s^\pm_{1}&=&\e^{\pm\sum_{i=1}^mZ_i/2}, \qquad s^\pm_{2}=\e^{\pm\sum_{j=m+1}^{m+k}Z_i/2}
\nonumber
\\
S_1&=&\sum_{q=2}^{m}\e^{+\sum_{r=1}^{q-1}Z_r/2-\sum_{r=q}^{m}Z_r/2},\quad
S_2=\sum_{q=m+2}^{m+k}\e^{-\sum_{r=m+1}^{q-1}Z_r/2+\sum_{r=q}^{m+k}Z_r/2},
\nonumber
\eea
and $A,B,D,P$ are the special ($2\times2$)-matrices (``{\em letters}''):
\bea
&&A=\left(%
\begin{array}{rc}
  +1 & 0 \\
  -1 & 0 \\
\end{array}%
\right),\qquad
B=\left(%
\begin{array}{cc}
  0 & ~0 \\
  0 & +1 \\
\end{array}%
\right)
\nonumber
\\
&&D=\left(%
\begin{array}{cr}
  0 & -1 \\
  0 & +1 \\
\end{array}%
\right),\qquad
P=\left(%
\begin{array}{rc}
  0 & 0 \\
  -1 & 0 \\
\end{array}%
\right).
\label{alpha}
\eea
These letters possess interesting multiplication properties which are summarized in the
next lemma, whose proof is a routine calculation.

\begin{lemma}\label{alphabet}
{\em The alphabet lemma.}
The multiplication table of letters {\rm(\ref{alpha})} reads:
\be
\begin{tabular}{c|cccc}
   & A & B & D & P \\
  \hline
  A & A & $0$ & D & $0$ \\
  B & P & B & B & P \\
  D & A & D & D & A \\
  P & P & $0$ & B & $0$
\end{tabular},
\label{AA2}
\ee
so the trace of any product of these matrices is either unity or zero.
In the product of $t$ matrices
of form {\rm(\ref{AA1})}, the only monomials that survive are
\be
\bigl[(A+D)^{i_\alpha}D\bigr]B^{j_\beta}(B+P)\cdots
\bigl[(A+D)^{i_\rho}D\bigr]
B^{j_\omega}(B+P)\ \hbox{and}\ (A+D)^{t},\ B^t.
\ee
\end{lemma}

The main point is that almost {all} cancellations of letters in
long words are due to the local multiplication rules (\ref{AA2}).
This means that, having a long sequence of letters, say, $L_{I+N}$ from Lemma~\ref{qlem2}, we can
split it into pieces depending on sequences of letters $L_I$, $L_{I-1}$, $\tilde L_I$, and $\tilde L_{I-1}$,
where the index $I$ is also assumed to be big enough. That is,
let $L_{I+N}=L_IL_{I-1}{\tilde L}_IL_I\dots L_I$ comprise $p_N$ entries $L_I$ and $\tilde L_I$ and
$q_N$ entries $L_{I-1}$ and $\tilde L_{I-1}$.
We have then the following estimate:\footnote{This estimate also follows from the
properties of long geodesic lines in hyperbolic geometry: for two lines of large lengths~$L_1$ and $L_2$
intersecting at angle $\alpha$, the length~$L_3$ of the third side of the resulting triangle
is $L_1+L_2+\log((1-\cos\alpha)/2) +O(1/L)$. This also shows that our estimate is very rough.}
\be
\label{AA3}
|\log\tr L_{I+N}-p_N\log\tr L_I-q_N\log\tr L_{I-1}|< \hbox{C}~(q_N+p_N),
\ee
where the constant~C depends only on the Teichm\"uller space coordinates $Z_\alpha$ and on the
genus and the number of holes of the Riemann surface, and we have also used (\ref{cf8}). The ratio of the
coefficients is given by the continued fraction
\be
\label{AA4}
q_N/p_N=\frac{1}{a_{I+1}+{\displaystyle\frac{1}{
  \renewcommand{\arraystretch}{0}
  \arraycolsep=0.05em
  \begin{array}{ccc}
   \strut a_{I+2}+{} & & \\
     & \ddots & \\
     & &  {}+{\displaystyle\frac{1}{a_{I+N-1}+{\displaystyle\frac1{a_{I+N}}}}}
     \end{array}
     }}},
\ee
and also has a definite limit as $N\to\infty$. Now the estimate follows: up to
exponential corrections, ${\rm p.l.}(L)$ coincides with the $\log\tr L$, so for any $\eps>0$, let us choose
the index~$I$ such that $\eps\, {\rm p.l.}(L_{I-1})/2>C$ and $\eps\, {\rm p.l.} (L_{I})/2>C$. Thus,
\be
\label{AA5}
\frac{{\rm p.l.}(L)(L_{I+N})}{{\rm g.l.}(L_{I+N})}=(1+O(\eps/2))
\frac{p_N{\rm p.l.}(L)(L_I)+q_N{\rm p.l.}(L)(L_{I-1})}{p_N{\rm g.l.}(L_I)+q_N{\rm g.l.}(L_{I-1})},
\ee
and because the ratio $q_N/p_N$ has a definite limit as $N\to\infty$, there exists $N_0$ such that the
relative error of this ratio times the sum of ratios of proper and graph lengths
of $L_I$ and $L_{I-1}$ will not exceed $\eps/2$. Thus, the collective relative error for such
fixed $I$ and for all
$N>N_0$ is less than $\eps$, proving the theorem.

\appendix{Degeneracy of the Poisson structure}\label{nondegeneracy}

We shall explicitly calculate the degeneracy of the Poisson brackets (\ref{WP-PB}) for a
special graph and choose the graph whose ``building blocks'' are depicted in Figure~15. Namely, we have a
line tree subgraph comprising edges $X_i$ with attached subgraphs as in Figure~15a and 15b. Attaching
a subgraph of type a corresponds to adding a handle (increasing $g$ by unity)
while a subgraph of type b corresponds to adding a hole (increasing $s$ by unity).  We shall assume
that
$2g+2s>5$ to avoid the once-punctured torus, which is already handled separately in
Section 3.

\centerline{
\epsffile{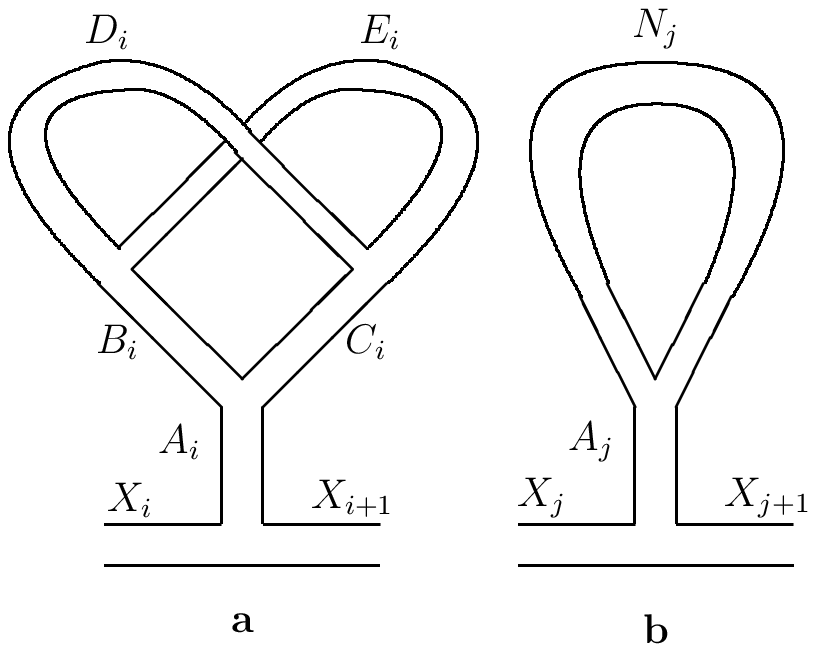}}

\centerline{\bf Figure 15-building blocks}

\vspace{10pt}

For the variables $A,B,C,D,E$ in Figure 15a, we have the Poisson bracket (sub)matrix
$$
\begin{tabular}{c|c|rrrr|}
    & $A_i$ & $B_i$ & $C_i$ & $D_i$ & $E_i$ \\
  \hline
  $A_i$ & 0 & 1 & $-1$ & 0 & 0 \\
  \hline
  $B_i$ & $-1$ & 0 & 1 & 1 & $-1$ \\
  $C_i$ & 1 & $-1$ & 0 & 1 & $-1$ \\
  $D_i$ & 0 & $-1$ & $-1$ & 0 & 2 \\
  $E_i$ & 0 & 1 & 1 & $-2$ & 0 \\
  \hline
\end{tabular}\ ,
$$
where the entries are the Poisson brackets between the corresponding variables.  Adding the last row
to the next-to-the-last row as well as adding the last column to the next-to-the-last column, then
adding the third row to the second row as well as the third column to the second column, we obtain the
matrix
$$
\left(%
\begin{array}{rrrrr}
  0 & 0 & -1 & 0 & 0 \\
  0 & 0 & 1 & 0 & -2 \\
  1 & -1 & 0 & 0 & -1 \\
  0 & 0 & 0 & 0 & 2 \\
  0 & 2 & 1 & -2 & 0 \\
\end{array}%
\right),
$$
which obviously has rank four and can be further
reduced (without adding the first column or row to any other)
to the form
$$
\left(%
\begin{array}{rrrrr}
  0 & 0 & 0 & 0 & 0 \\
  0 & 0 & +1 & 0 & 0 \\
  0 & -1 & 0 & 0 & 0 \\
  0 & 0 & 0 & 0 & +2 \\
  0 & 0 & 0 & -2 & 0 \\
\end{array}%
\right).
$$
Thus, erasing all columns and rows corresponding to the variables $B_i$, $C_i$, $D_i$,
and $E_i$ leaves invariant the rank of the Poisson bracket matrix.

Adjoining the subgraph in
Figure~15b creates exactly one degeneracy as the variable $N_j$ Poisson commutes with everything (as
it must be when adding a hole).

It remains only to calculate the rank of the matrix
corresponding to a tree graph with edges $X_i$ and $A_i$ remaining after erasing all $B$-, $C$-, $D$-,
$E$-, and
$N$-variable rows and columns. The corresponding Poisson bracket matrix has dimension $2g+2s-5>0$
and the simple block-diagonal form
$$
\begin{tabular}{r|rr|rr|rr|r}
  0 & 1 & $-1$ & 0 & 0 & $\cdot$ & $\cdot$ & $\cdot$\\
  \hline
  $-1$ & 0 & 1 & 0 & 0 &  &  & \\
  1 & $-1$ & 0 & 1 & $-1$ &  & & \\
  \hline
   & 0 & $-1$ & 0 & 1 & 0 & 0 &  \\
   & 0 & 1 & $-1$ & 0 & 1 & $-1$ & \\
  \hline
   &  &  & 0 & $-1$ & 0 & 1 & \\
   &  &  & 0 & 1 & $-1$ & 0 & $\ddots$\\
  \hline
   &  &  &  &  &  & $\ddots$ & $\ddots$\\
\end{tabular}\ .
$$
Adding each even-index row to its predecessor as well as adding each even-index column to
its predecessor, this reduces to the matrix whose only nonzero elements are
+1 on the main super-diagonal and -1
on the main sub-diagonal.  Since this matrix has full rank, the discussion is complete.

\end{document}